\documentclass[a4paper]{article}
\usepackage{amssymb,amscd}
\usepackage[all]{xy}

\def\goth{\mathfrak}
\def\double{\mathbb}

\def\cc{{\double C}}
\def\nn{{\double N}}
\def\zz{{\double Z}}

\def\rr{{\double R}}

\newtheorem{theorem}{Theorem}[section]

\newtheorem{corollary}[theorem]{Corollary}
\newtheorem{definition}[theorem]{Definition}
\newtheorem{proposition}[theorem]{Proposition}
\newtheorem{remark}[theorem]{Remark}

\newtheorem{example}[theorem]{Example}

\def\Uc{{\cal U}}

\def\limind{\mathop{\mathrm{lim}}\limits_{\longrightarrow}}

\def\si{\sigma}

\def\cinf{C^{\infty}}
\def\cinfc{C^{\infty}_c}

\newcommand{\be}{\begin{equation}}
\newcommand{\ee}{\end{equation}}
\newcommand{\beq}{\begin{eqnarray}}
\newcommand{\eeq}{\end{eqnarray}}
\newcommand{\om}{\omega}
\newcommand{\Om}{\Omega}
\newcommand{\al}{\alpha}
\def\nat{\natural}

\newcommand{\La}{\Lambda}
\newcommand{\la}{\lambda}
\newcommand{\Ec}{{\cal E}}
\newcommand{\Vc}{{\cal V}}
\newcommand{\Lc}{{\cal L}}
\newcommand{\non}{\nonumber}
\newcommand{\eps}{\varepsilon}

\newcommand{\Wc}{{\cal W}}

\newcommand{\Jc}{{\cal J}}
\newcommand{\Ind}{{\mathop{\mathrm{Ind}}}}
\def\ch{\mathrm{ch}}

\newcommand{\Tr}{{\mathop{\mathrm{Tr}}}}
\newcommand{\tr}{{\mathop{\mathrm{tr}}}}
\newcommand{\Ac}{{\cal A}}

\newcommand{\te}{\theta}

\newcommand{\cqfd}{\hfill\rule{1ex}{1ex}}

\def\dd{\,\mathrm{\bf d}}

\def\Hc{{\cal H}}

\def\Bc{{\cal B}}
\def\Cc{{\cal C}}
\def\Jc{{\cal J}}
\def\Kc{{\cal K}}

\def\Pc{{\cal P}}

\def\coker{\mathop{\mathrm{Coker}}}

\def\bb{\overline{b}}

\def\hom{{\mathop{\mathrm{Hom}}}}

\def\End{{\mathop{\mathrm{End}}}}

\def\hotimes{\hat{\otimes}}

\def\St{\widetilde{S}}
\def\Ut{\widetilde{U}}

\def\Sg{{\goth S}}
\def\Act{\widetilde{\cal A}}
\def\Tct{\widetilde{\cal T}}
\def\Bct{\widetilde{\cal B}}

\def\Ect{\widetilde{\cal E}}
\def\at{\widetilde{a}}
\def\gt{\widetilde{g}}

\def\Sge{{\goth S_{\epsilon}}}
\def\Sgan{{\goth S_{\mathrm{an}}}}

\def\chih{\widehat{\chi}}

\def\Ome{\Omega_{\epsilon}}

\def\Oman{\Omega_{\mathrm{an}}}

\def\Tc{{\cal T}}

\def\Uc{{\cal U}}
\def\Kc{{\cal K}}

\def\supp{\mathrm{supp}\,}
\def\diff{\mathrm{Diff}\,}
\def\eh{\hat{e}}

\begin{document}

\begin{center}

{\large\bf THE EQUIVARIANT INDEX THEOREM IN\\[1MM]
ENTIRE CYCLIC COHOMOLOGY\footnotemark[1]}
\vskip 1cm
{\bf Denis PERROT}
\vskip 0.5cm
Institut Camille Jordan, Universit\'e Claude Bernard Lyon 1,\\
21 av. Claude Bernard, 69622 Villeurbanne cedex, France \\[2mm]
{\tt perrot@math.univ-lyon1.fr}\\[2mm]
\today
\end{center}
\vskip 0.5cm
\begin{abstract}
Let $G$ be a locally compact group acting smoothly and properly by isometries on a complete Riemannian manifold $M$, with compact quotient $G\backslash M$. There is an assembly map $\mu: K_*^G(M)\to K_*(\Bc)$ which associates to any $G$-equivariant $K$-homology class on $M$, an element of the topological $K$-theory of a suitable Banach completion $\Bc$ of the convolution algebra of continuous compactly supported functions on $G$. The aim of this paper is to calculate the composition of the assembly map with the Chern character in entire cyclic homology $K_*(\Bc)\to HE_*(\Bc)$. We prove an index theorem reducing this computation to a cup-product in bivariant entire cyclic cohomology. As a consequence we obtain an explicit localization formula which includes, as particular cases, the equivariant Atiyah-Segal-Singer index theorem when $G$ is compact, and the Connes-Moscovici index theorem for G-coverings when $G$ is discrete. The proof is based on the bivariant Chern character introduced in previous papers.
\end{abstract}

\vskip 0.5cm

\footnotetext[1]{This work was partially made during the postdoctoral position of the author at the Mathematisches Institut, Universit\"at M\"unster, Germany.}


\section{Introduction}

Let $M$ be a smooth complete Riemannian manifold without boundary, on which a separable locally compact group $G$ acts smoothly and properly by isometries, with compact quotient $G\backslash M$. There is an analogue of the Baum-Connes assembly map \cite{BC, BCH}
\be
\mu: K_*^G(M)\to K_*(\Bc)\ ,
\ee
from the equivariant $K$-homology of $M$ to the (topological) $K$-theory of an admissible Banach completion $\Bc$ of the convolution algebra of continuous, compactly supported functions on $G$. The definition of an admissible completion is given in \ref{dadm}: roughly speaking, the elements of $\Bc$ are locally integrable functions on $G$ with certain decay properties at infinity. The construction of the assembly map is inspired from \cite{C1} and \cite{Kas}. The aim of this paper is to give an explicit formula calculating the composition of $\mu$ with the Chern character $\ch:K_*(\Bc)\to HE_*(\Bc)$ in entire cyclic homology \cite{C2,Me}. For this purpose we will use the bivariant Chern character and related techniques developped in \cite{P1,P2}. We insist on the fact that the same methods should in principle apply to more general situations, including for example locally compact groupoids, although we will only deal with groups and manifolds in the present work.\\

In many situations the equivariant $K$-homology classes $[D]\in K_*^G(M)$ of interest  are represented by $G$-invariant elliptic \emph{differential} operators of order one, hereafter denoted by $D$. Our main result is that the  Chern character of $\mu(D)$ in the entire cyclic homology $HE_*(\Bc)$ may be written as a certain composition product in bivariant cyclic cohomology. We first introduce the crossed-product algebra $\Ac=\cinfc(M)\rtimes G$. It is provided with a canonical $K$-theory element $[e]\in K_0(\Ac)$, which may be constructed via a cut-off function over $M$ as in \cite{Kas}. From the operator $D$, one can build an unbounded $\Ac$-$\Bc$-bimodule $(\Ec,\rho,D)$: it is a suitable smooth version of Kasparov bimodule, according to the construction of \cite{P1}. Its Chern character lies in the bivariant entire cyclic cohomology $HE_*(\Ac,\Bc)$. The statement of the equivariant index theorem \ref{tind} is the following:
\begin{theorem}
Let $D$ be a $G$-invariant elliptic differential operator of order one representing an equivariant $K$-homology class $[D]\in K_*^G(M)$. Consider the crossed-product algebra $\Ac=\cinfc(M)\rtimes G$ and its canonical $K$-theory class $[e]\in K_0(\Ac)$, and let $\Bc$ be any admissible completion of the convolution algebra $C_c(G)$. Then the Chern character of the image of $[D]$ under the analytic assembly map $\mu:K_*^G(M) \to K_*(\Bc)$ is given by the cup-product in bivariant entire cyclic cohomology
\be
\ch\circ\mu(D)= \ch(\Ec,\rho,D)\cdot \ch(e)\ \in HE_*(\Bc)\ ,
\ee
where $\ch(\Ec,\rho,D)\in HE_*(\Ac,\Bc)$ is the bivariant Chern character of the unbounded bimodule associated to $D$, and $\ch(e)\in HE_0(\Ac)$ is the Chern character of $[e]$.
\end{theorem}
This theorem should be viewed as the cohomological version of the construction of the assembly map using bivariant $K$-theory \cite{BCH,C1}. The Chern character $\ch(\Ec,\rho,D)$ is defined via an explicit formula of JLO-type \cite{JLO}, involving the heat operator $\exp(-t^2D^2)$, $t>0$. \\
Concretely the cyclic homology class $\ch\circ\mu(D)\in HE_*(\Bc)$ is represented, for any choice of parameter $t>0$, by an entire cycle
\be
\ch(tD)\in \Oman\Bc
\ee
in the complex $\Oman\Bc$ of entire chains over $\Bc$ (see section \ref{sent} and \cite{Me}). The cycle $\ch(tD)$ is given by the collection of its components $\ch_n(tD)$ in the spaces of non-commutative $n$-forms $\Om^n\Bc$, in any degree $n\in\nn$. Since $\Bc$ is an algebra of functions over the group $G$, the component $\ch_n(tD)$ is actually a function over the locally compact space $G^n\cup G^{n+1}$. In favourable circumstances, for example when $D$ is a Dirac-type operator acting on the sections of a $G$-equivariant Clifford module $E$ over $M$, the limit $t\to 0$ of the function $\ch_n(tD)$ evaluated at a point $\gt=(g_1,\ldots, g_n)\in G^{n}$ or $\gt=(g_0,\ldots, g_n)\in G^{n+1}$ may be explicitly computed. This leads to the following localization formula (Corollary \ref{cloc}), involving the submanifolds $M_g\subset M$ of fixed points for each element $g\in G$. 
\begin{corollary}
The $n$th degree component $\ch_n(tD)\in \Om^n\Bc$, viewed as a function on $G^n\cup G^{n+1}$, admits a pointwise limit when $t\to 0$ given by the localization formula
\be
\lim_{t\to 0}\ch_n(tD)(\gt)= \sum_{M_g}\frac{(-)^{q/2}}{(2\pi i)^{d/2}}\int_{M_g}\widehat{A}(M_g)\frac{\ch(E/S,g)}{\ch(S_N,g)}\, \ch_n(e)(\gt)\ ,
\ee
where $\gt\in G^n\cup G^{n+1}$, $g=g_n\ldots g_1 \in G$ (resp. $g=g_n\ldots g_0$) if $\gt=(g_1,\ldots, g_n)$ (resp. $\gt=(g_0,\ldots, g_n)$), and the sum runs over the fixed manifolds $M_g$ of all possible dimensions $d$ and codimensions $q=\mbox{dim}\, M-d$. 
\end{corollary}
One recognizes the usual ingredients of the Atiyah-Segal-Singer index theorem: the $\widehat{A}$-genus of the fixed submanifolds as well as some equivariant characteristic classes of vector bundles $\ch(E/S,g)$ and $\ch(S_N,g)$. See section \ref{sind} for details. The last ingredient is a noncommutative Chern character $\ch_n(e)$ associated to the canonical $K$-theory class $[e]\in K_0(\Ac)$. It is a function on $G^n\cup G^{n+1}$ with values in the space of differential forms with compact support on $M$, see Definition \ref{dchern}.\\

The above localization formula is a generalization of several known results. For example, when $G$ is a compact group one recovers the Atiyah-Singer equivariant index theorem \cite{AS}. The noncommutative Chern character $\ch_n(e)$ is trivial in this case. See also \cite{L1,L2} for a similar approach based on the JLO cocycle in the case of finite groups.  On the other hand, if $G$ is a discrete countable group acting freely and properly on $M$, one recovers the Connes-Moscovici index theorem for $G$-coverings \cite{CM90}. In the latter case, $\mu(D)$ is paired with periodic cyclic cocycles arising from group cohomology. Our result is a generalization in two directions. First, it is valid for general locally compact groups and proper actions, and second, it holds in the entire (as opposed to periodic) cyclic homology of certain Banach completions $\Bc$ of the convolution algebra of $G$. We believe that these points may be useful for the study of assembly maps and related conjectures \cite{BCH}.\\

The paper is organized as follows. In section \ref{sent} we recall the basic definitions of entire cyclic cohomology. The analytic assembly map is considered in section \ref{sass}. Section \ref{sbiv} contains the construction of the bivariant Chern character introduced in \cite{P1}. The equivariant index theorem is proved in section \ref{sind} together with the localization formula.\\

\section{Entire cyclic cohomology}\label{sent}

Entire cyclic cohomology was first designed by Connes for Banach algebras \cite{C2}, or locally convex topological algebras. More generally, the theory was adapted by Meyer \cite{Me} to the category of \emph{bornological} algebras. This is the framework we used in \cite{P1,P2} for the construction of a bivariant Chern character. We recall briefly in this section some basic definitions and facts.\\

A bornological vector space $\Vc$ is a vector space (over $\cc$) endowed with a collection $\Sg(\Vc)$ of subsets $S\subset \Vc$, called the bornology of $\Vc$, satisfying certain axioms reminiscent of the properties fulfilled by the bounded sets of a normed space \cite{HN, Me}:
\begin{itemize}
\item $\{x\}\in\Sg(\Vc)$ for any vector $x\in\Vc$.
\item $S_1+S_2\in\Sg(\Vc)$ for any $S_1,S_2\in\Sg(\Vc)$.
\item If $S\in \Sg(\Vc)$, then $T\in\Sg(\Vc)$ for any $T\subset S$.
\item $S^{\Diamond}\in\Sg(\Vc)$ for any $S\in\Sg(\Vc)$,
\end{itemize}
where $S^{\Diamond}$ is the circled convex hull of the subset $S$. The elements of $\Sg(\Vc)$ are called the {\it small} subsets of $\Vc$. If a small subset $S$ is a disk (i.e. circled and convex), then its linear span $\Vc_S$ carries a unique seminorm for which the closure of $S$ is the unit ball; $S$ is called completant iff $\Vc_S$ is a Banach space. $\Vc$ is a {\it complete bornological vector space} iff any small subset $T$ is contained in some completant small disk $S\in\Sg(\Vc)$. Typical examples of complete bornological spaces are provided by Banach or complete locally convex spaces, endowed with the bornology corresponding to the collection of all bounded subsets.\\
The interesting linear maps between two bornological spaces $\Vc$ and $\Wc$ are those which send $\Sg(\Vc)$ to $\Sg(\Wc)$. Such maps are called {\it bounded}. The vector space of bounded linear maps $\hom(\Vc,\Wc)$ is itself a bornological space, the small subsets corresponding to equibounded maps (see \cite{HN,Me}). It is complete if $\Wc$ is complete. Similarly, an $n$-linear map $\Vc_1\times\ldots\times \Vc_n\to\Wc$ is bounded if it sends an $n$-tuple of small sets $(S_1,\ldots, S_n)$ to a small set of $\Wc$. A {\it bornological algebra} $\Ac$ is a bornological vector space together with a bounded bilinear map $\Ac\times\Ac\to\Ac$. We will be concerned only with associative bornological algebras. \\
Given a bornological space $\Vc$, its {\it completion} $\Vc^c$ is a complete bornological space defined as the solution of a universal problem concerning the factorization of bounded maps $\Vc\to\Wc$ with complete target $\Wc$. This completion always exists \cite{HN,Me}. If $\Vc$ is a normed space endowed with the bornology of bounded subsets, then its bornological completion coincides with its Hausdorff completion. The {\it completed tensor product} $\Vc_1\hotimes\Vc_2$ of two bornological spaces $\Vc_1$ and $\Vc_2$ is the completion of their algebraic tensor product for the bornology generated by the bismall sets $S_1\otimes S_2$, $\forall S_i\in \Sg(\Vc_i)$. The completed tensor product is associative, whence the definition of the $n$-fold completed tensor product $\Vc_1\hotimes\ldots\hotimes\Vc_n$. If $\Ac$ and $\Bc$ are bornological algebras, their completed tensor product $\Ac\hotimes\Bc$ is a complete bornological algebra. $\hotimes$ is additive even for infinite direct sums, and commutes with inductive limits. The main examples used in this paper are Fr\'echet algebras and LF-algebras:
\begin{example}
\textup{A Fr\'echet space endowed with the bornology of all bounded subsets is a complete bornological space. A linear map $\Vc\to\Wc$ between two Fr\'echet spaces is bounded iff it is continuous for the Fr\'echet topology. In particular, a Fr\'echet algebra $\Ac$ is automatically a bornological algebra, i.e. the (jointly) continuous multiplication $\Ac\times\Ac\to\Ac$ is a bounded bilinear map. The situation is more delicate for completed bornological tensor products. As a vector space, the bornological tensor product $\Vc_1\hotimes\Vc_2$ of two Fr\'echet spaces coincides with the completed projective tensor product $\Vc_1\hotimes_{\pi}\Vc_2$. However, the induced bornology on $\Vc_1\hotimes\Vc_2$ may differ from the bornology of bounded subsets of the Fr\'echet space $\Vc_1\hotimes_{\pi}\Vc_2$. These bornologies nevertheless coincide in some particular cases (see \cite{Me}):\\
i) If $\Vc_1$ is nuclear \cite{T} and $\Vc_2$ is arbitrary;\\
ii) If $\Vc_1$ is the space of integrable functions over a locally compact space w.r.t. a Borel measure and $\Vc_2$ is arbitrary;\\
iii) If both $\Vc_1$ and $\Vc_2$ are Banach spaces.}
\end{example}
\begin{example}
\textup{An LF-space is a locally convex topological vector space obtained as an increasing countable union $\Vc=\cup_{i\in\nn}\Vc_i$, where $\Vc_i$ are Fr\'echet spaces with respect to the subspace topology. $\Vc$ is endowed with the finest topology making the inclusions $\Vc_i\to\Vc$ continuous. We gift $\Vc$ with the bornology of bounded subsets: a subset is small iff there exists a Fr\'echet subspace $\Vc_i$ containing $S$ as a bounded set. Then $\Vc$ is a complete bornological vector space. An LF-algebra $\Ac$ is an LF-space with a separately continuous multiplication $\Ac\times\Ac\to\Ac$. Then $\Ac$ is automatically a bornological algebra \cite{Me}. An example of LF-algebra is the convolution algebra $C_c(G)$ of continuous, compactly supported $\cc$-valued functions on a separable locally compact group $G$. $C_c(G)$ is actually the union, over compact subsets $K\subset G$, of the Banach spaces $C_K(G)=\{f\in C_c(G)\ |\ \supp f\subset K\}$ endowed with the supremum norm. It is easy to see that the multiplication is separately continuous but not jointly continuous.}
\end{example}

Let us now recall Meyer's formulation of entire cyclic homology, cohomology and bivariant cohomology for bornological algebras. There are in fact two equivalent ways to describe the entire cyclic cohomology of a complete bornological algebra $\Ac$. The first one is to use Connes' $(b,B)$-complex of non-commutative forms completed with respect to a certain bornology; we call this completion $\Ome \Ac$. The second one is to introduce the $X$-complex of Cuntz and Quillen \cite{CQ1,CQ2} for the analytic tensor algebra $\Tc\Ac$. These complexes are homotopy equivalent \cite{Me}, and give rise to the definition of entire cyclic cohomology. The construction of the bivariant Chern character proposed in \cite{P1} uses simultaneously the $(b,B)$-complex and $X$-complex approaches.\\
Let $\Ac$ be a complete bornological algebra. The algebra of non-commutative differential forms \cite{C1} over $\Ac$ is the direct sum $\Om\Ac=\bigoplus_{n\ge 0}\Om^n\Ac$ of the $n$-forms subspaces $\Om^n\Ac=\Act\hotimes\Ac^{\hotimes n}$ for $n\ge 1$ and $\Om^0\Ac=\Ac$, where $\Act=\Ac\oplus \cc$ is the unitalization of $\Ac$. It is customary to use the differential notation $a_0da_1\ldots da_n$ (resp. $da_1\ldots da_n$) for the string $a_0\otimes a_1\ldots\otimes a_n$ (resp. $1\otimes a_1\ldots\otimes a_n)$. The differential $d: \Om^n\Ac\to\Om^{n+1}\Ac$ is uniquely specified by $d(a_0da_1\ldots da_n)=da_0da_1\ldots da_n$ and $d^2=0$. The multiplication in $\Om\Ac$ is defined as usual and fulfills the Leibniz rule $d(\om_1\om_2)=d\om_1\om_2 +(-)^{|\om_1|}\om_1d\om_2$, where $|\om_1|$ is the degree of the (homogeneous) element $\om_1$. Each $\Om^n\Ac$ is a complete bornological space by construction, and we endow $\Om\Ac$ with the direct sum bornology. This turns $\Om\Ac$ into a complete bornological differential graded (DG) algebra, i.e. the multiplication map and $d$ are bounded. $\Om\Ac$ is the universal complete bornological DG algebra generated by $\Ac$. It is $\zz_2$-graded by the even/odd degree of differential forms.\\
On $\Om\Ac$ are defined various operators. First of all, the Hochschild boundary $b: \Om^{n+1}\Ac\to \Om^n\Ac$ is given by $b(\om da)=(-)^n[\om,a]$ for $\om\in\Om^n\Ac$, and $b=0$ on $\Om^0\Ac=\Ac$. One easily shows that $b$ is bounded and $b^2=0$. Then the Karoubi operator $\kappa:\Om^n\Ac\to \Om^n\Ac$ is constructed out of $b$ and $d$:
\be
1-\kappa=db+bd\ .
\ee
Therefore $\kappa$ is bounded and commutes with $b$ and $d$. The last operator is Connes' $B:\Om^n\Ac\to \Om^{n+1}\Ac$,
\be
B=(1+\kappa+\ldots+\kappa^n)d\quad \mbox{on} \ \Om^n\Ac\ ,
\ee
which is bounded and verifies $B^2=0=Bb+bB$ and $B\kappa=\kappa B=B$. Hence $\Om\Ac$ endowed with the total coboundary $b+B$ is a $\zz_2$-graded bornological complex. We denote by $\hom^f(\Om\Ac,\cc)$ the complex of cochains with finite dimension over $\Ac$, i.e. the complex of bounded linear maps $\Om\Ac\to\cc$ which vanish on $\Om^n\Ac$ for $n$ large. One has
\be
\hom^f(\Om\Ac,\cc)=\bigoplus_{n\geq 0}\hom(\Om^n\Ac,\cc)\ .
\ee
\begin{definition}
Let $\Ac$ be a complete bornological algebra. The $\zz_2$-graded complex of bounded cochains $\hom^f(\Om\Ac,\cc)$ \emph{with finite dimension} calculates the periodic cyclic cohomology of $\Ac$:
\be
HP^i(\Ac)=H^i(\hom^f(\Om\Ac,\cc),b+B)\ ,\qquad i=0,1\ .
\ee
\end{definition}
The periodic theory is not always sufficient for analytic applications, because it contains only finite-dimensional cocycles. The entire cyclic cohomology of Connes \cite{C2} provides a more general cyclic theory dealing with infinite-dimensional cocycles. Following Meyer \cite{Me}, we define two other bornologies on $\Om\Ac$:
\vskip 2mm
\noindent $\bullet$ The {\bf entire bornology} $\Sge(\Om\Ac)$ is generated by the sets
\be
\bigcup_{n\ge 0} [n/2]!\,\St (dS)^n\quad ,\ S\in\Sg(\Ac)\ ,\label{sge}
\ee 
where $[n/2]=k$ if $n=2k$ or $n=2k+1$, and $\St=S\cup\{1\}$. That is, a subset of $\Om\Ac$ is small iff it is contained in the circled convex hull of a set (\ref{sge}). We write $\Ome\Ac$ for the completion of $\Om\Ac$ with respect to this bornology. $\Ome\Ac$ will be the $(b,B)$-complex of entire chains.
\vskip 2mm
\noindent $\bullet$  The {\bf analytic bornology} $\Sgan(\Om\Ac)$ is generated by the sets $\bigcup_{n\ge 0} \St(dS)^n$, $S\in\Sg(\Ac)$. The corresponding completion of $\Om \Ac$ is denoted $\Oman\Ac$. It is related to the $X$-complex description of entire cyclic homology (see below).
\vskip 2mm
The multiplication in $\Om\Ac$ is bounded for the two bornologies above, as well as all the operators $d,b,\kappa,B$. Moreover, the $\zz_2$-graduation of $\Om\Ac$ given by even and odd forms is preserved by the completion process, so that $\Ome\Ac$ and $\Oman\Ac$ are $\zz_2$-graded differential algebras, endowed with the operators $b,\kappa,B$ fulfilling the usual relations. In particular, $\Ome\Ac$ is called the $(b,B)$-complex of {\it entire chains}. Note also that the linear map $c:\Oman\Ac\to\Ome\Ac$ 
\be
c(a_0da_1\ldots da_n)=(-)^{[n/2]}[n/2]!\, a_0da_1\ldots da_n\quad \forall n\in\nn\ .\label{res}
\ee
obviously provides a bornological isomorphism between the spaces $\Ome\Ac$ and $\Oman\Ac$.\\

Another algebra related to differential forms is the analytic tensor algebra. Let $\Om\Ac=\Om^+\Ac\oplus\Om^-\Ac$ be the $\zz_2$-graded algebra of differential forms. The even part $\Om^+\Ac$ is a trivialy graded subalgebra. Following Cuntz and Quillen \cite{CQ1}, we deform the usual product and endow $\Om^+\Ac$ with the {\it Fedosov product} 
\be
\om_1\odot\om_2 =\om_1\om_2 -d\om_1d\om_2\ ,\quad \om_{1,2}\in\Om^+\Ac\ .
\ee
Associativity is easy to check. In fact the algebra $(\Om^+\Ac,\odot)$ is isomorphic to the non-unital tensor algebra $T\Ac=\bigoplus_{n\ge 1}\Ac^{\hotimes n}$, under the correspondence
\be
\Om^+\Ac\ni a_0da_1\ldots da_{2n} \longleftrightarrow a_0\otimes\om(a_1,a_2)\otimes\ldots \otimes\om(a_{2n-1},a_{2n}) \in T\Ac\ ,\label{cor}
\ee
where $\om(a_i,a_j):=a_ia_j- a_i\otimes a_j \in \Ac\oplus \Ac^{\hotimes 2}$ is the {\it curvature} of the pair $(a_i,a_j)$. It turns out that the Fedosov product $\odot$ is bounded for the bornology $\Sgan$ restricted to $\Om^+\Ac$ \cite{Me}, and thus extends to the analytic completion $\Oman^+\Ac$. The complete bornological algebra $(\Oman^+\Ac,\odot)$ is also denoted by $\Tc\Ac$ and called the {\it analytic tensor algebra} of $\Ac$ in \cite{Me}.\\

We now turn to the description of the $X$-complex. It first appeared in the coalgebra context in Quillen's work \cite{Q2}, and subsequently was used by Cuntz-Quillen in their formulation of cyclic homology \cite{CQ1,CQ2}. Here we recall the $X$-complex construction for bornological algebras, following \cite{Me}. \\
Let $\Ac$ be a complete bornological algebra. The $X$-complex of $\Ac$ is the $\zz_2$-graded complex
\be
X(\Ac):\quad \Ac\ \xymatrix@1{\ar@<0.5ex>[r]^{\nat \dd} &  \ar@<0.5ex>[l]^{\bb}}\ \Om^1\Ac_{\nat}\ ,
\ee
where $\Om^1\Ac_{\nat}$ is the completion of the commutator quotient space $\Om^1\Ac/b\Om^2\Ac=\Om^1\Ac/[\Ac,\Om^1\Ac]$ endowed with the quotient bornology. The class of the generic element $(a_0\dd a_1\, \mbox{mod}\, [,])\in \Om^1\Ac_{\nat}$ is usually denoted by $\nat a_0\dd a_1$. We use the symbol $\dd$ for the universal derivation $\dd:\Ac\to\Om^1\Ac$ so that no confusion arises when dealing with the $X$-complex of the analytic tensor algebra (see below). The map $\nat \dd:\Ac\to \Om^1\Ac_{\nat}$ thus sends $a\in\Ac$ to $\nat \dd a$. Also, the Hochschild boundary $b:\Om^1\Ac\to \Ac$ vanishes on the commutators $[\Ac,\Om^1\Ac]$, hence passes to a well-defined map $\bb:\Om^1\Ac_{\nat}\to\Ac$. Explicitly the image of $\nat a_0\dd a_1$ by $\bb$ is the commutator $[a_0,a_1]$. These maps are bounded and satisfy $\nat \dd\circ \bb=0$ and $\bb\circ\nat \dd=0$, so that $X(\Ac)$ indeed defines a complete $\zz_2$-graded bornological  complex.\\
We now focus on the $X$-complex of the analytic tensor algebra $\Tc\Ac$. In that case, the quotient space $\Om^1\Tc\Ac_{\nat}=\Om^1\Tc\Ac/[\Tc\Ac,\Om^1\Tc\Ac]$ is always complete, and as a bornological vector space $X(\Tc\Ac)$ is canonically isomorphic to the analytic completion $\Oman\Ac$. This correspondence goes as follows \cite{CQ1,Me}. First, one knows that the even part $X_0(\Tc\Ac)=\Tc\Ac$ is isomorphic to the Fedosov algebra of differential forms $\Oman^+\Ac$. Second, one has a $\Tc\Ac$-bimodule isomorphism
\beq
\Om^1\Tc\Ac &\simeq& \Tct\Ac\hotimes \Ac\hotimes\Tct\Ac\\
x\, \dd a\, y &\leftrightarrow& x\otimes a\otimes y\qquad \mbox{for}\ a\in\Ac\ ,\ x,y\in\Tct\Ac\ , \non
\eeq
where $\Tct\Ac:=\cc\oplus\Tc\Ac$ is the unitalization of $\Tc\Ac$. This implies that the bornological space $\Om^1\Tc\Ac_{\nat}$ is isomorphic to $\Tct\Ac\hotimes\Ac$, which can further be identified with the analytic completion of odd forms $\Oman^-\Ac$, through the correspondence $x\otimes a\leftrightarrow xda$, $\forall a\in \Ac, x\in \Tct\Ac$. Thus collecting the even part $X_0(\Tc\Ac)=\Tc\Ac$ and the odd part $X_1(\Tc\Ac)=\Om^1\Tc\Ac_{\nat}$ together, yields a linear bornological isomorphism $X(\Tc\Ac)\simeq \Oman\Ac$. We still denote by $(\nat \dd,\bb)$ the boundaries induced on $\Oman\Ac$ through this isomorphism; Cuntz and Quillen explicitly computed them in terms of the usual operators on differential forms \cite{CQ1}:
\beq
\bb&=&b-(1+\kappa)d\quad \mbox{on}\ \Om^{2n+1}\Ac\ ,\\
\nat \dd&=&\sum_{i=0}^{2n}\kappa^id -\sum_{i=0}^{n-1}\kappa^{2i}b\quad \mbox{on}\ \Om^{2n}\Ac\ .\non
\eeq

The crucial result is that the bornological isomorphism (\ref{res}) induces a homotopy equivalence between the complex $(\Oman\Ac,\nat \dd,\bb)=X(\Tc\Ac)$ and the complex of entire chains $\Ome\Ac$ endowed with the differential $(b+B)$. See  \cite{CQ1,Me} for details. This leads to the definition of entire cyclic (co)homology:
\begin{definition}
Let $\Ac$ be a complete bornological algebra. \\
i) The entire cyclic homology of $\Ac$ is the homology of the $X$-complex of the analytic tensor algebra $\Tc\Ac$:
\be
HE_*(\Ac)=H_*(X(\Tc\Ac))\ ,
\ee
or equivalently, the $(b+B)$-homology of the $\zz_2$-graded complex of entire chains $\Ome\Ac$.\\
ii) Let $\hom(X(\Tc\Ac),\cc)$ be the $\zz_2$-graded complex of \emph{bounded} linear maps from $X(\Tc\Ac)$ to $\cc$, with differential the transposed of $(\nat \dd,\bb)$. Then the entire cyclic cohomology of $\Ac$ is the cohomology of this dual complex:
\be
HE^*(\Ac)=H^*(\hom(X(\Tc\Ac),\cc))\ .
\ee
iii) If $\Ac$ and $\Bc$ are complete bornological algebras, then $\hom(X(\Tc\Ac),X(\Tc\Bc))$ denotes the space of \emph{bounded} linear maps from $X(\Tc\Ac)$ to $X(\Tc\Bc)$. It is naturally a complete $\zz_2$-graded bornological complex, the differential of a map $f$ of degree $|f|$ is given by the graded commutator $(\nat \dd,\bb)\circ f-(-)^{|f|}f\circ(\nat \dd,\bb)$. The bivariant entire cyclic cohomology of $\Ac$ and $\Bc$ is the cohomology of this complex:
\be
HE_*(\Ac,\Bc)=H_*(\hom(X(\Tc\Ac),X(\Tc\Bc)))\ .
\ee
\end{definition}
In the case $\Ac=\cc$, one shows \cite{Me} that $X(\Tc\cc)$ is homotopically equivalent to $X(\cc): \cc\rightleftarrows 0$, thus the entire cyclic homology of $\cc$ is simply $HE_0(\cc)=\cc$ and $HE_1(\cc)=0$. This implies that for any complete bornological algebra $\Ac$, we get the usual isomorphisms $HE_*(\cc,\Ac)\simeq HE_*(\Ac)$ and $HE_*(\Ac,\cc)\simeq HE^*(\Ac)$. Furthermore, since the composition of bounded maps is bounded, there is a well-defined composition product on bivariant entire cyclic cohomology:
\be
HE_i(\Ac,\Bc)\times HE_j(\Bc,\Cc)\to HE_{i+j}(\Ac,\Cc)\ ,\quad i,j\in \zz_2
\ee
for complete bornological algebras $\Ac,\Bc,\Cc$. As an example of bivariant class, any bounded homomorphism $\rho: \Ac\to \Bc$ extends to a bounded homomorphism $\rho_*:\Tc\Ac\to \Tc\Bc$ by setting $\rho_*(a_1\otimes\ldots\otimes a_n)=\rho(a_1)\otimes\ldots\otimes\rho(a_n)$. The boundedness of $\rho_*$ becomes obvious once we rewrite it using the isomorphism $\Tc\Ac\simeq(\Oman^+\Ac,\odot)$, since $\rho_*(a_0da_1\ldots da_{2n})= \rho(a_1)d\rho(a_1)\ldots d\rho(a_{2n})$. The homomorphism $\rho_*$ gives rise to a bounded $X$-complex morphism $X(\rho_*): X(\Tc\Ac)\to X(\Tc\Bc)$:
\beq
x&\mapsto& \rho_*(x)\\
\nat x\dd y &\mapsto& \nat \rho_*(x)\dd\rho_*(y)\qquad \forall x,y\in\Tc\Ac\ ,\non
\eeq
and $X(\rho_*)$ defines a class $\ch(\rho)$ in $HE_0(\Ac,\Bc)$. It is the simplest example of bivariant Chern character induced by an homomorphism. In the particular case $\Ac=\cc$, an homomorphism $\rho:\cc\to \Bc$ is unambiguously specified by the idempotent $\rho(1)=e\in\Bc$. The homotopy equivalence $X(\Tc\cc)\sim \cc$ implies that $\ch(e)$ is an entire cyclic homology class in $HE_0(\Bc)$. It is represented by the following idempotent $\eh\in\Tc\Bc$, obtained by lifting (see \cite{CQ1,Me}):
\be
\eh=e+ \sum_{k\ge 1} \frac{(2k)!}{(k!)^2} (e-\frac{1}{2})\otimes (e-e\otimes e)^{\otimes k}\ \in \Tc\Bc\ ,
\ee
or equivalently in terms of differential forms,
\be
\eh=e+ \sum_{k\ge 1} \frac{(2k)!}{(k!)^2} (e-\frac{1}{2})(dede)^{k}\ \in \Oman^+\Bc\ .
\ee
At last, remark that unlike the periodic theory, the entire cyclic cohomology $HE^*(\Ac)$ contains infinite-dimensional cocycles. Indeed from the isomorphism $X(\Tc\Ac)\cong\Oman\Ac$, we see that an entire cochain $\varphi\in \hom(X(\Tc\Ac),\cc)$ is a collection of linear maps $\Om^n\Ac\to\cc$ for all $n$, which satisfy the following growth condition. For any small subset $S$ in the bornology of $\Ac$, and any elements $\at_0\in S\cup\{1\}$, $a_1,\ldots ,a_n\in S$, one has
\be
|\varphi(\at_0da_1\ldots da_n)|\leq C_S\ ,
\ee
where $C_S$ is a constant depending on $S$ but not on $n$. Hence the injection of the complex of finite-dimensional cochains into the entire complex induces a natural map $HP^*(\Ac)\to HE^*(\Ac)$.

\section{Analytic assembly map}\label{sass}

Let $G$ be a separable locally compact group. Denote by $C_c(G)$ be the convolution algebra of continuous, compactly supported $\cc$-valued functions on $G$. The product is given in terms of a right-invariant Haar measure $dh$:
\be
(b_1b_2)(g)=\int_G dh\, b_1(h)b_2(gh^{-1})\ ,\quad \forall\ b_i\in C_c(G)\ ,\ g\in G\ . \label{conv}
\ee
We endow $C_c(G)$ with its natural topology of LF-space. It is the union, over compact subsets $K\subset G$, of the Banach spaces 
\be
C_K(G)=\{b:G\to \cc \  \mbox{continuous}\ |\ \supp b\subset K\}
\ee
gifted with the supremum norm. The multiplication on $C_c(G)$ is separately continuous. $C_c(G)$ is also a complete bornological algebra for the bornology of bounded subsets: a subset $S\subset C_c(G)$ is small iff there is a compact $K\subset G$ such that $\supp b\subset K$ for any $b\in S$ and the supremum norm $\sup_{g\in G}|b(g)|$ is uniformly bounded over $S$. \\
Now let $M$ be a smooth complete Riemannian manifold, on which $G$ acts smoothly and properly by isometries, with compact quotient $G\backslash M$. We will introduce in this section the analytic assembly map 
\be
\mu:K_*^G(M) \to K_*(\Bc) \label{assem}
\ee 
from the equivariant $K$-homology of $M$ to the topological $K$-theory of a suitable Banach algebra completion $\Bc$ of the convolution algebra $C_c(G)$.\\

Let us first recall the definition of topological $K$-theory for arbitrary complete bornological algebras \cite{Me}. It is a generalization of the corresponding theory of Phillips for Fr\'echet $m$-algebras \cite{Ph}. Let $\Kc$ be the algera of ``smooth compact operators''. Its elements are infinite matrices $(A_{ij})_{i,j\in\nn}$ with rapidly decreasing entries in $\cc$, endowed with the family of submultiplicative norms
\be
\|A\|_n=\sup_{(i,j)\in\nn^2} (1+i+j)^n A_{ij} < \infty\qquad \forall n\in\nn\ .
\ee
$\Kc$ is a nuclear Fr\'echet algebra. In the following we shall identify $\Kc$ with a subalgebra of the algebra $\Lc(\Hc)$ of bounded operators on a separable Hilbert space $\Hc$, once a choice of basis is made. For any complete bornological algebra $\Bc$, denote by $\widetilde{\Bc\hotimes\Kc}$ the unitalization of the completed tensor product. Note that since $\Kc$ is nuclear, the bornological tensor product corresponds to the projective tensor product $\Bc\hotimes_{\pi}\Kc$ when $\Bc$ is a Fr\'echet algebra. In the general case, an element of the $K$-theory group $K_0(\Bc)$ is represented by an idempotent $e=e^2$ in the matrix algebra $ M_2(\widetilde{\Bc\hotimes\Kc})$, such that 
\be
e-\left(\begin{array}{cc}
                 0 & 0 \\
                 0 & 1 \\
                 \end{array} \right)\ \in M_2(\Bc\hotimes\Kc)
\ee
Two idempotents $e_0$ and $e_1$ are smoothly homotopic if they can be connected by a smooth path of idempotents. It means that there exists an idempotent relative to the tensor product $\Bc\hotimes \cinf[0,1]$ whose evaluation at $0$ and $1$ respectively coincide with $e_0$ and $e_1$. Smooth homotopy is an equivalence relation. The set of homotopy classes of such idempotents endowed with the natural operation of direct sum is the abelian group $K_0(\Bc)$.\\ 
Let $\cinf(S^1;\Bc):= \cinf(S^1)\hotimes\Bc$ be the algebra of smooth functions over the circle with values in $\Bc$. The product is given by pointwise multiplication and $\cinf(S^1)$ is gifted with its usual Fr\'echet topology. We define the odd $K$-theory group $K_1(\Bc)$ as the cokernel of the map $K_0(\Bc)\to K_0(\cinf(S^1;\Bc))$ induced by the constant homomorphism $\Bc\to\cinf(S^1;\Bc)$. These definitions coincide with the usual ones when $\Bc$ is a Banach algebra \cite{Bl}, in which case $K_0(\Bc)$ and $K_1(\Bc)$ are the only topological $K$-theory groups by Bott periodicity.\\
For any complete bornological algebra $\Bc$, the Chern character in entire cyclic homology is an additive map $\ch:K_i(\Bc)\to HE_i(\Bc)$, $i\in\zz_2$. It is constructed as follows \cite{Me}. An idempotent $e\in M_2(\widetilde{\Bc\hotimes\Kc})$ can be lifted to an idempotent $\eh$ of the analytic tensor algebra $\Tc(M_2(\widetilde{\Bc\hotimes\Kc}))$, see section \ref{sent}:
\be
\eh=e+ \sum_{k\ge 1} \frac{(2k)!}{(k!)^2} (e-\frac{1}{2})\otimes (e-e\otimes e)^{\otimes k}\ .
\ee
Then $\eh$ defines an entire cyclic homology class of even degree for the algebra $M_2(\widetilde{\Bc\hotimes\Kc})$. One has the isomorphisms
\be
HE_0(M_2(\widetilde{\Bc\hotimes\Kc}))\cong HE_0(\widetilde{\Bc\hotimes\Kc})\cong HE_0(\Bc\hotimes\Kc)\oplus \cc\ ,
\ee
and Morita invariance \cite{Me} implies $HE_0(\Bc\hotimes\Kc)\cong HE_0(\Bc)$. The image of $\eh$ in $HE_0(\Bc)$ is by definition the Chern character of the class $[e]\in K_0(\Bc)$. The Chern character on $K_1(\Bc)$ can be deduced, for example by excision \cite{Me2}. \\

The equivariant index theorem requires to complete the convolution algebra $C_c(G)$ of a locally compact group $G$ into a suitable Banach algebra. This will morally be the convolution algebra of integrable functions on $G$ with respect to an admissible measure, proportional to the Haar measure according to the following definition:
\begin{definition}\label{dadm}
Let $G$ be a locally compact group with right-invariant Haar measure $dg$. A measure $d\nu$ on $G$ is called \emph{admissible} if there is a real-valued, strictly positive and continuous function $\si$ on $G$ such that
\be
d\nu= \si\, dg \quad \mbox{and} \quad \si(gh)\leq \si(g)\si(h)\quad \forall g,h\in G\ .
\ee
It is easy to see that the $L^1$-norm $\|b\|=\int_Gd\nu\, |b(g)|$ associated to this measure is submultiplicative for the convolution product: $\|b_1b_2\|\leq \|b_1\|\|b_2\|$ for any $a,b\in C_c(G)$, and the corresponding Banach algebra $\Bc=L^1(G,d\nu)$ is called an \emph{admissible completion} of the convolution algebra.
\end{definition}
Note that the injection of bornological algebras $C_c(G)\hookrightarrow \Bc$ is bounded. The condition ``$\si$ is continuous'' is not essential but in practice interesting examples arise with this property:
\begin{example}\textup{For any locally compact group $G$, the Banach algebra $\Bc=L^1(G)$ of integrable functions with respect to the Haar measure $dg$ itself is an admissible completion of $C_c(G)$. As usual, the norm in $L^1(G)$ is given by $\|b\|_1=\int_G |b(g)|dg$.}
\end{example}
\begin{example}\label{edist}\textup{More generally, if $d: G\times G\to \rr_+$ is a right-invariant distance on $G$ (compatible with the topology), and $\al\in \rr_+$ is an arbitrary positive parameter, the admissible function $\si$ of definition \ref{dadm} may be chosen to be
\be
\si(g)=(1+d(g,1))^{\al}
\ee
hence, it grows as a power of the distance from $g$ to the unit element $1\in G$. The corresponding algebra $\Bc$ is therefore a subalgebra of $L^1(G)$, whose elements are integrable functions on $G$ with additional decay conditions at infinity. If moreover $G$ is abelian, the Fourier transform shows that $\Bc$ is isomorphic to an algebra of ``differentiable'' functions over the Pontrjagin dual $\widehat{G}$, with product given by pointwise multiplication. Increasing $\al$ improves the differentiability degree. The advantage of dealing with such a regular algebra of functions is that it is easy to construct cyclic cocycles, for example representing the fundamental class of $\widehat{G}$.}
\end{example}

We now turn to the description of the left hand side of the assembly map (\ref{assem}). We essentially follow Kasparov \cite{Kas}. Let $M$ be a smooth complete Riemannian manifold without boundary. Denote by $\diff(M)$ the group of diffeomorphisms $\phi:M\to M$, gifted with the topology of uniform convergence of $\phi$ and all its derivatives over compact subsets. An action of a locally compact group $G$ on $M$ is given by a continuous homomorphism $G\to\diff(M)$. For any $x\in M$ and $g\in G$, we denote by $gx\in M$ the (left) action of $g$ on $x$. The action is proper if the continuous map
\be
G\times M\to M\times M\ ,\qquad (g,x)\mapsto (x,gx)
\ee
is proper. In particular the isotropy subgroups of all points in $M$ are compact. The quotient space $X=G\backslash M$ is then a locally compact Hausdorff topological space. $M$ is called $G$-compact if $X$ is compact. When $M$ is a proper $G$-compact manifold, one can find a smooth, non-negative, compactly supported cut-off function $c$ over $M$, satisfying the identity
\be
\int_G dg \, c(gx)^2=1\qquad \forall x\in M\ . \label{cut}
\ee
Using the cut-off function we may obtain $G$-invariant objects over $M$, by an averaging procedure \cite{CM82}. In particular, $M$ can always be endowed with a $G$-invariant Riemannian metric. For if $\rho_0$ is any Riemannian metric, the averaged metric
\be
\rho(x)=\int_G dg \,c(gx)^2\rho_0^g(x)\ ,\qquad x\in M
\ee
is $G$-invariant. Here $\rho_0^g$ denotes the pullback of $\rho_0$ by the diffeomorphism induced by $g$. In what follows we will therefore always assume that $G$ acts by isometries on $M$. In the same way, a $G$-equivariant complex vector bundle of finite rank can always be endowed with a $G$-invariant hermitean structure.\\
Let $E_+\to M$ and $E_-\to M$ be two equivariant vector bundles. Denote by $\cinfc(E_{\pm})$ the space of smooth sections of $E_{\pm}$ with compact support. For any $m\in \rr$, let $\Psi^m_c(E_+,E_-)$ be the space of \emph{properly supported}, $G$-invariant pseudodifferential operators of order $m$, see \cite{H,Kas}. Such operators define linear maps $\cinfc(E_+)\to \cinfc(E_-)$ commuting with $G$. An operator is elliptic if its symbol is invertible outside the zero-section of the cotangent bundle $T^*M$. We will roughly define equivariant $K$-homology as the set of stable homotopy classes of $G$-invariant elliptic pseudodifferential operators of order zero. Stable homotopy is defined as follows. Two elliptic operators of order zero
\be
Q\in \Psi^0_c(E_+,E_-)\quad \mbox{and} \quad Q'\in \Psi^0_c(E_+',E_-')
\ee
are stably homotopic, if there exist $G$-equivariant vector bundles $F,F'$ with equivariant isomorphisms 
\be
\phi: E_+\oplus F\to E_+'\oplus F'\ ,\qquad \psi:E_-'\oplus F' \to E_-\oplus F\ ,
\ee
and a continuous homotopy of elliptic operators of order zero
\be
Q'\oplus 1_{F'} \sim \psi^*(Q\oplus 1_{F})\phi^*\ .
\ee
In other words, $Q$ and $Q'$ are homotopic modulo addition of a trivial operator. If $E$ is a $G$-equivariant vector bundle with $G$-invariant hermitean structure, we may consider elliptic operators $Q\in \Psi^0_c(E,E)$ which are \emph{selfadjoint} with respect to the $L^2$ inner product on $\cinfc(E)$. Stable homotopy of selfadjoint operators is defined in the obvious way.
\begin{definition} 
Let $G$ be a locally compact group, and $M$ a proper, $G$-compact complete Riemannian manifold without boundary on which $G$ acts by isometries. The equivariant K-homology group of even degree $K_0^G(M)$ is the set of stable homotopy classes of $G$-invariant elliptic operators $Q\in \Psi^0_c(E_+,E_-)$ between equivariant vector bundles $E_+,E_-$. The equivariant K-homology group of odd degree $K_1^G(M)$ is the set of stable homotopy classes of $G$-invariant selfadjoint elliptic operators $Q\in \Psi^0_c(E,E)$. The addition in $K_*^G(M)$ is induced by direct sum of vector bundles and operators.
\end{definition}
The zero element of $K_*^G(M)$ is given by the class of trivial operators, and the inverse of a class $[Q]$ corresponds to the class of a $G$-invariant parametrix for the elliptic operator $Q$. \\
Any $G$-invariant elliptic pseudodifferential operator $D_+\in \Psi^m_c(E_+,E_-)$ of order $m\neq 0$ also defines an element in $K_0^G(M)$ by taking the class of the elliptic operator of order zero
\be
Q=D_+\cdot \delta_m
\ee
where $\delta_m \in \Psi^{-m}_c(E_+,E_+)$ is the $G$-invariant elliptic operator with $G$-invariant symbol
\be
\si(x,p)=(1+\|p\|^2)^{-m/2}\ ,\quad \forall (x,p)\in T^*M\ .
\ee
Conventionally we will always consider such an operator $D_+$ together with its formal adjoint (in the sense of the $L^2$ inner product) $D_-\in \Psi^m_c(E_-,E_+)$, and the operator
\be
D=\left(\begin{array}{cc}
                                 0 & D_- \\
                                 D_+ & 0 \\
                                 \end{array} \right)
\ee
is a $G$-invariant, selfadjoint pseudodifferential operator of odd degree acting on the sections of the $\zz_2$-graded vector bundle $E=E_+\oplus E_-$. It thus determines a class $[D]\in K_0^G(M)$. Below we often consider that even-degree $K$-homology classes arise in this form.\\
In the case of $K_1^G(M)$, any selfadjoint elliptic operator $D\in \Psi_c^m(E,E)$ of order $m$ also defines an operator of order zero $Q=D\cdot\delta_m$, whence a class $[D]\in K_1^G(M)$. Here $E$ is of course trivially graded.
\begin{example}\label{ebott}\textup{Let $G=\zz^n$ act on $M=\rr^n$ by translations. The \emph{Bott element} is the $K$-homology class $\beta\in K_i^G(M)$ of degree $i\equiv n$ mod 2 corresponding to the Dirac operator $D$ acting on the sections of the trivial spinor bundle $S\to \rr^n$ ($S=S_+\oplus S_-$ is $\zz_2$-graded when the dimension $n$ is even and $D_+$ maps $\cinfc(S_+)$ to $\cinfc(S_-)$). }
\end{example}
If $G_1$ and $G_2$ are two groups acting respectively on the manifolds $M_1$ and $M_2$, there is an external product defined on $K$-homology. The construction is easier when the classes are represented by elliptic pseudodifferential operators of order one. If $D_1$ and $D_2$ are invariant elliptic operators of order one acting on the sections of two vector bundles on $E_1\to M_1$ and $E_2\to M_2$ respectively, their (graded) external product $D= D_1\otimes 1+ 1\otimes D_2$ is a $G_1\times G_2$ invariant elliptic operator of order one acting on the sections of the external tensor product $E_1\boxtimes E_2 \to M_1\times M_2$. The parity of $[D]$ is the sum of the parities of the $K$-homology classes $[D_1]$ and $[D_2]$. If both are odd, we may replace $E_1\boxtimes E_2$ by two copies of itself, endowed with its natural $\zz_2$-grading, and the graded sum of the differential operators reads
\be
D=\left(\begin{array}{cc}
                                 0 & D_1-iD_2 \\
                                 D_1+iD_2 & 0 \\
                                 \end{array} \right)\ .
\ee
This external product induces a bilinear map
\be
K_i^{G_1}(M_1)\boxtimes K_j^{G_2}(M_2)\to K_{i+j}^{G_1\times G_2}(M_1\times M_2)\ ,\quad i,j\in\zz_2
\ee
at the level of $K$-homology.\\

The analytic assembly map is then constructed as follows (see  \cite{C1}). We first consider the $K$-homology of even degree. So let $[Q]\in K_0^G(M)$ be represented by a $G$-invariant elliptic operator $Q\in \Psi^0_c(E_+,E_-)$ between the smooth sections of two vector bundles $E_{\pm}$ over $M$. The direct sum $E=E_+\oplus E_-$ is naturally a $\zz_2$-graded vector bundle. We will work in the $\zz_2$-graded algebra $\Pc=\cup_{m}\Psi^m_c(E,E)$ of all properly supported, $G$-invariant pseudodifferential operators on $\cinfc(E)$. Denote by $\Jc=\cap_{m}\Psi^m_c(E,E)$ the two-sided ideal of smoothing pseudodifferential operators in $\Pc$. Any operator $T\in \Jc$ has a smooth $G$-invariant kernel $T(x,y)$, such that its action on a section $\xi\in\cinfc(E)$ reads
\be
(T\xi)(x)=\int_M dy\, T(x,y)\xi(y)\ ,\quad x\in M\ ,\label{con}
\ee
where $dy$ is the $G$-invariant volume form associated to the Riemannian metric on $M$. We introduce the $2\times 2$ matrix notation in accordance with the $\zz_2$-grading on $E$:
\be
\Pc=\left(\begin{array}{cc}
                                 \Pc_{++} & \Pc_{+-} \\
                                 \Pc_{-+} & \Pc_{--} \\
                                 \end{array} \right)\ ,\quad \Jc=\left(\begin{array}{cc}
                                 \Jc_{++} & \Jc_{+-} \\
                                 \Jc_{-+} & \Jc_{--} \\
                                 \end{array} \right)\ .
\ee
Hence $Q\in \Pc_{-+}$. One can find a properly supported and $G$-invariant parametrix $P\in \Pc_{+-}$. Then $Q$ is an almost invertible operator, in the sense that
\be
PQ-1\in \Jc_{++}\ ,\qquad QP-1\in\Jc_{--}\ ,
\ee
and the matrix $F=\left(\begin{array}{cc}
                                 0 & P \\
                                 Q & 0 \\
                                 \end{array} \right)$ is an invertible element of the quotient algebra $\Pc/\Jc$. It follows from elementary algebraic $K$-theory \cite{M} that $Q$ has an index $\Ind(Q)$ in the algebraic group $K^{\mathrm{alg}}_0(\Jc)$, obtained as follows. There is a lift of $F$ to an invertible element $T\in\Pc$:
\be
T=\left(\begin{array}{cc}
                                 1 & 0 \\
                                 Q & -1 \\
                                 \end{array} \right)\left(\begin{array}{cc}
                                 1 & P \\
                                 0 & 1 \\
                                 \end{array} \right)\left(\begin{array}{cc}
                                 1 & 0 \\
                                 -Q & 1 \\
                                 \end{array} \right)=\left(\begin{array}{cc}
                                 1-PQ & P \\
                                 2Q-QPQ & QP-1 \\
                                 \end{array} \right)\ .
\ee
One has $T^{-1}=T$ and $T\equiv F$ mod $\Jc$. Denote by $\widetilde{\Jc}$ the algebra $\Jc$ augmented by adjoining units to the diagonal blocks $\Jc_{++}$ and $\Jc_{--}$. The idempotents
\be
e_1=T^{-1}\left(\begin{array}{cc}
  1 & 0 \\
  0 & 0 \\
\end{array} \right)T\in\widetilde{\Jc}\ ,\quad e_0=\left(\begin{array}{cc}
  0 & 0 \\
  0 & 1 \\
\end{array} \right)\in\widetilde{\Jc}
\ee
are such that the difference $e_1-e_0$ lies in $\Jc$, hence they define a $K$-theory class
\be
\Ind(Q)=[e_1]-[e_0]\in K^{\mathrm{alg}}_0(\Jc)\ .
\ee
The index $\Ind(Q)$ only depends on the image of $F$ in $\Pc/\Jc$, and not on the lift $T$. In particular, it is independent of the choice for the parametrix $P$.\\ 
Let $\Hc=L^2(E)$ be the $\zz_2$-graded Hilbert space of square-integrable sections of $E$, and let $\Lc(\Hc)$ denote the algebra of bounded operators on $\Hc$. Since the Riemannian metric and hermitean structure on $E$ are $G$-invariant, the elements of $G$, acting on the sections of $E$ by pullback, are represented by unitary operators on $\Hc$. We have a group (anti)homomorphism
\be
r: G\to \Uc(\Hc)\subset\Lc(\Hc)\ ,\qquad r(g)\xi=\xi^g \quad \forall g\in G\ ,\ \xi\in\Hc\ ,
\ee
and $r(gh)=r(h)r(g)$ for any $g,h\in G$. The superscript $^g$ denotes the pullback of sections by the diffeomorphism $g$. Note that $r$ is strongly continuous (i.e. the map $g\to r(g)\xi$ is continuous for any $\xi\in \Hc$), but not continuous for the operator norm on $\Lc(\Hc)$. Any smoothing operator $T\in\Jc$ being $G$-invariant, one has the commutation relation $Tr(g)=r(g)T, \forall g\in G$ in $\Lc(\Hc)$. Now remark that given a cut-off function $c$ on $M$ as defined in (\ref{cut}), the bounded operator on $\Hc$
\be
k_g=cTr(g)c \in \Lc(\Hc)\ ,
\ee
where $c$ acts by pointwise multiplication on $L^2(E)$, is a (non $G$-invariant) smoothing operator whose kernel $k_g(x,y)$ has compact support on $M\times M$, explicitly
\be
k_g(x,y)=c(x)T(gx,y)c(y)\quad \forall g\in G\ ,\ x,y\in M\ .
\ee
The algebra of (non $G$-invariant) smoothing kernels with compact support lies in an algebra of ``smooth compact operators'' $\Kc\subset\Lc(\Hc)$ as described before. Denote by $C_c(G;\Kc)$ the convolution algebra of continuous, compactly supported functions on $G$, with values in $\Kc$. We obtain an algebra homomorphism
\be
\te:\Jc\to C_c(G;\Kc)
\ee
by setting $\te(T)(g)=cTr(g)c$ for any $T\in\Jc$, $g\in G$. As a complete bornological vector space, $C_c(G;\Kc)$ is the union, over compact subsets $K\subset G$, of the Fr\'echet spaces $C_K(G;\Kc)$ of continuous functions with support in $K$. Since $\Kc$ is nuclear, $C_K(G;\Kc)$ is isomorphic to the projective tensor product $\Kc\hotimes_{\pi} C_K(G)$, and consequently the algebra $C_c(G;\Kc)$ corresponds to the bornological tensor product $\Kc\hotimes C_c(G)$. Therefore, we can take the push-forward of the index $\Ind(Q)$ by the induced homomorphism in $K$-theory
\be
\te_*:K^{\mathrm{alg}}_0(\Jc)\to K_0(C_c(G;\Kc))\cong K_0(C_c(G))\ .
\ee
The last isomorphism accounts for the stability relation $\Kc\hotimes\Kc\cong\Kc$. Finally, the image of the elliptic $G$-invariant pseudodifferential operator $Q$ by the assembly map $\mu$ is by definition
\be
\mu(Q)=\iota_*\te_* (\Ind(Q)) \ \in K_0(\Bc)\ ,
\ee
where $\iota_*:K_0(C_c(G))\to K_0(\Bc)$ is induced by the inclusion of $C_c(G)$ in any admissible completion $\Bc$.\\ 
In odd degree, we take advantage of the external product with the Bott element $\beta\in K_1^{\zz}(\rr)$,
\be
K_1^G(M)\stackrel{\beta\boxtimes}{\longrightarrow}K_0^{\zz\times G}(\rr\times M)\ee
to reduce the situation to the even-degree $K$-homology of the product manifold $\rr\times M$, on which the direct product of groups $\zz\times G$ acts. Now the convolution algebra of $\zz\times G$ is isomorphic to the (algebraic) tensor product $\cc\zz\otimes C_c(G)$, where $\cc\zz$ is the group ring of $\zz$. The above construction of the index for a representative $Q$ of an even degree class in $K_0^{\zz\times G}(\rr\times M)$ thus yields a $K$-theory element of the complete bornological algebra $\cc\zz\otimes C_c(G)$ (stabilized by $\Kc$)
\be
\te_*(\Ind(Q))\in K_0(\cc\zz\otimes C_c(G))\ .
\ee
Since $\cc\zz\otimes C_c(G)$ is a subalgebra of $\cinf(S^1;\Bc)=\cinf(S^1)\hotimes\Bc$, the image of $\te_*(\Ind(Q))$ in $K_1(\Bc)=\coker(K_0(\Bc)\to K_0(\cinf(S^1;\Bc)))$ is well-defined.
\begin{proposition}\label{dass}
The assembly map $\mu$ thus defined is compatible with stable homotopy of elliptic $G$-invariant pseudodifferential operators and descends to an additive map on equivariant $K$-homology 
\be
\mu:K_i^G(M)\to K_i(\Bc)\ ,\quad i\in\zz_2\ ,
\ee
for any admissible completion $\Bc$ of the convolution algebra of $G$. 
\end{proposition}
\begin{remark}\textup{When $G$ is a discrete countable group, the algebra $C_c(G;\Kc)$ is isomorphic to the algebraic tensor product $\Kc\otimes \cc G$. The analytic assembly map thus gives a topological $K$-theory element of the group ring $\cc G$, considered as an LF-algebra. When $G$ acts freely on $M$, this coincides with the construction of Connes and Moscovici \cite{CM90}.}
\end{remark}

\section{Bivariant Chern character}\label{sbiv}

Let $G$ be a locally compact group, and $M$ be a complete Riemannian manifold without boundary, endowed with a smooth, proper and $G$-compact action of $G$ by isometries. Since they are the main examples of interest, we will restrict from now on to equivariant $K$-homology classes $[D]\in K_*^G(M)$ represented by selfadjoint elliptic \emph{differential} operators of order one, acting on the smooth sections of some vector bundle $E\to M$. Such an operator $D$ has a unique extension to a selfadjoint (unbounded) operator on the Hilbert space $\Hc=L^2(E)$, from which we may obtain various functions like its modulus the heat operator $\exp(-sD^2)$, $s\geq 0$. In this section, we show that $D$ gives rise to an unbounded $\Ac$-$\Bc$-bimodule $(\Ec,\rho,D)$ according to the terminology of \cite{P1}, where $\Ac$ is the crossed-product algebra $\cinfc(M)\rtimes G$ and $\Bc$ is any admissible completion of the convolution algebra $C_c(G)$ (see definition \ref{dadm}). The bimodule $(\Ec,\rho,D)$ should be viewed as a ``smooth'' representative of a Kasparov bivariant $K$-theory class \cite{Bl}. Following the construction of \cite{P1}, we show that $(\Ec,\rho,D)$ has a Chern character in bivariant entire cyclic cohomology
\be
\ch(\Ec,\rho,D)\in HE_*(\Ac,\Bc)
\ee
based on the existence of the heat operator. This will play a central role in the equivariant index theorem.\\

Hence let $E\to M$ be a $G$-equivariant vector bundle endowed with an invariant hermitean structure, and let $D:\cinf(E)\to \cinf(E)$ be a $G$-invariant, selfadjoint elliptic differential operator of order one (for example a generalized Dirac operator \cite{BGV}). If the parity of the $K$-homology class $[D]$ is even, then $E$ is $\zz_2$-graded and $D$ is an operator of odd degree. We saw in section \ref{sass} that since $G$ acts isometrically on $M$, it is represented by unitary operators on the Hilbert space $\Hc=L^2(E)$ of square-integrable sections of $E$,
\be
r: G\to \Uc(\Hc)\ ,\qquad r(g)\xi=\xi^g\quad \forall g\in G\ ,\ \xi\in\Hc\ ,
\ee
and $r(gh)=r(h)r(g)$ $\forall g,h\in G$. Recall also that the function $r$ is strongly continuous (i.e. $\forall \xi\in\Hc$, $g\mapsto r(g)\xi$ is continuous), but not continuous for the operator norm in $\Lc(\Hc)$, unless $G$ is discrete. Also, $D$ extends to a self-adjoint unbounded operator on $\Hc$.\\
Consider the LF-space $C_c(G)$ of continuous, compactly supported functions over $G$, endowed with its bounded bornology: a subset $S\subset C_c(G)$ is small iff there is a compact subset $K\subset G$ such that all the functions $b\in S$ have support contained in $K$ and the supremum of $|b|$ over $K$ is uniformly bounded (see section \ref{sent}). Consider also the LF-space $\cinfc(M)$ of smooth functions with compact support on $M$. It is the union, over compact subsets $K\subset M$, of the Fr\'echet spaces $\cinf_K(M)$ of smooth fuctions with support contained in $K$. We provide $\cinfc(M)$ with its bounded bornology. Again, a subset $S\subset \cinfc(M)$ is small iff all $f\in S$ have support contained in a given compact $K$ in $M$ and all the derivatives of $f$ are uniformly bounded over $S$. $\cinfc(M)$ is a complete bornological algebra for the ordinary (pointwise) product of functions. By nuclearity of $\cinfc(M)$, the bornological tensor product of vector spaces $C_c(G)\hotimes \cinfc(M)$ is isomorphic to the inductive limit
\be
\Ac= \mathop{\limind}_{K\subset M} C_c(G; \cinf_K(M))\ ,
\ee
where $C_c(G;\cinf_K(M))$ is the space of continuous and compactly supported functions on $G$ with values in $\cinf_K(M)$. Endow $\Ac$ with the convolution product
\be
(a_1a_2)(g,x)=\int_G dh\, a_1(h,x)a_2(gh^{-1},h x)\ ,\quad \forall g\in G\ ,\ x\in M\ ,
\ee
for $a_i\in\Ac$. One checks that this product is bounded, hence turns $\Ac$ into a complete bornological algebra. We will note $\Ac$ as a crossed product 
\be
\Ac=\cinfc(M)\rtimes G\ .
\ee
Now let $\Bc$ be any admissible completion of the convolution algebra $C_c(G)$. Recall from definition \ref{dadm} that $\Bc=L^1(G,d\nu)$ for an admissible measure $d\nu=\si\, dg$. The elliptic operator $D$ gives rise to an unbounded $\Ac$-$\Bc$-bimodule according to the following definition (cf. \cite{P1}):
\begin{definition}\label{dunb}
Let $\Bc$ be any admissible completion of $C_c(G)$ and $\Ac$ be the crossed-product $\cinfc(M)\rtimes G$. Then any elliptic differential operator of order one $D:\cinfc(E)\to\cinfc(E)$ representing a $K$-homology class $[D]\in K_*^G(M)$ defines an \emph{unbounded $\Ac$-$\Bc$-bimodule} $(\Ec,\rho,D)$ as follows:
\begin{itemize}
\item $\Ec=\Hc\hotimes\Bc$ where $\Hc=L^2(E)$ is the Hilbert space of square-integrable sections of $E$. $\Ec$ is isomorphic to the Banach space $L^1(G,d\nu;\Hc)$ of $\Hc$-valued integrable functions on $G$ with respect to the admissible measure $d\nu$. It is endowed with an obvious right $\Bc$-module structure ($dh$ denotes the right-invariant Haar measure):
\be
(\xi b)(g)=\int_G dh\, \xi(h)b(gh^{-1})\ ,\quad \forall \xi\in\Ec\ ,\ b\in \Bc\ ,\ g\in G\ ,
\ee
and the module map $\Ec\times \Bc\to\Ec$ is bounded. We denote by $\End_{\Bc}(\Ec)$ the algebra of bounded endomorphisms of $\Ec$ commuting with the action of $\Bc$.
\item $\rho:\Ac\to\End_{\Bc}(\Ec)$ is the bounded algebra homomorphism given by
\be
(\rho(a)\xi)(g)=\int_G dh\, a(h)r(h)\cdot \xi(gh^{-1})\ ,\quad \forall a\in\Ac\ ,\ \xi\in\Ec\ ,\ g\in G\ ,
\ee
where for any $h\in G$, the smooth function $a(h)\in\cinfc(M)$ acting on the sections of $E$ by pointwise multiplication is considered as a bounded endomorphism of $\Hc$, and $r(h)\in\Uc(\Hc)$ is the unitary representation of $h$. Hence $\Ec$ is a left $\Ac$-module.
\item $D:\Ec\to\Ec$ is the unbounded operator with dense domain
\be
(D\xi)(g)=D(\xi(g))\ ,\quad \forall \xi\in\Ec\ ,\ g\in G\ ,
\ee
commuting with the right action of $\Bc$. The commutator $[D,\rho(a)]$ extends to an element of $\End_{\Bc}(\Ec)$ for any $a\in\Ac$.
\end{itemize}
\end{definition}
The last assertion comes from the fact that $D$, viewed as an operator on $\Hc$,  commutes with the representation of $G$ by hypothesis, and the commutator $[D,f]$ is bounded for any $f\in\cinfc(M)$. Also, note that the linear map $\Ac\to\End_{\Bc}(\Ec)$ induced by $a\mapsto [D,\rho(a)]$  is bounded. Finally, the triple $(\Ec,\rho,D)$ comes equipped with the same degree as the $K$-homology class $[D]$: it is even if $E$, $\Hc$ and $\Ec$ are $\zz_2$-graded (in which case $D$ is an odd operator and $\rho(a)$ is even for any $a\in\Ac$), and odd if $E$ is trivially graded.\\

The bimodule $(\Ec,\rho,D)$ has the required properties to apply the construction of the bivariant Chern character along the lines of \cite{P1}. The essential point is the existence of the heat operator $\exp(-tD^2)$ for any $t\geq 0$. We will give an explicit formula of JLO type \cite{JLO} yielding a bivariant entire cyclic cohomology class
\be
\ch(\Ec,\rho,D)\in HE_*(\Ac,\Bc)
\ee
 of the same degree as the class $[D]\in K_*^G(M)$. We first introduce for any $p\in[1,\infty)$ the Schatten ideal $\ell^p=\ell^p(\Hc)$ of $p$-summable operators on $\Hc$. It is a Banach algebra for the norm $\|x\|_p=(\Tr |x|^p)^{1/p}$, and the injections $\ell^p\to\ell^q$ for $p<q$ are continuous (bounded). Let us concentrate on the ideal $\ell^1$ of trace-class operators on $\Hc$. Denote by $C_c(G;\ell^1)$ the LF-space of compactly supported continuous functions on $G$ with values in $\ell^1$. Endow it with the convolution product (as usual $dh$ denotes the right-invariant Haar measure)
\be
(a_1a_2)(g)=\int_G dh\, a_1(h)a_2(gh^{-1})\ ,\quad \forall a_i\in C_c(G;\ell^1)\ .\label{toto}
\ee
Then $C_c(G;\ell^1)$ becomes a complete bornological algebra, and we obtain a bounded injective homomorphism $C_c(G;\ell^1)\hookrightarrow \End_{\Bc}(\Ec)$ by specifying the action on $\Ec$:
\be
(a\xi)(g)=\int_G dh\, a(h)\cdot \xi(gh^{-1})\ ,\quad \forall a \in C_c(G;\ell^1)\ ,\ \xi\in \Ec\ .\label{tata}
\ee
As an intermediate algebra between $C_c(G;\ell^1)$ and $\End_{\Bc}(\Ec)$, we may consider the space of $\ell^1$-valued integrable functions $L^1(G,d\nu;\ell^1)$ with respect to the admissible measure $d\nu$, endowed with a convolution product and an action on $\Ec$ also given by equations (\ref{toto},\ref{tata}). Now remark that this algebra is isomorphic to the projective tensor product of the Banach algebras $\ell^1$ and $L^1(G,d\nu)=\Bc$, which also corresponds to the bornological tensor product (see section \ref{sent})
\be
 L^1(G,d\nu;\ell^1)= \ell^1\hotimes \Bc\ ,
\ee
and a rapid inspection shows that its action on $\Ec=\Hc\hotimes\Bc$ decomposes as the representation of $\ell^1$ on $\Hc$ and the left multiplication of $\Bc$ on itself. We summarise with the following sequence of inclusions
\be
C_c(G;\ell^1)\hookrightarrow \ell^1\hotimes \Bc \hookrightarrow \End_{\Bc}(\Ec)\ .
\ee
Let us now have a closer look at the homomorphism $\rho:\Ac\to \End_{\Bc}(\Ec)$. For any $\xi\in \Ec$ and $a\in\Ac$, the left product $\rho(a)\xi\in\Ec$ reads
\be
(\rho(a)\xi)(g)=\int_G dh\, a(h)r(h)\cdot \xi(gh^{-1})\ ,
\ee
where the function $h\in G\mapsto a(h)r(h)\in\Lc(\Hc)$ is only strongly continuous (as $h\mapsto r(h)$ is) and not continuous for the operator norm in $\Lc(\Hc)$. However, the heat operator $\exp(-sD^2)$ is smoothing for $s>0$ and $a(h)$ is a smooth function with compact support on $M$ for any $h$. It follows that the product $\exp(-sD^2) a(h)r(h)\in \Lc(\Hc)$ is actually a trace-class operator when $s>0$, and the function 
\be
 h\in G \mapsto e^{-sD^2}a(h)r(h)\in \ell^1\ ,\quad s>0
\ee
is continuous with compact support on $G$, hence yields an element of $C_c(G;\ell^1)$. It is not hard to show that the operator $\exp(-sD^2)\rho(a)$ on $\Ec$ precisely corresponds to this element, hence we can write
\be
e^{-sD^2}\rho(a) \in C_c(G; \ell^1)\subset \End_{\Bc}(\Ec)\ ,\quad \forall a\in \Ac\ ,\ s>0\ .
\ee
In the same way, one has 
\be
e^{-sD^2}[D,\rho(a)]\in C_c(G;\ell^1) \subset \End_{\Bc}(\Ec)\quad \forall a\in \Ac\ ,\ s>0
\ee
since $e^{-sD^2}[D,\rho(a)]$ corresponds to the function $h\mapsto e^{-sD^2}[D,a(h)]r(h)$ (recall that $D$ is $G$-invariant and therefore commutes with the operator $r(h)$). Now we fix a parameter $t>0$ and rescale the unbounded operator $D$ as $tD$. Following \cite{P1}, the new bimodule $(\Ec,\rho,tD)$ gives rise to a chain map from the $(b,B)$-bicomplex $\Om\Ac$ to the $X$-complex $X(\Bc)$, with components
\be
\chi^n(\Ec,\rho,tD):\Om^n\Ac\to X(\Bc)
\ee
defined in any degree $n\in\nn$ (see section \ref{sent} for the definitions of these complexes). To that end, if $A_1,\ldots, A_n$ are (not necessarily bounded) operators on $\Hc$, we introduce the notation
\be
\langle A_1,\ldots, A_n\rangle_t = \int_{\Delta_{n}}ds_1\ldots ds_{n}\, e^{-s_0 t^2D^2}A_1e^{-s_1 t^2D^2}\ldots A_n e^{-s_{n} t^2D^2}\ ,
\ee
where $\Delta_{n}$ is the standard $n$-simplex with coordinates $s_i\geq 0$, $\sum_is_i=1$. The integration over $\Delta_n$ is supposed to make sense whenever we use it. For a typical example, if $a_0,\ldots,a_n$ denote elements of the algebra $\Ac$, the preceding discussion shows that
\be
\langle \rho(a_0),[D,\rho(a_{1})],\ldots ,[D,\rho(a_n)] \rangle_t \in C_c(G;\ell^1)\subset \End_{\Bc}(\Ec)\ .
\ee
The chain map $\chi(\Ec,\rho,tD)$ is constructed as follows. Recall that the $X$-complex of $\Bc$ splits as the direct sum of $\Bc$ in degree zero and $\Om^1\Bc_{\nat}$ in degree one. Each component $\chi^n$ therefore splits into two parts. The first one $\chi^n_0:\Om^n\Ac\to \Bc$ is defined whenever $n$ has the same parity as the bimodule $\Ec$, and given by its evaluation on a $n$-form $a_0da_1\ldots da_n\in \Om^n\Ac$:
\beq
\lefteqn{\chi^n_0(a_0da_1\ldots da_n) = (-t)^n\sum_{i=0}^n (-)^{i(n-i)}\times } \label{joe1}\\
&&\tau \langle [D,\rho_{i+1}],\ldots ,[D,\rho_n] ,\rho_0,[D,\rho_1],\ldots, [D,\rho_i]\rangle_t\ ,\non
\eeq
where we use the abbreviation $\rho_i=\rho(a_i)$. Here $\tau:\ell^1\to \cc$ is the (super)trace of operators depending on the parity of $(\Ec,\rho,D)$: when the Hilbert space $\Hc=\Hc_+\oplus\Hc_-$ is $\zz_2$-graded (even case), then $\tau=\Tr_s=\Tr_{\Hc_+}-\Tr_{\Hc_-}$ is the supertrace of operators in $\ell^1$. When $\Hc$ is trivially graded (odd case), then $\tau=\sqrt{2i}\, \Tr$ is proportional to the usual trace. The factor of $\sqrt{2i}$ is conventionally introduced for consistency with Bott periodicity \cite{P1}. After taking the trace, the right hand side of (\ref{joe1}) thus yields a continuous scalar-valued function with compact support on $G$, viewed as an element of the algebra $\Bc$ via the inclusion $C_c(G)\hookrightarrow\Bc$.\\ 
The other components $\chi^{n+1}_1:\Om^{n+1}\Ac\to \Om^1\Bc_{\nat}$ are also defined when $n$ has the same parity as $\Ec$. We first have to extend slightly the bimodule $\Ec$ by considering the unitalisation $\Bct=\Bc\oplus \cc$ of $\Bc$. Consider the obvious right $\Bct$-module
\be
\Ect=\Hc\hotimes \Bct\ .
\ee
Then as a bornological vector space, $\Ect=\Ec\oplus \Hc$. We want to endow $\Ect$ with a left $\Ac$-module structure. The action of $\Ac$ onto the first summand $\Ec$ is already defined, so we need to specify its action on the second summand $\Hc$. For any $a\in\Ac$, define the map $\rho(a):\Hc\to \Ec$ by
\be
(\rho(a)\xi)(g)=a(g)r(g)\cdot \xi\ ,\quad \forall \xi\in\Hc\ ,\ g\in G\ .
\ee
It is easy to check that $\rho(a)\in\End_{\Bct}(\Ect)$, so that $(\Ect,\rho,D)$ is an unbounded $\Ac$-$\Bct$-bimodule as required (the action of $D$ on the summand $\Hc$ is clear). Remark we have
\be
\Ec=\Ect\hotimes_{\Bct}\Bc
\ee
as an $\Ac$-$\Bc$-bimodule. Now consider the space $\Om^1\Bc$ of noncommutative one-forms over $\Bc$. It can be viewed as a left $\Bct$-module and right $\Bc$-module, together with the universal derivation $\dd:\Bc\to\Om^1\Bc$. Then define the following $\Ac$-$\Bc$-bimodule
\be
\Om^1\Ec=\Ect\hotimes_{\Bct}\Om^1\Bc\ .
\ee
As a bornological vector space, $\Om^1\Ec$ is isomorphic to the tensor product $\Hc\hotimes\Om^1\Bc$. This shows the existence of a canonical flat connection on $\Ec$
\be
\dd:\Ec\to \Om^1\Ec
\ee
induced by the universal derivation on $\Bc$. Moreover, any element of $\ell^1\hotimes \Om^1\Bc$ acts as a right $\Bc$-module map $\Ec\to \Om^1\Ec$ in the obvious way. For any $a\in\Ac$, the commutator
\be
\dd\rho(a)=[\dd,\rho(a)]: \Ec\to \Om^1\Ec
\ee
is also a right $\Bc$-module map, and after multiplication by the heat operator we get
\be
e^{-sD^2}\dd\rho(a) \in \ell^1\hotimes \Om^1\Bc\ ,\quad s>0\ .
\ee
It is possible to give a more precise description of this element. Recall that $\Om^1\Bc$ is by definition the space $\Bct\hotimes\Bc=\Bc\oplus\Bc^{\hotimes 2}$, and the first summand of the r.h.s. corresponds to the image of the derivation $\dd$. Since $\Bc=L^1(G,d\nu)$, the tensor product $\Bc\hotimes\Bc$ is isomorphic to the space of integrable functions over $G\times G$ with respect to the product measure $d\nu\times d\nu$. Therefore, $\Om^1\Bc$ canonically contains the subspace of continuous functions with compact support
\be
\Om^1_c(G):= C_c(G)\oplus C_c(G\times G)\subset \Om^1\Bc\ ,
\ee
and analogously
\be
\Om^1_c(G;\ell^1):=C_c(G;\ell^1)\oplus C_c(G\times G;\ell^1)\subset \ell^1\hotimes\Om^1\Bc\ .
\ee
Then $e^{-sD^2}\dd\rho(a)$ is an element of the first summand $C_c(G;\ell^1)$. All these notations  being fixed, the component $\chi^{n+1}_1:\Om^{n+1}\Ac\to \Om^1\Bc_{\nat}$ is given by
\beq
\lefteqn{\chi^{n+1}_1(a_0da_1\ldots da_{n+1}) = (-t)^n\sum_{i=1}^{n+1} (-)^{i(n-i+1)}\times} \label{joe2}\\
&&\nat \tau(\langle [D,\rho_{i+1}],\ldots , [D,\rho_{n+1}],\rho_0,[D,\rho_1],\ldots, [D,\rho_{i-1}]\rangle_t \dd \rho_i) \ .\non
\eeq
Here, the term under $\tau$ lies in $\Om^1_c(G;\ell^1)$. Hence taking the trace on $\ell^1$ yields an element of $\Om^1_c(G)$, or $\Om^1\Bc$, then composing with the projection $\nat:\Om^1\Bc\to\Om^1\Bc_{\nat}$ gives the desired map. Collecting all the components $\chi^n_0$ and $\chi^{n+1}_1$, we obtain a linear map from the $(b+B)$-complex $\Om\Ac$ to the $X$-complex of $\Bc$, whose degree coincide with the parity of $(\Ec,\rho,tD)$. The algebraic manipulations performed in \cite{P1} show that it is actually a chain map:
\be
\pm (\bb\oplus \nat\dd)\circ \chi^n= \chi^{n-1}\circ b+ \chi^{n+1}\circ B\ .
\ee
The sign in the l.h.s depends on the parity. The following proposition is a consequence of the crucial property that $\Bc$ is a Banach algebra.
\begin{proposition}\label{pchi}
Let $(\Ec,\rho,D)$ be an unbounded $\Ac$-$\Bc$-bimodule as in Definition \ref{dunb}. Then for any $t>0$, the collection of components $\chi^n(\Ec,\rho,tD):\Om^n\Ac\to X(\Bc)$ extends to an entire cyclic cocycle, i.e. a bounded chain map
\be
\chi(\Ec,\rho,tD):\Ome\Ac\to X(\Bc)
\ee
from the $(b+B)$-complex of entire chains over $\Ac$ to the $X$-complex of $\Bc$. The degree of $\chi(\Ec,\rho,tD)$ in the $\zz_2$-graded complex $\hom(\Ome\Ac,X(\Bc))$ coincides with the parity of $(\Ec,\rho,D)$.
\end{proposition}
{\it Proof:} We first show that the collection of components $\chi_0=(\chi^n_0)_{n}$ extends to a bounded linear map $\Ome\Ac\to \Bc$. Recall that $\Ome\Ac$ is the completion of $\Om\Ac$ for the bornology generated by the subsets
$$
\bigcup_{n\ge 0} [n/2]!\,\St (dS)^n\quad ,\ S\in\Sg(\Ac)\ ,
$$
with $\St=\{1\}\cup S$. Hence we must show that, given a small subset $S\subset\Ac$, $\chi^n_0$ maps $[n/2]!\,\St (dS)^n$ to a small subset of $\Bc$ which does not depend on $n$. So we fix a small $S\subset\Ac$. Since $\Ac$ is the inductive limit, over compact subsets $K\subset M$, of the spaces $C_c(G;\cinf_K(M))$, we have the following description of $S$: there exist two compact subsets $K\subset M$ and $L\subset G$, such that any function $a\in S$ has support contained in $L$ and the values $a(g)$, $g\in G$ are smooth functions on $M$ with support in $K$ and all derivatives bounded. In particular, one can find another compact subset $K'\subset M$ which contains $K$ and all its images $g\cdot K$ by the diffeomorphisms $g\in L$. It follows that, if $f\in \cinfc(M)$ is a real-valued function such that $f(x)=1$ for any $x\in K'$ and $0\leq f(x)\leq 1$ everywhere (we call $f$ a ``plateau'' function), we have 
$$\rho(a)=f\rho(a)=\rho(a)f\ \mbox{in}\ \End_{\Bc}(\Ec)$$ 
for any element $a\in S$ (use $\rho(a)(g)=a(g)r(g)$). Note that $f$ viewed as an operator in $\Lc(\Hc)$ by pointwise multiplication, is selfadjoint with spectrum contained in the interval $[0,1]$. The plateau function $f$ will be useful for dealing with the non-compacity of the manifold $M$. Now, given $n+1$ elements $a_i\in S$, the function $\chi^n_0(a_0da_1\ldots da_n)\in C_c(G)$ is given by Equation (\ref{joe1}). We have to evaluate the supremum norm, over $G$, of each of the $n+1$ terms appearing in the sum of the r.h.s. Let us consider for example the last term of this sum, up to a sign:
$$
I=t^n\int_{\Delta_{n+1}}ds\, \tau(e^{-s_0t^2D^2}\rho_0 e^{-s_1t^2D^2} [D,\rho_1] e^{-s_2t^2D^2}\ldots [D,\rho_n] e^{-s_{n+1}t^2D^2})\ ,
$$
where $\rho_i=\rho(a_i)$. The integrand at a point $s\in\Delta_{n+1}$ is an element of $C_c(G)$, and its evaluation at a point $g\in G$ can be estimated using the trace norm $\|\cdot\|_1$ on $\ell^1(\Hc)$ as follows:
\beq
\lefteqn{|\tau(e^{-s_0t^2D^2}\rho_0 e^{-s_1t^2D^2} [D,\rho_1] e^{-s_2t^2D^2}\ldots [D,\rho_n] e^{-s_{n+1}t^2D^2})|(g)}\non\\
& \leq & \int_{G^n} dh_0\ldots dh_{n-1} \, \| e^{-s_0t^2D^2}\rho_0(h_0) e^{-s_1t^2D^2} [D,\rho_1](h_1) e^{-s_2t^2D^2}\non\\
&&\qquad \qquad \qquad \qquad \ldots [D,\rho_n](gh_0^{-1}\ldots h_{n-1}^{-1}) e^{-s_{n+1}t^2D^2}\|_1\ .\non
\eeq
(When the trace $\tau$ corresponds to a bimodule of odd degree, we should take care of an overall factor of $\sqrt{2i}$). Taking into account that $\rho_i=f\rho_i f$ and $[D,\rho_i]=f[D,\rho_i]f$, the H\"older inequality implies
\beq
\lefteqn{\| e^{-s_0t^2D^2}\rho_0(h_0) e^{-s_1t^2D^2} [D,\rho_1](h_1)  \ldots [D,\rho_n](gh_0^{-1}\ldots h_{n-1}^{-1}) e^{-s_{n+1}t^2D^2}\|_1}\non\\
&\leq&  \|\rho_0(h_0)\|_{\infty} \| f e^{-s_1t^2D^2}f\|_{s_1^{-1}} \|[D,\rho_1](h_1)\|_{\infty} \| f e^{-s_2t^2D^2}f\|_{s_2^{-1}} \ldots\non\\
&& \qquad \ldots \|[D,\rho_n](gh_0^{-1}\ldots  h_{n-1}^{-1})\|_{\infty} \|f e^{-(s_{n+1}+s_0)t^2D^2}f\|_{(s_{n+1}+s_0)^{-1}} \ ,\non
\eeq
where for any $i$ the operator $f e^{-s_it^2D^2} f$ is an element of the Schatten ideal $\ell^{s_i^{-1}}$ with norm $\|\cdot\|_{s_i^{-1}}$, and $\|\cdot \|_{\infty}$ is the operator norm. The compression by the plateau function $f\in\cinfc(M)$ is necessary because $e^{-s_it^2D^2}$ alone is a priori not a compact operator. Now observe that since the heat operator and $f$ are selfadjoint and positive, one has
$$
\|f e^{-s_it^2D^2}f\|_{s_i^{-1}} = \big(\Tr \big((f e^{-s_i t^2D^2}f)^{s_i^{-1}}\big)\big)^{s_i}\ .
$$
Moreover, the function $x\mapsto x^{s_i^{-1}}$ being convex for $s_i^{-1}\geq 1$, and $f=f^*\leq 1$, Proposition 4.6 (ii) of \cite{FK} allows to write\footnotemark[1]\footnotetext[1]{I wish to thank T. Fack for pointing ref. \cite{FK} to my attention.}$$
\Tr \big((f e^{-s_i t^2D^2}f)^{s_i^{-1}}\big)\leq \Tr (f e^{-t^2D^2}f)
$$
\noindent so that $\| f e^{-s_it^2D^2}f\|_{s_i^{-1}}\leq \big(\Tr (f e^{-t^2D^2}f)\big)^{s_i}$. Introduce now the following positive functions over $G$
$$
b_0=\| \rho_0\|_{\infty}\ ,\quad b_i=\|[D,\rho_i]\|_{\infty}\ ,\ i\geq 1\ .
$$
They are compactly supported and continuous, because $b_0(g)=\|a_0(g)r(g)\|_{\infty}=\|a_0(g)\|_{\infty}$ and $b_i(g)=\|[D,a_i(g)]r(g)\|_{\infty}=\|[D,a_i(g)]\|_{\infty}$ (recall that $r(g)\in \Uc(\Hc)$). Also, all the $b_i$'s are contained in a small subset $T_S$ in the bornology of the LF-space $C_c(G)$. Remark that $T_S$ depends only on $S\subset \Ac$. We can write
\beq
\lefteqn{|\tau(e^{-s_0t^2D^2}\rho_0 e^{-s_1t^2D^2} [D,\rho_1] e^{-s_2t^2D^2}\ldots [D,\rho_n] e^{-s_{n+1}t^2D^2})|(g)}\non\\
&\leq& \Tr (f e^{-t^2D^2}f)\int_{G^n} dh_0\ldots dh_{n-1} \, b_0(h_0)b_1(h_1)\ldots b_n(gh_0^{-1}\ldots  h_{n-1}^{-1})\non\\
&\leq& \Tr (f e^{-t^2D^2}f) (b_0\ldots b_n)(g)\ ,\non
\eeq
so that taking the integral over the simplex $\Delta_{n+1}$ brings a factor of $1/(n+1)!$ in the estimate
$$
|I|\ \leq\ t^n\, \frac{\Tr (f e^{-t^2D^2}f)}{(n+1)!}\, b_0\ldots b_n\ .
$$
This is an inequality of continuous, compactly supported functions over $G$. Since the inclusion $C_c(G)\to \Bc$ is bounded, $T_S$ is also a small subset of the Banach algebra $\Bc$. Let $\|\cdot\|$ denote the norm of $\Bc$, and put $\la_S=\|T_S\|$. One has $\|b_0\ldots b_n\|\leq \la_S^{n+1}$, so that the above inequality implies
$$
\|I\|\ \leq\ t^n\, \frac{\Tr (f e^{-t^2D^2}f)}{(n+1)!}\, \la_S^{n+1}\ .
$$
The same estimate holds for all the $n+1$ terms of the sum (\ref{joe1}), therefore 
$$
\|[n/2]!\, \chi^n_0(a_0da_1\ldots da_n)\|\ \leq\ \Tr (f e^{-t^2D^2}f)\frac{[n/2]!}{n!}\,t^n \la_S^{n+1}\ .
$$
The function $n\mapsto \frac{[n/2]!}{n!}\,t^n \la_S^{n+1}$ being bounded as $n\to \infty$, the collection $\chi_0$ maps the union $\bigcup_{n\ge 0} [n/2]!\,\St (dS)^n$ to a bounded subset of $\Bc$, as required. It follows that $\chi_0$ is bounded for the entire bornology on $\Om\Ac$ and extends to a bounded map on the completion $\Ome\Ac$.\\
Let us now concentrate on the odd components $\chi_1^{n+1}:\Om^{n+1}\Ac\to \Om^1\Bc_{\nat}$. We will perform the same kind of estimates, with the only difference that now the target space is the direct sum $\Om^1_c(G)=C_c(G)\oplus C_c(G\times G)$. The function of two variables $\chi_1^{n+1}(a_0da_1\ldots da_{n+1})\in C_c(G\times G)$ is given by formula (\ref{joe2}), and taking for example the last term of this sum one gets
$$
J=t^n\int_{\Delta_{n+1}}ds\, \tau(e^{-s_0t^2D^2}\rho_0 e^{-s_1t^2D^2} [D,\rho_1] e^{-s_2t^2D^2}\ldots [D,\rho_n] e^{-s_{n+1}t^2D^2}\dd\rho_{n+1})\ ,
$$
with $\rho_i=\rho(a_i)$, $a_i\in S\subset \Ac$. Proceeding as in the previous case, the integrand of $J$ can be estimated by evaluation on two points $g_0,g_1\in G$:
\beq
\lefteqn{|\tau(e^{-s_0t^2D^2}\rho_0 e^{-s_1t^2D^2} [D,\rho_1] e^{-s_2t^2D^2}\ldots [D,\rho_n] e^{-s_{n+1}t^2D^2}\dd\rho_{n+1})|(g_0,g_1)}\non\\
&\leq& \int_{G^n} dh_0\ldots dh_{n-1} \, \| e^{-s_0t^2D^2}\rho_0(h_0) e^{-s_1t^2D^2} [D,\rho_1](h_1) e^{-s_2t^2D^2}\non\\
&&\qquad \qquad \qquad \qquad \ldots [D,\rho_n](g_0h_0^{-1}\ldots h_{n-1}^{-1}) e^{-s_{n+1}t^2D^2}\rho_{n+1}(g_1)\|_1 \non\\
&\leq& \Tr (f e^{-t^2D^2}f)\, (b_0\ldots b_n)(g_0)\, b_{n+1}(g_1)\non\\
&\leq& \Tr (f e^{-t^2D^2}f)\, (b_0\ldots b_n\dd b_{n+1})(g_0,g_1)\ ,\non
\eeq
with $b_0=\| \rho_0\|_{\infty}$, $b_i=\|[D,\rho_i]\|_{\infty}$ for $1\le i\le n$ and $b_{n+1}=\| \rho_{n+1}\|_{\infty}$. Again, all the $b_i$'s are contained in a small subset $T_S\subset C_c(G)$ depending only on $S$. We denote by $\|\cdot \|^{\hotimes 2}$ the projective norm on the tensor product of Banach spaces $\Bc\hotimes\Bc\subset \Om^1\Bc$. Its evaluation on the function of two varables $b_0\ldots b_n \dd b_{n+1}\in C_c(G\times G)$ gives
$$
\|b_0\ldots b_n \dd b_{n+1}\|^{\hotimes 2}= \|b_0\ldots b_n \|\cdot \| b_{n+1}\|\leq \la_S^{n+2}\ ,
$$
with $\la_S=\|T_S\|$. Therefore, the norm of $J$ verifies the following bound:
$$
\|J\|^{\hotimes 2}\ \leq\ t^n\, \frac{\Tr (f e^{-t^2D^2}f)}{(n+1)!}\,\la_S^{n+2}\ .
$$
The same estimate holds for all the terms of the sum (\ref{joe2}), so that
$$
\|[(n+1)/2]! \, \chi^{n+1}_1(a_0da_1\ldots da_{n+1})\|^{\hotimes 2}\ \leq\ \Tr (f e^{-t^2D^2}f)\frac{[(n+1)/2]!}{n!}\,t^n\la_S^{n+2},
$$
which shows that the collection $\chi_1$ maps $\bigcup_{n\geq 0}[n/2]!\,\St (dS)^{n}$ to  a small subset of $\Om^1\Bc_{\nat}$. Consequently, $\chi$ extends to a bounded linear map $\Ome\Ac\to X(\Bc)$. The fact that it is still a cocycle is a consequence of the universal properties of bornological completions. \cqfd\\

The cohomology class of the cocycle $\chi(\Ec,\rho,tD)$ enjoys good properties. It is in particular invariant under suitable ``smooth'' homotopies of the homomorphism $\rho$ and the operator $D$ \cite{P1}. For our purpose we need only to establish the invariance with respect to changes of the parameter $t$. This is achieved by introducing the Chern-Simons transgressions $cs^n(\Ec,\rho,tD):\Om^n\Ac\to X(\Bc)\hotimes\Om^1\rr_+^*$. Here $\Om^1\rr_+^*$ is the space of smooth one-forms with respect to $t\in\rr_+^*$. The transgressions are given by formulas analogous to (\ref{joe1},\ref{joe2}), with insertions of the operator-valued one-form $dt\,D$. In each degree $n$ of parity opposite to $(\Ec,\rho,D)$, the component $cs^n_0:\Om^n\Ac\to \Bc\hotimes\Om^1\rr_+^*$ reads
\beq
\lefteqn{cs^n_0(a_0da_1\ldots da_n) = (-t)^n\sum_{i, j}(-)^{i(n-i)} \tau \langle [D,\rho_{i+1}],\ldots ,[D,\rho_j], dt D, \ldots }\non\\
&&\qquad\qquad\qquad\qquad\qquad \ldots ,[D,\rho_n],\rho_0,[D,\rho_1],\ldots, [D,\rho_{i}]\rangle_t\ ,\label{eoj1}
\eeq
whereas the component $cs^{n+1}_1:\Om^{n+1}\Ac\to \Om^1\Bc_{\nat}\hotimes\Om^1\rr_+^*$ is
\beq
\lefteqn{cs^{n+1}_1(a_0da_1\ldots da_{n+1}) = (-t)^n\sum_{i, j} (-)^{i(n-i+1)} \nat \tau(\langle [D,\rho_{i+1}],\ldots ,[D, \rho_j],dt D, } \non\\
&& \qquad \qquad\qquad\qquad\qquad \ldots ,[D,\rho_{n+1}],\rho_0,[D,\rho_1],\ldots, [D,\rho_{i-1}]\rangle_t \dd \rho_i) \ .\label{eoj2}
\eeq
The Chern-Simons forms define a linear map from the $(b+B)$-complex $\Om\Ac$ to the $X$-complex $X(\Bc)\hotimes\Om^1\rr_+^*$, of degree opposite to the parity of $(\Ec,\rho,D)$. It is shown in \cite{P1} that the derivative of $\chi(\Ec,\rho,tD)$ with respect to $t$ is the coboundary of $cs(\Ec,\rho,tD)$ in the complex $\hom(\Om\Ac,X(\Bc)\hotimes\Om^1\rr_+^*)$:
\be
d_t \chi^n=-(\bb\oplus\nat\dd)\circ  cs^n\pm (cs^{n-1}\circ b+ cs^{n+1}\circ B)\ .
\ee
The sign $\pm$ depends on the parity of the cochains. Again, the properties of the Banach completion $\Bc$ ensure that $cs(\Ec,\rho,tD)$ extends to an entire cochain:
\begin{proposition}\label{pcs}
Let $(\Ec,\rho,D)$ be an unbounded $\Ac$-$\Bc$-bimodule as in Definition \ref{dunb}. Then for any $t>0$, the collection of components $cs^n(\Ec,\rho,tD):\Om^n\Ac\to X(\Bc)$ extends to an entire cochain, i.e. a bounded linear map
\be
cs(\Ec,\rho,tD):\Ome\Ac\to X(\Bc)\hotimes\Om^1\rr_+^*
\ee
whose coboundary in the complex $\hom(\Ome\Ac, X(\Bc)\hotimes\Om^1\rr_+^*)$ equals the derivative of the cocycle $\chi(\Ec,\rho,tD)$ with respect to $t$. Hence, the entire cyclic cohomology class of $\chi(\Ec,\rho,tD)$ is independent of $t$.
\end{proposition}
{\it Proof:} It follows from the same estimates used with the cocycles $\chi(\Ec,\rho,tD)$, the only difference being the insertion of the operator $D$ in formulas (\ref{eoj1}, \ref{eoj2}) defining the Chern-Simons transgressions. For example, let us pick up the following term in the sum (\ref{eoj1}), where $D$ is at position $i$:
$$
t^n\int_{\Delta_{n+2}}ds\, \tau(e^{-s_0t^2D^2}\rho_0 e^{-s_1t^2D^2} [D,\rho_1] \ldots D e^{-s_it^2D^2}\ldots [D,\rho_n] e^{-s_{n+2}t^2D^2})\ ,
$$
$\rho_j=\rho(a_j)$ and all the $a_j$'s belong to a small subset $S\in\Ac$. We estimate the integrand as follows:
\beq
\lefteqn{|\tau(e^{-s_0t^2D^2}\rho_0 e^{-s_1t^2D^2}[D,\rho_1] \ldots D e^{-s_it^2D^2}\ldots [D,\rho_n] e^{-s_{n+2}t^2D^2})|(g)}\non\\
&\leq& \int_{G^{ n}} dh_0\ldots dh_{n-1} \, \| e^{-s_0t^2D^2}\rho_0(h_0) e^{-s_1t^2D^2}[D,\rho_1] \ldots D e^{-s_it^2D^2}\non\\
&&\qquad \qquad \qquad \qquad \ldots [D,\rho_n](gh_0^{-1}\ldots h_{n-1}^{-1}) e^{-s_{n+2}t^2D^2}\|_1\ ,\non
\eeq
and then, including a plateau function $f\in \cinfc(M)$, $0\leq f\leq 1$ such that $\rho_j=f\rho_j f$ for any $j$, one has
\beq
\lefteqn{ \| e^{-s_0t^2D^2}\rho_0(h_0) e^{-s_1t^2D^2}\ldots  D e^{-s_it^2D^2} \ldots [D,\rho_n](gh_0^{-1}\ldots h_{n-1}^{-1}) e^{-s_{n+2}t^2D^2}\|_1}\non\\
&&\qquad \leq  \|\rho_0(h_0)\|_{\infty} \|f e^{-s_1t^2D^2}f\|_{s_1^{-1}}\ldots  \|f  D e^{-(s_{i-1}+s_i)t^2D^2}f\|_{(s_{i-1}+s_i)^{-1}} \ldots\non\\
&& \qquad\qquad  \ldots \|[D,\rho_n](gh_0^{-1}\ldots  h_{n-1}^{-1})\|_{\infty} \|f e^{-(s_{n+2}+s_0)t^2D^2}f\|_{(s_{n+2}+s_0)^{-1}} \ .\non
\eeq
For any $j\neq i$ and $i-1$, the estimates of \cite{FK} and the fact that the operator $e^{-t^2D^2/2}$ has spectrum $\leq 1$ show
$$\|f e^{-s_jt^2D^2}f\|_{s_j^{-1}}\leq (\Tr\, f e^{-t^2D^2}f)^{s_j} \leq (\Tr\, f e^{-t^2D^2/2}f)^{s_j}$$ 
whereas for the factor containing $D$ we write, with $s=s_{i-1}+s_i$:
\beq
\|f D e^{-st^2D^2}f\|_{s^{-1}} &=& \frac{1}{t\sqrt{s}}\, \|f t\sqrt{s}D e^{-st^2D^2}f\|_{s^{-1}}\non\\
&\leq& \frac{1}{t\sqrt{s}}\, \big(\Tr\, f (t\sqrt{s}D e^{-st^2D^2})^{s^{-1}}f\big)^s\ .\non
\eeq
Then write $ f (t\sqrt{s}D e^{-st^2D^2})^{s^{-1}}f= f (t\sqrt{s}D e^{-\frac{st^2D^2}{2}})^{s^{-1}}e^{-\frac{t^2D^2}{2}}f$ and observe that the function $x\mapsto x e^{-x^2/2}$ attains its maximum for $x=1$ with value $1/\sqrt{e} < 1$. Therefore the spectrum of $(t\sqrt{s}D e^{-st^2D^2/2})^{s^{-1}}$ is contained in $[0,1]$ and we have 
$$
\Tr\,( f (t\sqrt{s}D e^{-st^2D^2})^{s^{-1}}f) \leq \Tr\,( fe^{-t^2D^2/2}f)
$$
so that finally
$$
\|f D e^{-st^2D^2}f\|_{s^{-1}} \leq \frac{1}{t\sqrt{s}}\, (\Tr\,( fe^{-t^2D^2/2}f))^s\ .
$$
Then, by proceeding as in the proof of Proposition \ref{pchi} we introduce the functions $b_0=\| \rho_0\|_{\infty}$ and $b_i=\|[D,\rho_i]\|_{\infty}$ included in a small subset $T_S\subset C_c(G)$ depending only on $S$, so that
\beq
\lefteqn{|\tau(e^{-s_0t^2D^2}\rho_0 e^{-s_1t^2D^2} \ldots D e^{-s_it^2D^2}\ldots [D,\rho_n] e^{-s_{n+2}t^2D^2})|(g)}\non\\
&&\qquad\qquad\qquad \leq \ \frac{\Tr(fe^{-t^2D^2/2}f)}{t\sqrt{s_{i-1}+s_i}}\ (b_0\ldots b_n)(g) \non
\eeq
at any point $g\in G$. Hence taking the integral over $\Delta_{n+2}$ yields a factor
$$
\int_{\Delta_{n+2}}ds\, \frac{1}{\sqrt{s_{i-1}+s_i}}=\frac{1}{n!}\int_0^1 du_0 \int_{u_0}^1du_1\, \frac{(1-u_1)^n}{\sqrt{u_1}} \leq \frac{2}{(n+1)!}\ .
$$
The same bound holds for all the $(n+1)(n+2)$ terms of the sum (\ref{eoj1}), so that finally the norm of $cs_0^n(a_0da_1\ldots da_n)$ in the Banach algebra $\Bc$ may be estimated as
$$
\|cs_0^n(a_0da_1\ldots da_n)\| \leq  \Tr(fe^{-t^2D^2/2}f)\,2 \frac{(n+2)}{n!}\,t^{n-1} \la_S^{n+1}\ ,
$$
with a parameter $\la_S=\|T_S\|$ depending only on $S$. The end of the argument exactly follows the proof of Proposition \ref{pchi}. \cqfd\\

The cocycle $\chi(\Ec,\rho,tD)$ cannot represent the bivariant Chern character of the unbounded bimodule $(\Ec,\rho,D)$, since the $X$-complex $X(\Bc)$ does not calculate the entire cyclic homology of $\Bc$ in general. A satisfactory bivariant Chern character is built in \cite{P1} by a slight modification of the construction above. The point is to replace $\Ac$ and $\Bc$ by their analytic tensor algebras $\Tc\Ac$ and $\Tc\Bc$, and the bimodule $(\Ec,\rho,D)$ must be lifted to a suitable bimodule $(\Oman^+\Ec,\rho_*,D)$ over tensor algebras. This yields a bounded chain map
\be
\chi(\Oman^+\Ec,\rho_*,tD):\Ome\Tc\Ac\to X(\Tc\Bc)
\ee
for any $t>0$. By the generalized Goodwillie theorem \cite{Me, P1}, the complex $\Ome\Tc\Ac$ calculates the entire cyclic homology of $\Ac$. Since $X(\Tc\Bc)$ calculates the entire cyclic homology of $\Bc$, the cohomology class of $\chi(\Oman^+\Ec,\rho_*,tD)$, independent of $t$, is an element of the bivariant entire cyclic cohomology $HE_*(\Ac,\Bc)$, whose degree coincides with the parity of $(\Ec,\rho,D)$. This is by definition the bivariant Chern character.\\
Let us explain, in our context, how to obtain the lifted bimodule $(\Oman^+\Ec,\rho_*,D)$. With $\Bct$ denoting the unitalization of $\Bc$, we already introduced the $\Ac$-$\Bct$-bimodule $\Ect=\Hc\hotimes\Bct$. On the other hand, the analytic completion of the DG algebra of noncommutative differential forms $\Oman\Bc$ (see section \ref{sent}) may be considered as a left $\Bct$-module and a right $\Oman\Bc$-module, using the obvious multiplications. Therefore, the space
\be
\Oman\Ec= \Ect\hotimes_{\Bct}\Oman\Bc
\ee
is an $\Ac$-$\Oman\Bc$-bimodule. $\Oman\Ec$ is isomorphic to $\Hc\hotimes\Oman\Bc$ as a bornological space, and the splitting of $\Oman\Bc$ as the direct sum of even/odd forms $\Oman^{\pm}\Bc$ leads to the subspaces $\Oman^{\pm}\Ec=\Hc\hotimes\Oman^{\pm}\Bc$. Finally the differential $d$ on $\Oman\Bc$ extends to a flat connection
\be
d:\Oman^{\pm}\Ec\to\Oman^{\mp}\Ec\ .
\ee
To make the link with analytic tensor algebras, we shall deform the module structures using a Fedosov-type product. Recall that the analytic tensor algebra $\Tc\Bc$ over $\Bc$ is isomorphic to the algebra of even-degree differential forms $(\Oman^+\Bc,\odot)$ endowed with the Fedosov product
\be
\om_1\odot\om_2=\om_1\om_2-d\om_1d\om_2\ ,\quad \forall \om_i\in\Oman^+\Bc\ .
\ee
Indeed the isomorphism is induced by the correspondence 
\be
\Oman^+\Bc\ni b_0db_1db_2\ldots db_{2k-1}db_{2k} \leftrightarrow b_0\otimes \om(b_1,b_2)\otimes\ldots\otimes \om(b_{2k-1},b_{2k})\in \Tc\Bc\ ,
\ee
with $\om(b_1,b_2)=b_1b_2-b_1\otimes b_2$. One then gets a right $\Tc\Bc$-module structure on $\Oman^+\Ec$ via the Fedosov deformation
\be
\odot: \Oman^+\Ec\times \Tc\Bc \to\Oman^+\Ec\ ,\quad \xi\odot \om=\xi\om - d\xi d\om
\ee
for any $\xi\in\Oman^+\Ec$ and $\om\in\Oman^+\Bc=\Tc\Bc$. Let $\End_{\Tc\Bc}(\Oman^+\Ec)$ be the algebra of bounded $\Tc\Bc$-module endomorphisms. We are looking for a homomorphism from the non-analytic tensor algebra over $\Ac$
\be
\rho_*:T\Ac \to \End_{\Tc\Bc}(\Oman^+\Ec)
\ee
lifting the homomomorphism $\rho:\Ac\to \End_{\Bc}(\Ec)$. First, we define a \emph{linear map} $\rho_*:\Ac\to \End_{\Tc\Bc}(\Oman^+\Ec)$ via a Fedosov deformation (recall that $\Oman\Ec$ is a left $\Ac$-module):
\be
\rho_*(a)\xi= \rho(a)\xi- d\rho(a)d\xi\ ,\quad \forall a\in\Ac\ ,\ \xi\in \Oman^+\Ec\ .
\ee
Here $d\rho(a)$ means the commutator $[d,\rho(a)]$. This linear map then extends uniquely to a homomorphism $\rho_*$ on $T\Ac$ by the universal property of the tensor algebra. However, it is not clear that this homomorphism extends to the \emph{analytic} tensor algebra $\Tc\Ac$. One of the aims of the proposition below is to show that the construction of the chain map $\chi(\Oman^+\Ec,\rho_*,tD)$ is indeed compatible with analytic completions. Hence, what we obtained so far is a $T\Ac$-$\Tc\Bc$-bimodule $\Oman^+\Ec$. Finally, the action of the unbounded operator $D$ on $\Hc$ yields an unbounded operator on $\Oman^+\Ec=\Hc\hotimes\Oman^+\Bc$ commuting with the right action of $\Tc\Bc$. Note that the commutator of $D$ with any element $x\in T\Ac$ is a bounded endomorphism
\be
[D,\rho_*(x)]\in \End_{\Tc\Bc}(\Oman^+\Ec)\ ,
\ee
because this property can be checked on the generators $a\in\Ac$ of $T\Ac$. Conventionally we assume that $D$ is an operator of odd degree and therefore anticommutes with the differential $d$, see \cite{P1}. We introduce now some useful subalgebras of trace-class endomorphisms of $\Oman^+\Ec$. Since $\ell^1$ acts on $\Hc$ and $\Oman^+\Bc$ acts on itself (from the left) by Fedosov multiplication, the tensor product of algebras $\ell^1\hotimes \Oman^+\Bc$ (or $\ell^1\hotimes\Tc\Bc$) clearly acts on $\Oman^+\Ec=\Hc\hotimes \Oman^+\Bc$ by bounded endomorphisms. Also, remark that in any degree $2n$ the subspace of $2n$-forms $\Om^{2n}\Bc\subset \Oman^+\Bc$ is isomorphic to
\be
\Om^{2n}\Bc=\Bc^{\hotimes 2n}\oplus \Bc^{\hotimes (2n+1)}=L^1(G^{2n})\oplus L^1(G^{2n+1})\ ,
\ee
the $L^1$-spaces being taken with respect to the admissible measure $d\nu$ on $G$. Therefore we define the subspaces of continuous $2n$-forms with compact support \be
\Om_c^{2n}(G)=C_c(G^{2n})\oplus C_c(G^{2n+1})\subset \Om^{2n}\Bc\ ,
\ee
and the direct sum
\be
\Om_c^+(G)=\bigoplus_{n\geq 0}\Om^{2n}_c(G)=\underbrace{C_c(G)}_{\Om^0_c(G)}\oplus \underbrace{C_c(G\times G)\oplus C_c(G\times G\times G)}_{\Om^2_c(G)} \oplus \ldots
\ee
is a subspace of $\Oman^+\Bc$. Moreover, the Fedosov product on $\Oman^+\Bc$ restricts to a product
\be
\odot: \Om_c^+(G)\times \Om_c^+(G) \to \Om_c^+(G)\ .
\ee
In the same way, we note that 
\be
\ell^1\hotimes \Om^{2n}\Bc= L^1(G^{2n};\ell^1)\oplus L^1(G^{2n+1};\ell^1)\ ,
\ee
hence the subspaces of continuous $2n$-forms with values in $\ell^1$
\be
\Om_c^{2n}(G;\ell^1)=C_c(G^{2n};\ell^1)\oplus C_c(G^{2n+1};\ell^1)\subset \ell^1\hotimes\Om^{2n}\Bc
\ee
yield a subalgebra $\Om_c^+(G;\ell^1)=\bigoplus_{n\geq 0}\Om^{2n}_c(G;\ell^1)$ of $\ell^1\hotimes\Oman^+\Bc$. We have thus obtained a sequence of inclusions 
\be
\Om_c^+(G;\ell^1)\hookrightarrow \ell^1\hotimes\Oman^+\Bc\hookrightarrow \End_{\Tc\Bc}(\Oman^+\Ec)\ .
\ee
As a crucial result, the heat operator $\exp(-sD^2)\in \End_{\Tc\Bc}(\Oman^+\Ec)$ acts as a regulator for $s>0$, in the sense that the following endomorphism lies in the smallest possible algebra:
\be
e^{-sD^2}\rho_*(x)\in \Om_c^+(G;\ell^1)\quad \forall x\in T\Ac\ ,\ s>0\ .
\ee
The same is true for $e^{-sD^2}[D,\rho_*(x)]$. Indeed, this can be seen by an explicit computation when $x$ is, say, a $2n$-form $a_0da_1\ldots da_{2n}\in \Om^{2n}\Ac$ (we use the identification $T\Ac=\Om^+\Ac$). The product $e^{-sD^2}\rho_*(x)\in \Om^{2n}_c(G;\ell^1)$ is then a function in $C_c(G^{2n+1};\ell^1)$, whose evaluation at a point $(g_0,\ldots,g_{2n})$ is the trace-class operator
\be
e^{-sD^2}\rho_*(a_0da_1\ldots da_{2n})(g_0,\ldots,g_{2n})=e^{-sD^2}a_0(g_0)r(g_0)\ldots a_{2n}(g_{2n})r(g_{2n})\ .
\ee
The situation with the triple $(\Oman^+\Ec,\rho_*,D)$ is therefore completely analogous to the previous situation with $(\Ec,\rho,D)$. The chain map $\chi(\Oman^+\Ec,\rho_*,tD)$ is given by formulas (\ref{joe1},\ref{joe2}), replacing everywhere $\Ec$, $\Ac$, $\Bc$, $C_c(G)$, $\rho$ respectively by $\Oman^+\Ec$, $T\Ac$, $\Tc\Bc$, $\Om_c^+(G)$, $\rho_*$. Hence, for any $t>0$ and any integer $n\in\nn$ whose parity equals the degree of the bimodule, the component $\chi^n_0: \Om^nT\Ac\to \Tc\Bc$ evaluated on a $n$-form $x_0\dd x_1\ldots \dd x_n\in \Om^nT\Ac$ (we use the script $\dd$ for differential forms over the algebra $T\Ac$, in order to avoid confusion with differential forms over $\Ac$) is
\beq
\lefteqn{\chi^n_0(x_0\dd x_1\ldots \dd x_n) = (-t)^n\sum_{i=0}^n (-)^{i(n-i)}\times } \label{joe3}\\
&&\tau \langle [D,\rho_*^{i+1}],\ldots ,[D,\rho_*^n] ,\rho_*^0,[D,\rho_*^1],\ldots, [D,\rho_*^i]\rangle_t\ ,\non
\eeq
with the abbreviation $\rho_*^i=\rho_*(x_i)$. The right-hand-side is actually an element of the algebra $\Om_c^+(G)\subset\Tc\Bc$. The other components $\chi^n_1: \Om^{n+1}T\Ac\to \Om^1\Tc\Bc_{\nat}$ are obtained along the lines leading to formula (\ref{joe2}): one first has to consider the $T\Ac$-$\Tc\Bc$-bimodule $\Om^1\Oman^+\Ec$ and work with suitable algebras of endomorphisms. Details are left to the reader. We end up with the formula
\beq
\lefteqn{\chi^{n+1}_1(x_0\dd x_1\ldots \dd x_{n+1}) = (-t)^n\sum_{i=1}^{n+1} (-)^{i(n-i+1)}\times} \label{joe4}\\
&&\nat \tau(\langle [D,\rho_*^{i+1}],\ldots , [D,\rho_*^{n+1}],\rho_*^0,[D,\rho_*^1],\ldots, [D,\rho_*^{i-1}]\rangle_t \dd \rho_*^i) \ .\non
\eeq
Here the right-hand side must be interpreted as an element of the subspace of continuous forms of odd degree $\Om^-_c(G)\subset \Oman^-\Bc\cong \Om^1\Tc\Bc_{\nat}$. One thus gets the components of a chain map $\chi(\Oman^+\Ec,\rho_*,tD):\Om T\Ac\to X(\Tc\Bc)$, and similarly for the Chern-Simons transgressions
\be
cs^n(\Oman^+\Ec,\rho_*,tD): \Om^nT\Ac\to X(\Tc\Bc)\hotimes\Om^1\rr_+^*\ .
\ee
The delicate point is that these maps extend to the space of entire forms over the analytic tensor algebra $\Tc\Ac$. 
\begin{proposition}\label{pana}
Let $(\Ec,\rho,D)$ be an unbounded $\Ac$-$\Bc$-bimodule as in Definition \ref{dunb}. Then for any $t>0$, the collection of components $\chi^n(\Oman^+\Ec,\rho_*,tD)$ extends to an entire cyclic cocycle, i.e. a bounded chain map
\be
\chi(\Oman^+\Ec,\rho_*,tD):\Ome\Tc\Ac\to X(\Tc\Bc)\ ,
\ee
and similarly the Chern-Simons transgression extends to a bounded linear map $cs(\Oman^+\Ec,\rho_*,tD)$, whose coboundary is $d_t\chi(\Oman^+\Ec,\rho_*,tD)$. The bivariant Chern character $\ch(\Ec,\rho,D)$ is the class of the cocycle $\chi(\Oman^+\Ec,\rho_*,tD)$ in the bivariant entire cyclic cohomology $HE_*(\Ac,\Bc)$, independent of $t>0$. The degree of the Chern character coincides with the parity of the unbounded bimodule.
\end{proposition}
{\it Proof:}  We shall only show that the components landing to the even degree subspace of the $X$-complex $\chi^n_0:\Om^n T\Ac\to \Tc\Bc$ extend to an entire cochain. The odd case is treated similarly, as well as the Chern-Simons transgressions.\\
Recall that the tensor algebra $T\Ac$ is isomorphic to the algebra of differential forms of even degree $\Om^+\Ac$ endowed with the Fedosov product
$$
x_1\odot x_2=x_1x_2-dx_1dx_2\ .
$$
A small subset of the analytic bornology on $\Om^+\Ac$ is contained in the convex hull of $\bigcup_{k\ge 0} \St (dS)^{2k}$ for some small $S\subset \Ac$, and the analytic tensor algebra $\Tc\Ac$ is the completion of $T\Ac$ with respect to this bornology. Also, $\Ome\Tc\Ac$ is the completion of $\Om \Tc\Ac$ with respect to the entire bornology, generated by the subsets $\bigcup_{n\ge 0} [n/2]!\,\Ut (\dd U)^n$, for all small $U\subset\Tc\Ac$. As already mentioned, $\dd$ denotes the differential of forms over $\Tc\Ac$, to avoid confusion with the differential $d$ of forms over $\Ac$. It is possible to obtain $\Ome\Tc\Ac$ directly from one completion of the space $\Om T\Ac$, by taking the bornology generated by the subsets \cite{P1}
$$
\bigcup_{n\ge 0} [n/2]!\,\Ut (\dd U)^n\ ,\ \mbox{with}\  U=\bigcup_{k\ge 0} \St (dS)^{2k}\ ,
$$
for all small $S\subset\Ac$. So we fix $S$ and consider $n+1$ elements $x_0,\ldots x_n$ contained in the small subset $\bigcup_{k\ge 0} \St (dS)^{2k}\subset \Tc\Ac$, so that we can write without loss of generality
$$
x_i=a^i_0da^i_1\ldots da^i_{2k_i}\in \Om^{2k_i}\Ac\ ,\qquad \forall i\in\{0,\ldots, n\}
$$
for some elements $a^i_j\in S$. Therefore, we must show that the element 
$$\chi^n_0(x_0\dd x_1\ldots\dd x_n) \in \Om_c^+(G)$$ 
lies in a small subset of $\Oman^+\Bc=\Tc\Bc$, depending only on $S$. We will use the following notation. If $b_0,\ldots, b_{2k}$ are $2k+1$ elements in $C_c(G)=\Om_c^0(G)$, we denote by $b_0db_1\ldots db_{2k}$ the corresponding differential form in $\Om^+_c(G)$: it is a function in $C_c(G^{2k+1})$ whose evaluation at a point $(g_0,\ldots,g_{2k})$ reads
$$
b_0db_1\ldots db_{2k}(g_0,\ldots,g_{2k})=b_0(g_0)b_1(g_1)\ldots b_{2k}(g_{2k})\ .
$$
Also, an exact form $db_1\ldots db_{2k}$ is an element of $C_c(G^{2k})$:
$$
db_1\ldots db_{2k}(g_1,\ldots,g_{2k})=b_1(g_1)\ldots b_{2k}(g_{2k})\ .
$$
Then the product on $\Om^+_c(G)$ induced by the inclusion $\Om^+_c(G)\subset \Oman^+\Bc$ is of course of Fedosov type. The differential form $\chi^n_0(x_0\dd x_1\ldots\dd x_n)$ is given by the sum (\ref{joe3}). Let us isolate for example the last term, 
$$
I=t^n\int_{\Delta_{n+1}}ds\, \tau(e^{-s_0t^2D^2}\rho^0_* e^{-s_1t^2D^2}\odot [D,\rho^1_*] e^{-s_2t^2D^2}\ldots\odot [D,\rho^n_*] e^{-s_{n+1}t^2D^2})\ ,
$$
with the notation $\rho^i_*:=\rho_*(x_i)=\rho(a_0^i)d\rho(a_1^i)\ldots d\rho(a_{2k_i}^i)$. The products between the $\rho_*^i$'s and $[D,\rho_*^i]$'s are again of Fedosov type. To become familiar with these notations, let us calculate a little example where $n=1$ and
$$
x_0=a_0da_1da_2\in \Om^2\Ac\ ,\qquad x_1=a_3\in\Om^0\Ac\ .
$$
We abbreviate $\rho(a_i)=\rho_i$. Thus one wants to estimate the supremum norm of 
$$
I=t \int_{\Delta_2}ds\, \tau(e^{-s_0t^2D^2}(\rho_0d\rho_1d\rho_2) e^{-s_1t^2D^2}\odot [D,\rho_3] e^{-s_2t^2D^2})\ .
$$
It is an element of $\Om^2_c(G)\oplus\Om^4_c(G)$. The operator under the trace is computed as follows:
\beq
\lefteqn{e^{-s_0t^2D^2}(\rho_0d\rho_1d\rho_2) e^{-s_1t^2D^2}\odot [D,\rho_3] e^{-s_2t^2D^2}}\non\\ 
&=& e^{-s_0t^2D^2}\rho_0d\rho_1d(\rho_2 e^{-s_1t^2D^2} [D,\rho_3]) e^{-s_2t^2D^2}\non\\
&& -e^{-s_0t^2D^2}\rho_0d(\rho_1\rho_2) e^{-s_1t^2D^2} d[D,\rho_3] e^{-s_2t^2D^2}\non\\
&& + e^{-s_0t^2D^2}\rho_0\rho_1d\rho_2 e^{-s_1t^2D^2} d[D,\rho_3] e^{-s_2t^2D^2}\non\\
&& -e^{-s_0t^2D^2}d\rho_0d\rho_1d\rho_2 e^{-s_1t^2D^2} d[D,\rho_3] e^{-s_2t^2D^2}\ .\non
\eeq
The first three terms of the r.h.s. lie in the summand $C_c(G^{ 3};\ell^1)$ of the two-form space $\Om^2_c(G;\ell^1)$, whereas the last term is in the summand $C_c(G^{ 4};\ell^1)$ of the four-form space $\Om^4_c(G;\ell^1)$. Therefore, we can evaluate the integrand of $I$ either on three points:
\beq
\lefteqn{|\tau(e^{-s_0t^2D^2}(\rho_0d\rho_1d\rho_2) e^{-s_1t^2D^2}\odot [D,\rho_3] e^{-s_2t^2D^2})|(g_0,g_1,g_2)}\non\\
&\leq& \int_G dh\,  \|e^{-s_0t^2D^2}\rho_0(g_0)\rho_1(g_1)\rho_2(h) e^{-s_1t^2D^2} [D,\rho_3](g_2h^{-1}) e^{-s_2t^2D^2}\non\\
&&\qquad  - e^{-s_0t^2D^2}\rho_0(g_0)\rho_1(h)\rho_2(g_1h^{-1}) e^{-s_1t^2D^2} [D,\rho_3](g_2) e^{-s_2t^2D^2}\non\\
&&\qquad  + e^{-s_0t^2D^2}\rho_0(h)\rho_1(g_0h^{-1})\rho_2(g_1) e^{-s_1t^2D^2} [D,\rho_3](g_2) e^{-s_2t^2D^2}\|_1\non\\
&\leq& \Tr(f e^{-t^2D^2}f)\, \int_G dh\,(b_0(g_0)b_1(g_1)b_2(h)b_3(g_2h^{-1}) \non\\
&& \qquad\qquad\qquad + b_0(g_0)b_1(h)b_2(g_1h^{-1})b_3(g_2)\non\\
&& \qquad\qquad\qquad + b_0(h)b_1(g_0h^{-1})b_2(g_1)b_3(g_2))\non\\
&\leq& \Tr(f e^{-t^2D^2}f)\, (b_0db_1d(b_2b_3) + b_0d(b_1b_2)db_3 + b_0b_1db_2db_3)(g_0,g_1,g_2)\ ,\non
\eeq
or on four points:
\beq
\lefteqn{|\tau(e^{-s_0t^2D^2}(\rho_0d\rho_1d\rho_2) e^{-s_1t^2D^2}\odot [D,\rho_3] e^{-s_2t^2D^2})|(g_0,g_1,g_2,g_3)}\non\\
&\leq& \|e^{-s_0t^2D^2}\rho_0(g_0)\rho_1(g_1)\rho_2(g_2) e^{-s_1t^2D^2} [D,\rho_3](g_3) e^{-s_2t^2D^2}\|_1\non\\
&\leq& \Tr(f e^{-t^2D^2}f)\, b_0(g_0)b_1(g_1)b_2(g_2)b_3(g_3)\non\\
&\leq& \Tr(f e^{-t^2D^2}f)\, (db_0db_1db_2db_3)(g_0,g_1,g_2,g_3)\ .\non
\eeq
We have adopted the notations and method of the proof of Proposition \ref{pchi}: $b_i\in C_c(G)$ is the positive function $\|\rho_i\|_{\infty}$ for $i=0,1,2$, $b_3=\|[D,\rho_3]\|_{\infty}$, and $f\in \cinfc(M)$, $0\leq f\leq 1$ is a plateau function verifying $\rho_i=f\rho_i f$. All the $b_i$'s belong to a small subset $T_S$ in the bornology of the LF-space $C_c(G)$, depending only on $S$. We estimate the component of $I$ in $\Om^2_c(G)$ as follows. Consider it as an element of the subspace $\Bc\hotimes\Bc\hotimes\Bc\subset \Om^2\Bc$. Denote by $\|\cdot\|^{\hotimes 3}$ the projective norm on the Banach space $\Bc^{\hotimes 3}$ coming from the $L^1$-norm $\|\cdot \|$ of $\Bc$. If $p_{2m}:\Om^+_c(G)\to\Om^{2m}_c(G)$ denotes the projection onto $2m$-forms, we have
$$
\|p_2 I\|^{\hotimes 3}\leq \frac{t}{2}\Tr(f e^{-t^2D^2}f)\|b_0db_1d(b_2b_3) + b_0d(b_1b_2)db_3 + b_0b_1db_2db_3\|^{\hotimes 3}\ ,
$$
where the factor $1/2$ comes from the integration over $\Delta_2$. Each of the three terms in the r.h.s can be controlled, for instance
$$
\|b_0db_1d(b_2b_3)\|^{\hotimes 3}=\|b_0\|\, \|b_1\|\, \|b_2b_3\|\leq \la_S^4
$$
with the parameter $\la_S=\| T_S\|$ depending only on $S$. Consequently,
$$
\|p_2 I\|^{\hotimes 3}\leq \frac{t}{2}\Tr(fe^{-t^2D^2}f)3 \la_S^4\ .
$$
In the same way, the component of $I$ in $\Om^4\Bc$ may be estimated by the projective norm $\|\cdot\|^{\hotimes 4}$ on $\Bc^{\hotimes 4}$. One has $\|db_0db_1db_2db_3\|^{\hotimes 4}=\|b_0\|\,\|b_1\|\,\|b_2\|\,\|b_3\|$, so that
$$
\|p_4I\|^{\hotimes 4} \leq \frac{t}{2}\Tr(fe^{-t^2D^2}f) \la_S^{4}\ .
$$
Remark that we could control not only the Banach norm of $I$ as a $L^1$ function over $G^3$ and $G^4$, but also its supremum norm as a continuous function with compact support. Although we do not need it for the moment, this estimate will be useful for proving Corollary \ref{cloc}. Let us look for example at the supremum norm of the sum $b_0db_1d(b_2b_3) + b_0d(b_1b_2)db_3 + b_0b_1db_2db_3 \in \Om^2_c(G)$, as a function over $G^3$. We know that the $b_i$'s belong to a small subset $T_S\subset C_c(G)$, hence there is a compact subset $K_S\subset G$ and a number $\zeta_S$ such that $b_i$ has support contained in $K_S$ and $\sup_{g\in K_S}|b_i(g)|\leq \zeta_S$. Moreover, the norm of a product $b_ib_j$ may be controlled as follows:
$$
\sup_{g\in G}|(b_ib_j)(g)|\leq \sup_{g\in G}\int_Gdh\, |b_i(h)b_j(gh^{-1})|\leq \zeta_S^2\int_{K_S}dh=\zeta_S^2|K_S|\ ,
$$
where $|K_S|$ denotes the volume of $K_S$ with respect to the Haar measure $dh$. Therefore, we deduce the following bound for the supremum norm of $p_2 I$ over $G^3$:
$$
\sup_{G^3}|p_2I|\leq \frac{t}{2}\Tr(fe^{-t^2D^2}f)3 |K_S|\zeta_S^4\ .
$$
Now turn back to the general case. The last term of the sum (\ref{joe3}),
$$
I=t^n\int_{\Delta_{n+1}}ds\, \tau(e^{-s_0t^2D^2}\rho^0_* e^{-s_1t^2D^2}\odot [D,\rho^1_*] e^{-s_2t^2D^2}\ldots\odot [D,\rho^n_*] e^{-s_{n+1}t^2D^2})\ ,
$$
with $\rho^i_*=\rho(a_0^i)d\rho(a_1^i)\ldots d\rho(a_{2k_i}^i)$, is a differential form in $\Om^+_c(G)$ of degree $\geq 2k$, for $k=\sum_{i=0}^nk_i$. First, remark that the commutator $[D,\rho^i_*]$ is a differential form of degree $2k_i$ (recall $D$ is odd):
\beq
\lefteqn{[D,\rho^i_*]=[D,\rho(a_0^i)]d\rho(a_1^i)\ldots d\rho(a_{2k_i}^i)}\non\\
&&\qquad\qquad + \sum_{j=1}^{2k_i}(-)^j \rho(a_0^i)d\rho(a_1^i)\ldots d[D,\rho(a^i_j)]\ldots  d\rho(a_{2k_i}^i)\ .\non
\eeq
It is a sum of $2k_i+1$ terms, so that $I$ may actually be written as a sum of $(2k_1+1)\ldots(2k_n+1)$ terms. Let us concentrate on one of these. If we neglect the presence of the heat operators for a while, we are led to calculate the Fedosov product
\beq
\lefteqn{\rho_*(a^0_0da^0_1\ldots da^0_{2k_0})\odot \rho_*(a^1_0da^1_1\ldots da^1_{2k_1})\odot \ldots \odot \rho_*(a^n_0da^n_1\ldots da^n_{2k_n})}\non\\
&&\qquad\qquad = \rho_*(a^0_0da^0_1\ldots da^0_{2k_0} a^1_0da^1_1\ldots da^1_{2k_1} \ldots  a^n_0da^n_1\ldots da^n_{2k_n})\non\\
&& \qquad\qquad\qquad +\ \mbox{terms of degree $>2k$,}\non
\eeq
where some $\rho(a^i_j)$'s may be replaced by the commutator $[D,\rho(a^i_j)]$ as well. We find an upper bound for this quantity via the introduction of the functions $b^i_j=\|\rho(a^i_j)\|_{\infty}$ or $\|[D,\rho(a^i_j)]\|_{\infty}$ over $G$, as in the above little example. All the $b^i_j$'s are contained in a small subset $T_S\subset C_c(G)$ depending only on $S$. It is not difficult to see, using the Leibniz rule, that the term of degree $2k$ in the r.h.s. of the above equality is bounded by a sum of $(2k_0+1)\ldots(2k_{n-1}+1)$ terms (in the example, it corresponds to $3$ terms $b_0db_1d(b_2b_3) + b_0d(b_1b_2)db_3 + b_0b_1db_2db_3$). Therefore, by proceeding as before, the component of $I$ belonging to the Banach subspace $\Bc^{\hotimes (2k+1)}\subset \Om^{2k}\Bc$ may be estimated by the projective norm $\|\cdot\|^{\hotimes (2k+1)}$ as follows:
$$
\|p_{2k}I\|^{\hotimes(2k+1)} \leq  \frac{t^n}{(n+1)!} \Tr(fe^{-t^2D^2}f)\, (2k_0+1)\prod_{i=1}^{n-1}(2k_i+1)^2(2k_n+1) \la_S^{2k+n+1}\ .
$$
The same estimate holds for all the $(n+1)$ terms of the sum (\ref{joe3}), so that finally the component of $\chi^n_0(x_0\dd x_1\ldots\dd x_n)$ in $\Om^{2k}\Bc$ verifies the inequality
\beq
\lefteqn{\|p_{2k}\chi^n_0(x_0\dd x_1\ldots\dd x_n)\|^{\hotimes (2k+1)}} \non\\
&\leq& \frac{t^n}{n!} \Tr(fe^{-t^2D^2}f)\, (2k_0+1)\prod_{i=1}^{n-1}(2k_i+1)^2(2k_n+1) \la_S^{2k+n+1} \ .\non
\eeq
Next, the components of $I$ of degree $2m>2k$ arise from the Fedosov products between the $\rho^i_*$'s, and the number of terms of a given degree is strictly less than $(2k_0+1)\ldots(2k_{n-1}+1)$ (in the little example appeared the only term $db_0db_1db_2db_3$). Thus when $2m>2k$ one has the same bound for the projective norms $\|\cdot\|^{\hotimes 2m}$ or $\|\cdot\|^{\hotimes(2m+1)}$ on $\Om^{2m}\Bc$:
\beq
\lefteqn{\|p_{2m}\chi^n_0(x_0\dd x_1\ldots\dd x_n)\|^{\hotimes(2m,2m+1)}} \non\\
&\leq& \frac{t^n}{n!} \Tr(fe^{-t^2D^2}f)\, (2k_0+1)\prod_{i=1}^{n-1}(2k_i+1)^2(2k_n+1) \la_S^{2k+n+1}\ .\non
\eeq
Suppose now $\la_S\geq 1$ (if it is less than 1, the following discussion is trivial), so that $\la_S^{2m}\geq\la_S^{2k}$ whenever $m\geq k$. Since the number $(2k_0+1)(2k_1+1)^2\ldots(2k_{n-1}+1)^2(2k_n+1)$ is always bounded by $16^k$, the component of the form $\chi^n_0(x_0\dd x_1\ldots\dd x_n)$ in {\it any} of the subspaces $\Om^{2m}\Bc$, $m\in\nn$, can be estimated as
\beq
\|p_{2m}\chi^n_0(x_0\dd x_1\ldots\dd x_n)\|^{\hotimes(2m,2m+1)} &\leq& \frac{t^n}{n!} \Tr(fe^{-t^2D^2}f)\, 16^k  \la_S^{2k+n+1}\non\\
& \leq&  \Tr(fe^{-t^2D^2}f)\,\frac{(t\la_S)^n}{n!}\, (4\la_S)^{2m+1} \ .\non
\eeq
Finally, remark that we may also estimate as in the little example the supremum norm of $p_{2m}\chi^n_0(x_0\dd x_1\ldots\dd x_n)\in \Om^{2m}_c(G)$ as a continuous function with compact support on $G^{2m}\cup G^{2m+1}$. Indeed, for any elements $b_0,\ldots,b_k$ contained in the small subset $T_S\subset C_c(G)$, one has by induction on $k$
$$
\sup_{g\in G}|(b_0\ldots b_k)(g)|\leq \zeta_S^{k+1}|K_S|^k\ ,
$$
where the numbers $\zeta_S$ and $|K_S|$ depend only on $S$. Therefore
$$
\sup_{G^{2m}\cup G^{2m+1}}|p_{2m}\chi^{n}_0(x_0\dd x_1\ldots\dd x_{n})| \leq  \Tr(fe^{-t^2D^2}f)\, \frac{(t\zeta_S\,|K_S|)^{n}}{n!}\, (4\zeta_S)^{2m+1}\ .
$$
We are ready to conclude. The function $n\mapsto (t\la_S)^n[n/2]!/n!$ being bounded, the element 
$$
[n/2]!\chi^n_0(x_0\dd x_1\ldots\dd x_n)\in \Oman^+\Bc
$$
is a finite sum of differential forms whose norm in each subspace $\Om^{2m}\Bc$ is bounded by a constant times $\la_S^{2m+1}$, uniformly in $n$. By the very definition of the analytic bornology on $\Oman^+\Bc$ (section \ref{sent}), this shows that the collection of components $\chi_0$ sends the small subset of $\Om T\Ac$
$$
\bigcup_{n\ge 0} [n/2]!\,\Ut (\dd U)^n\ ,\ \mbox{with}\ U=\bigcup_{k\ge 0} \St (dS)^{2k}\ ,
$$
to a small subset of the analytic tensor algebra $\Oman^+\Bc\cong \Tc\Bc$. Hence $\chi_0$ extends to a bounded map $\Ome\Tc\Ac\to\Tc\Bc$. One proceeds in the same manner with $\chi_1:\Ome\Tc\Ac\to \Om^1\Tc\Bc_{\nat}$ and the Chern-Simons transgressions.\cqfd\\

The cocycle $\chi(\Oman^+\Ec,\rho_*,tD)$ is a lift of the cocycle $\chi(\Ec,\rho,tD)$ in the sense that the diagram
\be
\begin{CD}
\Ome\Tc\Ac @>{\chi(\Oman^+\Ec,\rho_*,tD)}>> X(\Tc\Bc) \\
@V{\sim}VV  @VVV \\
\Ome\Ac @>{\chi(\Ec,\rho,tD)}>> X(\Bc)\\
\end{CD}
\ee
commutes. The vertical arrows are induced by the (bounded) multiplication homomorphisms $\Tc\Ac\to \Ac$ and $\Tc\Bc\to\Bc$. The left vertical arrow is an homotopy equivalence by virtue of the Goodwillie theorem. Since the entire cyclic homology of $\Ac$ is also computed by the $X$-complex $X(\Tc\Ac)$, one can construct an explicit chain map $\gamma: X(\Tc\Ac)\to \Ome\Tc\Ac$ realizing this equivalence \cite{P1}. The bivariant Chern character is thus represented by the composition
\be
\ch(\Ec,\rho,D)=\chi(\Oman^+\Ec,\rho_*,tD)\circ \gamma: X(\Tc\Ac)\to X(\Tc\Bc)\ .
\ee

\section{Index theorem}\label{sind}

In this section we will state the main result of this paper, namely, an index theorem in entire cyclic cohomology. Given a locally compact group $G$ acting properly on a manifold $M$, and $\Bc$ any admissible completion of the convolution algebra $C_c(G)$, our aim is to compute the composition of the assembly map with the Chern character in entire cyclic homology
\be 
\ch\circ \mu: K_*^G(M) \to K_*(\Bc)\to HE_*(\Bc)\ .
\ee
Actually, for practical purposes we are only interested in the $K$-homology classes $[D]\in K_*^G(M)$ represented by Dirac-type operators $D$ on $M$. The computation of $\ch\circ \mu (D)$ will be related to the bivariant Chern character introduced in section \ref{sbiv}. As a corollary, we find a localization formula generalizing the results of Atiyah-Singer \cite{AS}, Connes-Moscovici \cite{CM90}, or Mishchenko-Fomenko \cite{MF}.\\

As in the previous sections, $M$ is a complete $G$-compact Riemannian manifold without boundary, on which $G$ acts properly by isometries. Let $\Ac$ be the crossed-product algebra $\cinfc(M)\rtimes G$ (section \ref{sbiv}). It is an LF-algebra, and its topological $K$-theory group $K_0(\Ac)$ is provided with a canonical element $[e]$. Indeed, since $M$ is proper and $G$-compact, we can find a cut-off function $c\in\cinfc(M)$. This is a smooth, compactly supported and non-negative function $c$ such that
\be
\int_G c(gx)^2=1\qquad \forall x\in M\ .
\ee
Then the following formula defines an idempotent $e$ in the crossed product algebra $\Ac$:
\be
e(g,x)=c(x)c(gx)\quad \forall g\in G\ ,\ x\in M\ , \qquad e^2=e\ .
\ee
The $K$-theory class of $e$ does not depend on the choice of the cut-off function, two such functions being related via a smooth homotopy. The Chern character $\ch(e)$ in the entire cyclic homology of even degree $HE_0(\Ac)$ is represented by a cycle in the even part of the $X$-complex of the analytic tensor algebra $\Tc\Ac$. By \cite{CQ1,Me} (see also section \ref{sent}), it is given by the following idempotent $\eh \in \Tc\Ac$:
\be
\eh=e+ \sum_{k\ge 1} \frac{(2k)!}{(k!)^2} (e-\frac{1}{2})\otimes (e-e\otimes e)^{\otimes k}\ .\label{lif}
\ee
Now let $D$ be a $G$-invariant elliptic differential operator of order one representing a class $[D]\in K_*^G(M)$. In section \ref{sbiv} we attached to $D$ an unbounded $\Ac$-$\Bc$-bimodule $(\Ec,\rho,D)$ and obtained its Chern character as a bivariant entire cyclic cohomology class
\be 
\ch(\Ec,\rho,D)\in HE_*(\Ac,\Bc)\ .
\ee
The latter is represented by a chain map $X(\Tc\Ac)\to X(\Tc\Bc)$, composite of the homotopy equivalence $\gamma: X(\Tc\Ac)\stackrel{\sim}{\longrightarrow} \Ome\Tc\Ac$ with the character $\chi(\Oman^+,\rho_*,tD):\Ome\Tc\Ac\to X(\Tc\Bc)$ of proposition \ref{pana}. Therefore the cup-product of the class $\ch(e)\in HE_0(\Ac)$ with the bivariant Chern character of $(\Ec,\rho,D)$
\be
\ch(\Ec,\rho,D)\cdot \ch(e) \in HE_*(\Bc)
\ee
may be explicitly computed as an entire cyclic homology class of the admissible completion $\Bc$, of degree equal to the parity of the $K$-homology class $[D]$. First, the image of the idempotent $\eh \in \Tc\Ac$ under the homotopy equivalence $\gamma$ is given by the following $(b+B)$-cycle in $\Ome\Tc\Ac$ (see \cite{P1}):
\be
\gamma(\eh)= \eh+ \sum_{n\ge 1} (-)^n\frac{(2n)!}{n!} (\eh-\frac{1}{2}) (\dd\eh \dd\eh)^n\ ,
\ee
where $\dd$ is the differential of forms over $\Tc\Ac$. When the parity of $[D]$ is even, the product $\ch(\Ec,\rho,D)\cdot \ch(e)$ lies in $HE_0(\Bc)$ and is represented by the entire cycle $\chi_0(\Oman^+,\rho_*,tD)\circ \gamma(\eh) \in \Tc\Bc$. When the parity is odd, the product $\ch(\Ec,\rho,D)\cdot \ch(e)$ is an element of $HE_1(\Bc)$ represented by the one-form $\chi_1(\Oman^+,\rho_*,tD)\circ \gamma(\eh)\in \Om^1\Tc\Bc_{\nat}$. We will use these formulas to establish the local form of the equivariant index theorem. Our main result is the following:
\begin{theorem}\label{tind}
Let $D$ be a $G$-invariant elliptic differential operator of order one representing an equivariant $K$-homology class $[D]\in K_*^G(M)$. Consider the crossed-product algebra $\Ac=\cinfc(M)\rtimes G$ and its canonical $K$-theory class $[e]\in K_0(\Ac)$, and let $\Bc$ be any admissible completion of the convolution algebra $C_c(G)$. Then the Chern character of the image of $[D]$ under the analytic assembly map $\mu:K_*^G(M) \to K_*(\Bc)$ is given by the cup-product in bivariant entire cyclic cohomology
\be
\ch\circ\mu(D)= \ch(\Ec,\rho,D)\cdot \ch(e)\ \in HE_*(\Bc)\ ,
\ee
where $\ch(\Ec,\rho,D)\in HE_*(\Ac,\Bc)$ is the bivariant Chern character of the unbounded bimodule associated to $D$ by definition \ref{dunb}.
\end{theorem}
\noindent{\it Proof:} We will prove the theorem using the results of \cite{P2}. The point is that the unbounded bimodule $(\Ec,\rho,D)$ is finitely summable. As a consequence, its bivariant Chern character in entire cyclic cohomology retracts on a finite-dimensional cocycle directly related to the assembly map. We first observe that because $D$ is a differential operator of order one on the manifold $M$, the triple $(\Ec,\rho,D)$ is $p$-summable for any real number $p>\mbox{dim}\, M$, in the sense that
$$
(1+D^2)^{-1/2}\rho(a)\in \ell^p(\Hc)\hotimes\Bc\quad \forall a\in\Ac\ .
$$
As before $\Hc$ denotes the Hilbert space of square-integrable sections of the vector bundle $E\to M$ on which $D$ acts. This property allows to show that the entire chain map $\chi(\Oman^+\Ec,\rho_*,tD):\Ome\Tc\Ac\to X(\Tc\Bc)$ is cohomologous to a chain map of finite degree, i.e. vanishing on $\Om^n\Tc\Ac$ for $n$ sufficiently large. This retraction property was explained at a formal level in \cite{P2}, and we only have to check that it actually works in our present case. For convenience, let us sketch the main steps of this process. It is a bivariant generalization of the retraction presented in \cite{CM93} for finite-dimensional $K$-cycles over an algebra $\Ac$. Recall that in the proof of Proposition \ref{pana}, we established the following estimate for the component $\chi^n_0:\Om^nT\Ac\to T\Bc$ of the chain map $\chi(\Oman^+\Ec,\rho_*,tD)$. For any small subset $S\subset\Ac$ and $n+1$ elements of the (non completed) tensor algebra $x_0,\ldots ,x_n\in T\Ac=\Om^+\Ac$ given by
$$
x_i=a^i_0da^i_1\ldots da^i_{2k_i}\quad a^i_j\in S\ ,
$$
we can find a constant $\la_S$ depending only on $S$, such that the component of the non-commutative form $\chi^n_0(x_0\dd x_1\ldots\dd x_n)$ in $\Om^{2m}\Bc$ is bounded in norm by
$$
\|p_{2m}\chi^n_0(x_0\dd x_1\ldots\dd x_n)\|^{\hotimes (2m,2m+1)}\leq \Tr(f e^{-t^2D^2}f)\,\frac{(t\la_S)^n}{n!}\, (4\la_S)^{2m+1}
$$
for any $m\in\nn$. Recall that $f\in\cinfc(M)$ is a plateau function, $0\leq f\leq 1$, which also depends on $S$. The odd components $\chi^{n+1}_1:\Om^{n+1}T\Ac\to \Om^1T\Bc_{\nat}$ verify the same kind of bound, also with the factor $\Tr(fe^{-t^2D^2}f)t^n$. Moreover, proceeding as in the proof of Proposition \ref{pcs}, one shows the following bound for the component of the Chern-Simons form $cs^n_0:\Om^nT\Ac\to T\Bc\hotimes\Om^1\rr^*_+$,
$$
\|p_{2m}cs^n_0(x_0\dd x_1\ldots\dd x_n)\|^{\hotimes (2m,2m+1)}\leq \Tr(fe^{-t^2D^2/2}f)\,2 \frac{(n+2)}{n!}\,\frac{(t\la_S)^n}{t}(4\la_S)^{2m+1}
$$
and similarly for the odd components $cs^{n+1}_1:\Om^{n+1}T\Ac\to \Om^1T\Bc_{\nat}\hotimes\Om^1\rr^*_+$. Now, $p$-summability implies that the operator $(1+D^2)^{-p/2}f$ is trace-class for any $p>\dim M$. Hence
\beq
\Tr(fe^{-t^2D^2}f) &\leq& \|f(1+D^2)^{p/2}e^{-t^2D^2}(1+D^2)^{-p/2}f\|_1\non\\
&\leq& \|(1+D^2)^{p/2}e^{-t^2D^2}\|_{\infty}\|(1+D^2)^{-p/2}f\|_1\non\\
&\leq& \frac{C}{t^p} \, \|(1+D^2)^{-p/2}f\|_1\non
\eeq
for small values of $t$, where the constant $C$ is the maximum of the function $x\to (1+x)^{p/2}e^{-x}$ over $[0,\infty]$. Consequently, the components $\chi^n_0$ and $\chi^{n+1}_1$ behave like $t^{n-p}$ when $t\to 0$, whereas the Chern-Simons components $cs^n_0$ and $cs^{n+1}_1$ behave like $t^{n-p-1}$. Therefore, if $n$ is sufficiently large, one has $\lim_{t\to 0}\chi^n(\Oman^+\Ec,\rho_*,tD)=0$, and the one-form $cs^n(\Oman^+\Ec,\rho_*,tD)$ is integrable over any interval $[0,t_0]$ w.r.t. the parameter $t$. Then following exactly the proof of \cite{P2} Proposition 4.2, for any choice $t_0>0$ we can cut the tail of the cocycle $\chi(\Oman^+\Ec,\rho_*,t_0D)$ by adding an entire coboundary, which yields for any integer $n$ large enough the following chain map $\Ome\Tc\Ac\to X(\Tc\Bc)$:
\beq
\lefteqn{\chih^n_{t_0}(\Oman^+\Ec,\rho_*,D):=\sum_{k=0}^{n+1}\chi^k(\Oman^+\Ec,\rho_*,t_0D)}\non\\
&&+ \int_{t=0}^{t=t_0}\big((\bb\oplus\nat\dd)cs^{n+1}(tD)-(-)^i(cs^{n+1}(tD)+cs^{n+2}(tD))\circ B\big)\ ,\non
\eeq
where $i=0,1$ is the parity of the bimodule, and $cs^n(tD)=cs^n(\Oman^+\Ec,\rho_*,tD)$. The property of this new cocycle is that it vanishes on $\Om^k\Tc\Ac$ for $k>n+1$. Its cohomology class in the complex $\hom(\Ome\Tc\Ac, X(\Tc\Bc))$ does not depend on $n$ nor $t_0$, and coincides with the Chern character of $(\Ec,\rho,D)$. The next step is then to take the limit $t_0\to \infty$ of the above cocycle. This requires to modify slightly the triple $(\Ec,\rho,D)$ so that the operator $D$ becomes invertible. Consider the $\cc$-$\cc$-bimodule of even degree $(\cc^2,\al,H_m)$ where $\cc^2$ is given its natural $\zz_2$-graduation, and the homomorphism $\al:\cc\to M_2(\cc)$ and the Fredholm operator $H_m$ read
$$
\al(1)=\left( \begin{array}{cc}
          1 & 0 \\
          0 & 0 \\
     \end{array} \right)\ ,\qquad H_m=\left( \begin{array}{cc}
          0 & m \\
          m & 0 \\
     \end{array} \right)\ ,
$$
for a given mass term $m\in\rr$. Then the graded tensor product of bimodules $(\Ec',\rho',D')=(\Ec,\rho,D)\hotimes(\cc^2,\al,H_m)$ defined by
$$
\Ec'=\Ec\hotimes\cc^2\ ,\quad \rho'=\rho\hotimes\al\ ,\quad D'=D\otimes 1+1\otimes H_m\ ,
$$
is an $\Ac$-$\Bc$-bimodule homotopic to $(\Ec,\rho,D)$, modulo addition of a degenerate bimodule (i.e. for which the homomorphism $\rho$ is zero). Therefore, the homotopy invariance of the bivariant Chern character with respect to the parameter $m$ shows that the chain map $\chi(\Oman^+\Ec,\rho_*,tD)$ is cohomologous to $\chi(\Oman^+\Ec',\rho'_*,tD')$ for any choice of $m\in\rr$. Furthermore, for $m\neq 0$, the operator $D'$ is invertible since $D'^2=D^2+m^2$. Using the estimates established in the proof of Proposition \ref{pana}, it is not hard to show that $\lim_{t_0\to\infty}\chi^k(\Oman^+\Ec',\rho'_*,t_0D')=0$ for any $k$, and the limit 
$$
\chih^n_{\infty}(\Oman^+\Ec',\rho'_*,D') = \int_0^{\infty}\big((\bb\oplus\nat\dd)cs^{n+1}(tD')-(-)^i(cs^{n+1}(tD')+cs^{n+2}(tD))\circ B\big)
$$
is a well-defined chain map $\Ome\Tc\Ac\to X(\Tc\Bc)$ representing the bivariant Chern character of $(\Ec,\rho,D)$, provided $n$ is large enough. As expected, the cohomology class of $\chih^n_{\infty}(\Oman^+\Ec',\rho'_*,D')$ enjoys good invariance properties, in particular with respect to suitable homotopies of the operator $D'$. For any $u\in(0,1]$ the negative powers of the selfadjoint operator $|D'|$ can be defined as
$$
|D'|^{-u}=C(u)\int_0^{\infty}d\la\,\frac{\la^{-u/2}}{\la+|D'|^2}\ ,
$$
where $C(u)$ is a normalization factor, and for $u=0$ we set $|D'|^0=1$. For $u\in [0,1]$, define the one-parameter family of operators $D'_u=D'|D'|^{-u}$. It connects homotopically the operator $D'$ for $u=0$, to its phase $F':=D'/|D'|$ for $u=1$. The operator $F'$ is bounded on the Hilbert space $\Hc\hotimes\cc^2$ and verifies the $p$-summability condition 
$$ 
[F',\rho'(a)] \in \ell^p(\Hc\hotimes\cc^2)\hotimes\Bc
$$
for any $a\in\Ac$ and $p>$ dim $M$. Here $\ell^p$ is the Schatten $p$-class on the Hilbert space $\Hc\hotimes\cc^2$. The homotopy of operators $D'_u$ and Proposition 4.4 of \cite{P2} imply that $D'$ can be replaced by its phase $F'$ in the formula of the cocycle $\chih^n_{\infty}(\Oman^+\Ec',\rho'_*,D')$ without changing the cohomology class. This invariance property is based on a double transgression formula involving the derivative $\frac{d}{du}D'_u= -\ln |D'|\, D'_u$. This yields new cocycles
$$
\chih^n_{\infty}(\Oman^+\Ec',\rho'_*,F'):\Ome\Tc\Ac\to X(\Tc\Bc)
$$
in any degree $n$ sufficiently large (it actually suffices to choose $n>p-1$ where $p$ is the summability degree). They all are cohomologous in bivariant entire cyclic cohomology and represent the Chern character of the $p$-summable Fredholm $\Ac$-$\Bc$-bimodule $(\Ec',\rho',F')$ (the bounded version of $(\Ec,\rho,D)$). Proposition 4.6 of \cite{P2} gives an explicit formula for the cocycle $\chih^n_{\infty}(\Oman^+\Ec',\rho'_*,F')$, when the degree $n$ is chosen with the same parity as the bimodule. It vanishes on $\Om^k\Tc\Ac$ if $k$ is different from $n$ or $n+1$, and for any elements $x_0,\ldots ,x_{n+1}\in\Tc\Ac$ one has
\beq
\lefteqn{\chih^n_{\infty}(\Oman^+\Ec',\rho'_*,F')(x_0\dd x_1\ldots \dd x_n)=}\non\\
&&\qquad\qquad (-)^n\frac{\Gamma(1+\frac{n}{2})}{(n+1)!}\, \frac{1}{2}\sum_{\la\in S_{n+1}}\eps(\la)\tau(F'[F',\rho_{\la(0)}]\ldots[F',\rho_{\la(n)}])\non
\eeq
as an element of $\Tc\Bc$, and
\beq
\lefteqn{\chih^n_{\infty}(\Oman^+\Ec,\rho'_*,F')(x_0 \dd x_1\ldots \dd x_{n+1})=}\non \\
&&\qquad (-)^n\frac{\Gamma(1+\frac{n}{2})}{(n+1)!}\, \frac{1}{2}\nat\tau\Big( \dd(\rho_0F'[F',\rho_{1}]\ldots[F',\rho_{n+1}])\non\\
&&\qquad\qquad\qquad\qquad +\sum_{\la\in S_{n+2}}\eps(\la) F'[F',\rho_{\la(0)}]\ldots[F',\rho_{\la(n)}]\dd \rho_{\la(n+1)}\Big)\non
\eeq
as an element of $\Om^1\Tc\Bc_{\nat}$. Here, $S_{n+1}$ denotes the group of cyclic permutations on $n+1$ elements, and we have adopted the notation $\rho_i=\rho'_*(x_i)$. It is instructive to look at the commutators. One sees that if $x\in \Om^+\Ac=T\Ac$ is an element of the non-completed tensor algebra, the commutator $[F',\rho(x)]$ is an endomorphism of $\Oman^+\Ec'$ which belongs to the subalgebra of $p$-summable operators
$$
\ell^p(\Hc\hotimes\cc^2)\hotimes\Om^+\Bc \subset \End_{\Tc\Bc}(\Oman^+\Ec')\ .
$$
Because $n>p-1$, the product of commutators $[F,\rho_i]$ in the above cocycle yields a trace-class endomorphism (recall by the way that the product on $\Om^+\Bc$ is Fedosov),
$$
[F',\rho_{\la(0)}]\ldots[F',\rho_{\la(n)}]\in \ell^1(\Hc\hotimes\cc^2)\hotimes\Om^+\Bc\ ,
$$
so that taking the (super)trace $\tau$ of operators on $\Hc\hotimes\cc^2$ is well-defined. One shows that $\chih^n_{\infty}(\Oman^+\Ec',\rho'_*,F')$ extends to a cocycle on the analytic tensor algebra $\Tc\Ac$ using the H\"older inequality and the kinds of estimates established in the proof of Proposition \ref{pana}. The bivariant cyclic cohomology class of the cocycles $\chih^n_{\infty}(\Oman^+\Ec',\rho'_*,F')$ also enjoys an invariance property with respect to $p$-summable homotopies of the homomorphism $\rho'$ and the operator $F'$, see \cite{P2} for details.\\ 
We now restrict to the case where the parity of the bimodule is even. Then the Hilbert space $\Hc=L^2(E)$ is $\zz_2$-graded, and in usual $2\times 2$ matrix notation we have
$$
\Ec= \left( \begin{array}{c}
          \Ec_+ \\
          \Ec_- \\
     \end{array} \right)\ ,\quad
\rho(a)=\left( \begin{array}{cc}
          \rho_+(a) & 0 \\
          0 & \rho_-(a) \\
     \end{array} \right)\ ,\quad  D=\left( \begin{array}{cc}
          0 & D_- \\
          D_+ & 0 \\
     \end{array} \right)\ .
$$
The ``massive amplification'' $(\Ec',\rho',D')$ is also a $\zz_2$-graded bimodule  given in $2\times 2$ matrix notation by
$$
\Ec'= \left( \begin{array}{c}
          \Ec \\
          \Ec \\
     \end{array} \right)\ ,\quad
\rho'(a)=\left( \begin{array}{cc}
          \rho'_+(a) & 0 \\
          0 & \rho'_-(a) \\
     \end{array} \right)\ ,\quad  D'=\left( \begin{array}{cc}
          0 & D_m \\
          D_m & 0 \\
     \end{array} \right)\ ,
$$
with $\rho'_+=\left( \begin{array}{cc}
          \rho_+ & 0 \\
          0 & 0 \\
     \end{array} \right)$, $\rho'_-=\left( \begin{array}{cc}
          0 & 0 \\
          0 & \rho_- \\
     \end{array} \right)$, $D_m=D+m\Gamma$ and $\Gamma=\left( \begin{array}{cc}
          1 & 0 \\
          0 & -1 \\
     \end{array} \right)$. The phase of $D'$ is therefore
$$
F'=\left( \begin{array}{cc}
          0 & F_m \\
          F_m & 0 \\
     \end{array} \right)\ ,\quad \mbox{with}\quad F_m=D_m/|D_m|\ ,\ F_m^2=1\ .
$$
Now, remark that the cocycle $\chih^n_{\infty}(\Oman^+\Ec',\rho'_*,F')$ is clearly equal to the cocycle $\chih^n_{\infty}(\Oman^+\Ec',\rho''_*,F)$ associated to the following bounded $\Ac$-$\Bc$-bimodule $(\Ec',\rho'',F)$, where the Fredholm operator $F$ is put into a canonical form:
$$
\rho''(a)=\left( \begin{array}{cc}
          F_m\rho'_+(a)F_m & 0 \\
          0 & \rho'_-(a) \\ 
     \end{array} \right)\ , \qquad F=\left( \begin{array}{cc}
          0 & 1 \\
          1 & 0 \\
     \end{array} \right)\ .
$$
We are now ready to evaluate these cocycles on the Chern character of the canonical $K$-theory class $[e]\in K_0(\Ac)$. Let $e\in\Ac$ be the idempotent obtained via a cut-off function $c\in \cinfc(M)$,
$$
e(g,x)=c(x)c(gx)\ .
$$
The image of $e$ under the homomorphism $\rho:\Ac\to \End_{\Bc}(\Ec)$ is the following idempotent operator
$$
(\rho(e)\xi)(g)=\int_Gdh\, c\, r(h)\, c\cdot \xi(gh^{-1}) \ ,\quad \forall \xi\in\Ec\ ,\ g\in G\ ,
$$
where the cut-off function $c$ acting on the sections of the vector bundle $E$ by pointwise multiplication, is viewed as an element of $\Lc(\Hc)$, and $r(h)\in\Uc(\Hc)$ is the unitary representation of the group. The idempotent $\rho(e)$ is related to the construction of the analytic assembly map of section \ref{sass}. To see this, let $\Pc$ be the algebra of properly supported, $G$-invariant pseudodifferential operators acting on $\cinfc(E)$, $\Pc_0\subset\Pc$ the subalgebra of operators of order $\leq 0$, and $\Jc\subset\Pc$ the ideal of smoothing operators. We regard $\Pc_0$ and $\Jc$ as subalgebras of $\Lc(\Hc)$. Recall that one has an homomorphism
$$
\te : \Jc\to C_c(G;\Kc)\ ,
$$
defined by $\te(T)(g)=c\,T\,r(g)\,c$ for any $T\in\Jc$. $\Kc$ denotes the algebra of smooth compact operators and $C_c(G;\Kc)$ is the convolution algebra of $\Kc$-valued continuous functions with compact support on $G$, see section \ref{sass}. In fact $C_c(G;\Kc)$ is a subalgebra of the endomorphisms $\End_{\Bc}(\Ec)$, and we can extend $\te$ to an homomorphism on the algebra of all pseudodifferential operators of order $\leq 0$
$$
\te: \Pc_0\to \End_{\Bc}(\Ec)
$$
by setting
$$
(\te(T)\xi)(g)= \int_Gdh\, c Tr(h)c\cdot \xi(gh^{-1})\qquad \forall\ T\in \Pc_0\ ,\ \xi\in\Ec\ .
$$
Since $\Pc_0$ has a unit (the identity operator on $\cinfc(E)$), a direct computation shows the equality
$$
\te(1)=\rho(e)
$$
of idempotents in $\End_{\Bc}(\Ec)$. Introduce the matrices 
$$
e_+=\left( \begin{array}{cc}
          1 & 0 \\
          0 & 0 \\
     \end{array} \right)\ ,\quad e_-=\left( \begin{array}{cc}
          0 & 0 \\
          0 & 1 \\
     \end{array} \right)\quad \in\Lc(\Hc)\ ,
$$
and let us split the idempotent $\te(1)=\rho(e)$ into two parts $\te(e_+),\te(e_-)$ according to the $\zz_2$-grading of $\End_{\Bc}(\Ec)$:
$$
\te(e_+)=\rho'_+(e)=\left( \begin{array}{cc}
          \rho_+(e) & 0 \\
          0 & 0 \\
     \end{array} \right)\ ,\qquad \te(e_-)=\rho'_-(e)=\left( \begin{array}{cc}
          0 & 0 \\
          0 & \rho_-(e) \\
     \end{array} \right)\ .
$$
Consider now the bounded $\cc$-$\Bc$-bimodule $(\Ec',e',F)$ defined by
$$
\Ec'= \left( \begin{array}{c}
          \Ec \\
          \Ec \\
     \end{array} \right)\ ,\quad
e'=\rho''(e)=\left( \begin{array}{cc}
          F_m\te(e_+)F_m & 0 \\
          0 & \te(e_-) \\
     \end{array} \right)\ ,\quad  F=\left( \begin{array}{cc}
          0 & 1 \\
          1 & 0 \\
     \end{array} \right)\ .
$$
The idempotent $e'$ is viewed as an homomorphism $\cc\to \End_{\Bc}(\Ec')$. The character $\chih^n_{\infty}(\Oman^+\Ec',e'_*,F)$ of this bounded bimodule is a cocycle in the complex $\hom(\Ome\Tc\cc,X(\Tc\Bc))$, that is, taking into account the homotopy equivalence $\Ome\Tc\cc\sim \cc$, a cycle in $X(\Tc\Bc)$. By definition the Chern character of $[e]\in K_0(\Ac)$ is represented by the entire cycle $\gamma(\eh)\in \Ome\Tc\Ac$, and because $e'=\rho''(e)$ one has the equality of entire cycles over $\Bc$:
$$
\chih^n_{\infty}(\Oman^+\Ec',\rho''_*,F)\cdot \gamma(\eh)=\chih^n_{\infty}(\Oman^+\Ec',e'_*,F)\ \in\Tc\Bc\ ,
$$
for any even integer $n>p-1$. The invariance of $\chih^n_{\infty}$ with respect to smooth homotopies allows to deform $(\Ec',e',F)$ into a bimodule representing $\mu([D])$. First, remark that the operator $F_m$ is $G$-invariant but not properly supported; hence $F_m\notin \Pc_0$. However, we can find a small deformation of $F_m$ which lies in $\Pc_0$. Indeed, let $Q:\cinfc(E_+)\to\cinfc(E_-)$ be the pseudodifferential operator of order $0$ defined by
$$
Q=D_+\cdot \delta_1\ ,
$$
where $\delta_1\in \Psi_c^{-1}(E_+,E_-)$ is the elliptic operator with symbol $\si(x,p)=(m^2+\|p\|^2)^{-1/2}$, as in section \ref{sass}. One can find a properly supported and $G$-invariant parametrix $P:\cinfc(E_-)\to\cinfc(E_+)$ for $Q$. Then $PQ-1$ and $QP-1$ are smoothing operators. We introduce as in section \ref{sass} the invertible operator $T\in\Pc_0$
$$
T=\left( \begin{array}{cc}
          1-PQ & P \\
          2Q-QPQ & QP-1 \\
     \end{array} \right)\ ,\qquad T=T^{-1}\ .
$$
Then $T$ is a small deformation of $F_m$, in the sense that the operators $(T-F_m)f$ and $f(T-F_m)$ are in $\ell^p(\Hc)$ for any compactly supported function $f\in\cinfc(M)$. We want to show that, upon stabilization by $2\times 2$ matrices, the idempotent $e'$ is $p$-summably homotopic to 
$$
e''=\left( \begin{array}{cc}
          \te(T^{-1}e_+T) & 0 \\
          0 & \te(e_-) \\
     \end{array} \right)\ .
$$
To this end, let us replace $e'$ by its stabilization defined on its diagonal elements as follows
$$
F_m\te(e_+)F_m \to \left( \begin{array}{cc}
          F_m\te(e_+)F_m & 0 \\
          0 & 0 \\
     \end{array} \right)\ ,\qquad \te(e_-)\to \left( \begin{array}{cc}
          \te(e_-) & 0 \\
          0 & 0 \\
     \end{array} \right)\ .
$$
For any parameter $t\in [0,1]$, introduce the invertible matrix
$$
U_t=R_t^{-1}\left( \begin{array}{cc}
          F_m & 0 \\
          0 & \te(T)+1-\te(1) \\
     \end{array} \right) R_t\ ,
$$
where $R_t$ is the rotation matrix
$$R_t=\left( \begin{array}{cc}
          \cos(\pi t/2) &  -\sin(\pi t/2)\\
           \sin(\pi t/2)&  \cos(\pi t/2)\\
     \end{array} \right)\ .
$$
We next observe that the operators $\te(T)-F_m\te(1)$ and $\te(T)-\te(1)F_m$ belong to the subalgebra of $p$-summable endomorphisms $\ell^p(\Hc)\hotimes\Bc$. Indeed by evaluation on a point $g\in G$,
$$
(\te(T)-F_m\te(1))(g)=cTr(g)c - F_mcr(g)c=-[F_m,c]r(g)c+ c(T-F_m)r(g)c
$$
is a $p$-summable operator on $\Hc$, and similarly with $\te(T)-\te(1)F_m$. This property allows to show by a direct computation that the difference of idempotents 
$$
U_t^{-1}\left( \begin{array}{cc}
          \te(e_+) & 0 \\
          0 & 0 \\
     \end{array} \right)U_t - \left( \begin{array}{cc}
          \te(e_-) & 0 \\
          0 & 0 \\
     \end{array} \right)
$$
lies in the subalgebra $\ell^p(\Hc)\hotimes M_2(\cc)\hotimes\Bc$ for any $t\in [0,1]$. This shows that after stabilization, the idempotent $e'$ corresponding to the value $t=0$ is smoothly and $p$-summably homotopic to the idempotent $e''$ corresponding to $t=1$. Hence homotopy invariance implies that the cocycle $\chih^n_{\infty}(\Oman^+\Ec',e'_*,F)$ is cohomologous to $\chih^n_{\infty}(\Oman^+\Ec',e''_*,F)$. Finally, it remains to observe that the $K$-theory class represented by the $\cc$-$\Bc$-bimodule $(\Ec',e'',F)$ is by definition the image of $[D]\in K_0^G(M)$ under the assembly map. Indeed, the idempotents $\te(T^{-1}e_+T)$ and $\te(e_-)$ lie in the subalgebra of endomorphisms $C_c(G;\Kc)= \Kc\hotimes C_c(G)$ (augmented by adjoining units to the diagonal blocks), and one has by definition
$$
\mu(D)=[\te(T^{-1}e_+T)]-[\te(e_-)]\in K_0(\Bc)\ .
$$
We already know that the product $\ch(\Ec,\rho,D)\cdot\ch(e)\in HE_0(\Bc)$ is represented by the entire cycle $\chih^n_{\infty}(\Oman^+\Ec',e''_*,F)$ for any even integer $n>p-1$. We claim that the latter also represents the Chern character of $\mu(D)$. In effect, since the difference $\te(T^{-1}e_+T)-\te(e_-)$ lies in a subalgebra of trace-class operators $\Kc\hotimes C_c(G)\subset \ell^1(\Hc)\hotimes\Bc$, the cycles $\chih^n_{\infty}(\Om^+\Ec',e''_*,F)$ are well-defined and cohomologous for any even degree $n$, including the smallest possible value $n=0$. One has
$$
\chih^0_{\infty}(\Om^+\Ec',e''_*,F)=\frac{1}{2}\Tr_s(F[F,e''_*])=\Tr(\te(T^{-1}e_+T)_*-\te(e_-)_*)=\ch(\mu(D))\ ,
$$
where more precisely $\te(T^{-1}e_+T)_*$ and $\te(e_-)_*$ are idempotents of the algebra $\ell^1\hotimes\Tc\Bc$ (augmented by adjoining units to the diagonal blocks), lifting the idempotents $\te(T^{-1}e_+T)$ and  $\te(e_-)$ of $\ell^1\hotimes\Bc$ (augmented). This proves the theorem in the even case.\\
In the odd case, the image of $[D]\in K_1^G(M)$ under the analytic assembly map $\mu(D)\in K_1(\Bc)$ is obtained by first taking the external product with the Bott class $\beta\in K_1^{\zz}(\rr)$, then composing with the assembly map in even degree:
$$
K_1^G(M)\stackrel{\beta\boxtimes}{\longrightarrow}K_0^{\zz\times G}(\rr\times M) \to K_0(\cc\zz\otimes C_c(G))\to K_1(\Bc)\ .
$$
Hence the equality $\ch(\mu(D))=\ch(\Ec,\rho,D)\cdot\ch(e)$ in the odd case stems from the above study of the even case together with the compatibility of the Chern character with Bott periodicity \cite{P1}. \cqfd\\

Theorem \ref{tind} enables to find a local index formula, after taking the limit $t\to 0$ of the cocycles $\chi(\Oman^+\Ec,\rho_*,tD)$. This formula involves some classical equivariant Chern classes for the various bundles defined over the manifold $M$, and a kind of ``noncommutative'' Chern character form associated to the canonical $K$-theory element $[e]\in K_0(\Ac)$. We first have to define the crossed-product algebra $\cinfc(M)\rtimes T\Bc$. Recall that the non-completed tensor algebra $T\Bc$ is the direct sum
\be
T\Bc=\Bc\oplus \Bc\hotimes\Bc\oplus \Bc\hotimes\Bc\hotimes\Bc\oplus \ldots
\ee
Since an $n$-fold tensor product $\Bc^{\hotimes n}$ is isomorphic, as a bornological vector space, to the space of integrable functions $L^1(G^n,d\nu^n)$ with respect to the admissible measure $d\nu$, we see that the product of two homogeneous tensors $x\in \Bc^{\hotimes n}$ and $y\in \Bc^{\hotimes m}$ is a function of $n+m$ points on $G$:
\be
(xy)(g_1,\ldots,g_{n+m})=x(g_1,\ldots,g_n)y(g_{n+1},\ldots,g_{n+m})\ .
\ee
We form the crossed product $\cinfc(M)\rtimes T\Bc$ as follows. As a vector space, it is the completed tensor product 
\be
\cinfc(M)\hotimes T\Bc=\bigoplus_{n\ge 1}\cinfc(M)\hotimes \Bc^{\hotimes n}=\bigoplus_{n\geq 1}\mathop{\limind}_{K}L^1(G^n;\cinf_K(M))\ ,
\ee
where the inductive limit is taken over the compact subsets $K\subset M$ and $\cinf_K(M)$ is the space of smooth functions with support in $K$. The multiplication of two elements $x\in L^1(G^n;\cinf_{K_1}(M))$ and $y\in L^1(G^m;\cinf_{K_2}(M))$ is defined by
\be
(xy)(g_1,\ldots,g_{n+m})=x(g_1,\ldots,g_n)y(g_{n+1},\ldots,g_{n+m})^{g_n\ldots g_1}\ ,\label{cross}
\ee
where the superscript $^{g_n\ldots g_1}$ denotes the pullback of the smooth function with compact support $y(g_{n+1},\ldots,g_{n+m})\in \cinfc(M)$ by the diffeomorphism $g_n\ldots g_1 \in G$. We obtain an algebra homomorphism $\psi:T\Ac\to \cinfc(M)\rtimes T\Bc$ by setting
\be
\psi(a_1\otimes\ldots\otimes a_n)(g_1,\ldots,g_n)=a_1(g_1)a_2(g_2)^{g_1}\ldots a_n(g_n)^{g_{n-1}\ldots g_1}\ .\label{homo}
\ee
Alternatively, we may describe this homomorphism in terms of the Fedosov algebras of differential forms $(\Om^+\Ac,\odot)=T\Ac$ and $(\Om^+\Bc,\odot)=T\Bc$. Observe that restricted to the $2k$-forms over $\Ac$, $\psi$ yields a linear map
\be
\psi:\Om^{2k}\Ac \to \cinfc(M)\hotimes\Om^{2k}\Bc\ .
\ee
Moreover, the space $\Om^{2k}\Bc=\Bc^{\hotimes 2k}\oplus \Bc^{\hotimes (2k+1)}$ may be identified with a space of functions over the locally compact space $G^{2k}\cup G^{2k+1}$, and one has
\be
\psi(a_0da_1\ldots da_{2k})(g_0,\ldots,g_{2k})=a_0(g_0)a_1(g_1)^{g_0}\ldots a_{2k}(g_{2k})^{g_{2k-1}\ldots g_0}\label{medor1}
\ee
at any point $(g_0,\ldots,g_{2k}) \in G^{2k+1}$, and similarly
\be
\psi(da_1\ldots da_{2k})(g_1,\ldots,g_{2k})=a_1(g_1)a_2(g_2)^{g_1}\ldots a_{2k}(g_{2k})^{g_{2k-1}\ldots g_1}\label{medor2}
\ee
at any point $(g_1,\ldots,g_{2k}) \in G^{2k}$. Next, let $\Om_c(M)$ denote the differential graded (DG) algebra of smooth differential forms with compact support on $M$. The supremum norm of a differential form and all its derivatives over compact subsets is defined with the help of the Riemannian metric. Hence $\Om_c(M)$ is an LF-algebra as a union of Fr\'echet spaces, and we gift it with the corresponding bounded bornology. Denote by $\delta$ the de Rham coboundary on $\Om_c(M)$. Replacing $\cinfc(M)$ by $\Om_c(M)$ in the above construction, we define the crossed product algebra $\Om_c(M)\rtimes T\Bc$. The multiplication of two elements is given by a formula analogous to (\ref{cross}), involving the pullback of differential forms by diffeomorphisms. Moreover, the de Rham coboundary induces a differential 
\be
\delta : \Om^n_c(M)\hotimes T\Bc \to \Om^{n+1}_c(M)\hotimes T\Bc
\ee
by acting only on the first factor, which turns $\Om_c(M)\rtimes T\Bc$ into a DG algebra. $\cinfc(M)\rtimes T\Bc$ is of course its subalgebra of degree zero forms. Consider now the space of noncommutative one-forms $\Om^1T\Bc$. It is a bimodule over $T\Bc$, isomorphic to $\widetilde{T}\Bc\hotimes\Bc\hotimes\widetilde{T}\Bc$ (see section \ref{sent}, $\widetilde{T}\Bc$ is the unitalization of $T\Bc$). Now it is clear how to endow the tensor product $\Om_c(M)\hotimes \Om^1T\Bc$ with a bimodule structure over the algebra $\Om_c(M)\rtimes T\Bc$: the left and right module maps are given again by formulas like (\ref{cross}); the elements of $\Om_c(M)$ are pulled-back by diffeomorphisms whenever they cross a tensor product $\Bc^{\hotimes n}$. Finally, there is a derivation
\be
\dd:\Om_c(M)\rtimes T\Bc\to \Om_c(M)\hotimes \Om^1 T\Bc
\ee
induced by the universal derivation $\dd:T\Bc\to \Om^1T\Bc$ acting only on the second factor (one has nevertheless to put a minus sign when $\dd$ crosses a differential form of odd degree in $\Om_c(M)$). We will use the DG algebra $\Om_c(M)\rtimes T\Bc$ and its bimodule $\Om_c(M)\hotimes \Om^1 T\Bc$ to construct a linear map
\be
\Psi: \Om T\Ac \to \Om_c(M)\hotimes X(T\Bc)\ ,
\ee
where $\Om T\Ac$ is the DG algebra of noncommutative forms over the tensor algebra $T\Ac$. Here we ignore the various boundary maps on the complexes $\Om T\Ac$, $\Om_c(M)$ or $X(T\Bc)$, because $\Psi$ will not exactly be a chain map. For any $n$-form $x_0\dd x_1\ldots\dd x_n \in \Om^nT\Ac$, set
\beq
\Psi(x_0\dd x_1\ldots\dd x_n) &=&  \frac{1}{n!}\, \psi(x_0)\delta\psi(x_1)\ldots \delta\psi(x_n) \\
 &+& \frac{1}{n!}\sum_{i=1}^n \nat \psi(x_0)\delta\psi(x_1)\ldots \dd \psi(x_i)\ldots \delta\psi(x_n)\ .\non
\eeq
The first term of the r.h.s. is an element of $\Om^n_c(M)\hotimes T\Bc$, whereas the second lies in $\Om_c^{n-1}(M)\hotimes \Om^1T\Bc_{\nat}$. Next, using the identification of $X(T\Bc)$ with the space of noncommutative forms $\Om\Bc$, the map $\Psi: \Om T\Ac \to \Om_c(M)\hotimes \Om\Bc$ may be decomposed into components of any degree $n\in\nn$:
\be
\Psi_n: \Om T\Ac \to \Om_c(M)\hotimes \Om^n\Bc\ .
\ee
The important remark is that each map $\Psi_n$ immediately extends to the completed space $\Ome\Tc\Ac$. This is due to the fact that the degree of differential forms in $\Om_c(M)$ and $\Om^n\Bc$ cannot exceed a given number (respectively, dim $M$ and $n$). This leads to the following definition.
\begin{definition}\label{dchern}
Let $e\in\Ac$ be an idempotent, $\eh\in \Tc\Ac$ its idempotent lift and $\gamma(\eh) \in\Ome\Tc\Ac$ its Chern character in the complex of entire forms over $\Tc\Ac$. For any $n\in\nn$, we define the mixed differential form
\be
\ch_n(e):= \Psi_n(\gamma(\eh))\in \Om_c(M)\hotimes \Om^n\Bc\ .
\ee
The isomorphism $\Om^n\Bc=\Bc^{\hotimes n}\oplus \Bc^{\hotimes (n+1)}$ shows that $\ch_n(e)$ may be identified with a continuous and compactly supported function on the locally compact space $G^n\cup G^{n+1}$, with values in $\Om_c(M)$.
\end{definition}
The other ingredients of the local index formula are well-known. They already appear in the Atiyah-Segal-Singer index theorem \cite{AS} or the Lefschetz fixed point theorem \cite{Gi}. Let us recall below some of the geometric apparatus we need.\\
For simplicity, we will restrict ourselves to $K$-homology classes $[D]\in K_*^G(M)$ of \emph{Dirac type}. This means that $M$ must be an oriented manifold, $G$ acts on $M$ by orientation-preserving diffeomorphisms, the vector bundle $E$ is a $G$-equivariant Clifford module, and $D$ is a generalized Dirac operator, in the following sense. Since $G$ acts on $M$ by diffeomorphisms, it acts on the smooth sections of the cotangent bundle $T^*M$ (i.e. one-forms on $M$) by pullback. Because $M$ is endowed with a $G$-invariant Riemannian metric $(\ ,\ ):\cinf(T^*M)\times\cinf(T^*M)\to \cinf(M)$, $G$ acts by automorphisms on the algebra of smooth sections of the Clifford bundle $Cl(T^*M)$. The latter is the bundle of Clifford algebras generated by the linear inclusion $c:T^*M_x\to Cl(T^*M)_x$ at any point $x\in M$ and relations
\be
c(\al)c(\beta)+c(\beta)c(\al)=-2( \al,\beta)\qquad \forall\al,\beta\in T^*M_x\ .
\ee
The vector bundle $E\to M$ is a Clifford module if it is endowed with a left action $Cl(T^*M)\times E\to E$ of the Clifford bundle. We say that $E$ is a $G$-equivariant Clifford module if the action of $Cl(T^*M)$ is equivariant:
\be
(c(\al)\xi)^g=c(\al^g)\xi^g\qquad\forall \al\in\cinf(T^*M)\ ,\ \xi\in\cinf(E)\ ,\ g\in G\ .
\ee
A Clifford connection $\nabla^E:\cinf(E)\to\cinf(T^*M\otimes E)$ is a connection on $E$ compatible with the Levi-Civita connection $\nabla^{LC}$ on $T^*M$, in the sense that the equality 
\be
[\nabla^E_X,c(\al)]\xi=c(\nabla^{LC}_X\al)\xi
\ee
holds for any $\xi\in \cinf(E)$, $\al\in\cinf(T^*M)$, and $X\in \mbox{Vect}(M)$. The connection $\nabla^E$ is said to be $G$-invariant if it is equivariant as a map $\cinf(E)\to\cinf(T^*M\otimes E)$. From a $G$-invariant Clifford connection, we may form a $G$-invariant generalized Dirac operator
\be
D:\cinf(E)\to\cinf(E)\ ,\qquad D=\sum_ic(dx^i)\nabla^E_i\ ,
\ee
expressed in terms of a local coordinate system $\{x^i\}_{i=1,\ldots,n}$, $n=$ dim $M$. The $K$-homology class of the operator $D$ thus obtained is by definition of Dirac type.\\
We now describe the various characteristic classes involved in the equivariant index theorem. Choose a local coordinate chart $\{x^i\}_{i=1,\ldots,n}$ and a local field of orthonormal frames $\{e^a\}_{a=1,\ldots,n}$. The components of the Levi-Civita connection along the direction $\frac{\partial}{\partial x^i}$ may be expressed as
\be
\nabla^{LC}_ie^a=-\sum_b\om_{iab}e^b\ ,
\ee
where $\om_{iab}$ is antisymmetric in $a,b$. We deduce the components of the Riemannian curvature through the commutator 
\be
[\nabla^{LC}_i,\nabla^{LC}_j]e^a=-\sum_b R_{ijab}e^b\ ,
\ee
which yields
\be
R_{ijab}= \partial_i\om_{jab}-\partial_j\om_{iab} + \sum_c(\om_{iac}\om_{jcb} - \om_{jac}\om_{icb})\ .
\ee
We denote by $R$ the curvature two-form, determined by the antisymmetric matrix with two-form coefficients $R_{ab}=\frac{1}{2}R_{ijab}dx^i\wedge dx^j$. Thus, $R$ defines a $G$-invariant section of the vector bundle $\La^2T^*M\otimes \End(T^*M)$.\\
Next, introduce the spin connection $\nabla^S:\cinf(E)\to\cinf(T^*M\otimes E)$, given in the local coordinate chart $\{x^i\}$ by
\be
\nabla^S_i=\partial_i-\frac{1}{4}\sum_{a,b}\om_{iab}c(e^a)c(e^b)\ .
\ee
One checks that $\nabla^S$ is a $G$-invariant Clifford connection. The \emph{relative connection} is the difference
\be
A^{E/S}:= \nabla^E-\nabla^S\ .
\ee
It defines a $G$-equivariant bundle map from $E$ to $T^*M\otimes E$, commuting with the Clifford action:
\be
[A^{E/S},c(\al)]=0\qquad \forall \al\in\cinf(T^*M)\ .
\ee
The relative curvature is the two-form on $M$ with values in the endomorphisms of $E$ given by
\be
F^{E/S}=[\nabla^S,A^{E/S}]+(A^{E/S})^2\ ,
\ee
and defines a $G$-invariant section of the vector bundle $\La^2T^*M\otimes \End(E)$.\\ 
Now choose an element $g\in G$. The set of fixed points for the diffeomorphism $g\in \diff(M)$ is an union of isolated submanifolds $M_g\subset M$, which may be of different dimensions. The centralizer $Z_g\subset G$ of $g$ is the subgroup
\be
Z_g=\{z\in G| zg=gz\}\ .
\ee
It leaves the fixed point set globally invariant and acts by diffeomorphisms on each fixed manifold $M_g$ of a given dimension. The restriction of the Riemannian curvature $R$ to $M_g$ may therefore be considered as a $Z_g$-invariant section of the vector bundle $\La^2T^*M_g\otimes \End(T^*M)$ over $M_g$. Using the Riemannian metric on $M$, we may identify the tangent and cotangent bundles $TM$ and $T^*M$. Let $N$ be the orthogonal complement of the tangent bundle of $M_g$ in $TM$:
\be
TM=TM_g\oplus N\ .
\ee
$N$ is the normal bundle over $M_g$. Since $g$ is an orientation-preserving isometry, the tangent space $TM_x$ over a fixed point $x\in M_g$ carries a linear representation of $g$ via an element $r^T(g)_x$ of the orthogonal group $SO(n)$, $n=\mbox{dim}\, M$. Using the orthogonal decomposition $TM_x= (TM_g)_x\oplus N_x$, the normal fiber $N_x$ is actually stable by $r^T(g)_x$ and $(TM_g)_x$ is the fixed subspace. It follows that the curvature $R$, as an antisymmetric matrix with two-form coefficients commuting with $g$, splits according to the above decomposition $TM_g\oplus N$:
\be
R=\left( \begin{array}{cc}
    R^0 & 0 \\
    0 & R^1 \end{array} \right)\ ,
\ee
where $R^0\in \La^2T^*M_g\otimes \End(TM_g)$ is the curvature of the Levi-Civita connection of the submanifold $M_g$, and $R^1\in \La^2T^*M_g\otimes \End(N)$ is the curvature of the normal bundle $N$. In a similar way, the restriction of the relative curvature $F^{E/S}$ to the submanifold $M_g$ may be considered as a $Z_g$-invariant section of the vector bundle $\La^2T^*M_g\otimes \End(E)$. We introduce the Atiyah-Hirzebruch genus of $M_g$ as the closed $Z_g$-invariant differential form over $M_g$
\be
\widehat{A}(M_g)={\det}^{1/2}\left(\frac{R^0/2}{\sinh(R^0/2)}\right)\ .
\ee 
Next, because any point $x\in M_g$ is fixed by the action of $g$, the fiber $E_x$ also carries a linear representation $r^E(g)_x$ of $g$. Let $r^S(g)$ be a lift of $r^T(g)$ into the group $\mbox{Spin}(n)$, in a neighborhood $U\subset M_g$ of a point $x_0$. It is defined up to a sign. For any $x\in U$, $r^S(g)$ acts on the fiber $E_x$ via Clifford multiplication. We define the \emph{relative action} $r^{E/S}(g)$ on $E_x$ as
\be
r^{E/S}(g)=r^S(g)^{-1}r^E(g)\ .
\ee
It commutes with the Clifford action: $[r^{E/S}(g),c(\al)]=0$ for any $\al\in T^*M_x$. Now, if $U$ is small enough, the bundle $E$ over $U$ is decomposed into a tensor product of a spinor bundle $S$ (i.e. an irreducible faithful representation of the Clifford bundle) with a relative bundle $E/S$ on which $Cl(T^*M)$ acts by the identity:
\be
E=S\otimes E/S\ .\label{loc}
\ee
Hence, $r^S(g)$ acts on $S$ while $r^{E/S}(g)$ acts on $E/S$. Also, the relative curvature $F^{E/S}$ acts on the factor $E/S$ only. Both bundles $S$ and $E/S$ can be $\zz_2$-graded, depending on the dimension of $M$ and other conventions. We define the equivariant Chern character of the relative bundle $E/S$ as the closed differential form over $U$
\be
\ch(E/S,g)=\tr_s\big(r^{E/S}(g)\exp(-F^{E/S})\big)\ ,
\ee
where $\tr_s$ is the (super)trace of endomorphisms of $E/S$. This expression is not always globally defined on $M_g$, because of the sign chosen in the spin lift $r^S(g)$. Now consider the normal bundle $N$ over $U$. Let $q$ be the dimension of the fibers. Then at each point of $U$, the tangent representation $r^T(g)$ is actually an element $g_N$ of the orthogonal group $SO(q)\subset SO(n)$ of the normal fiber, and its lift $r^S(g)\in \mbox{Spin}(q)$ lies in the normal Clifford bundle $Cl(N)\subset Cl(TM)$. It should be noticed that, since $g_N\in SO(q)$ has no fixed eigenvector, the codimension $q$ of $M_g$ is necessarily \emph{even}. We let $S_N$ be the ($\zz_2$-graded) spinor bundle over $U$ associated to $Cl(N)$. The curvature $R^1$ of the normal bundle is a two-form with values in the Lie algebra of the orthogonal group $SO(q)$, hence it acts on the normal spinor bundle via the spin representation $r^S(R^1)$. Define the equivariant Chern character of the normal spinor bundle $S_N$ as the closed differential form over $U$
\be
\ch(S_N,g)=\tr_s\big(r^S(g\exp(-R^1))\big)\ ,
\ee
where $\tr_s$ is the (super)trace of endomorphisms of $S_N$, defined once an orientation of the fibers of $N$ is chosen. $\ch(S_N,g)$ can be expressed directly in terms of the representations of $g$ and $R^1$ on the normal bundle $N$:
\be
\ch(S_N,g)=i^{q/2}{\det}^{1/2}(1-g_N\exp(-R^1))\ .
\ee
In particular, the degree zero component of this differential form is the locally constant function $i^{q/2}{\det}^{1/2}(1-g_N)$ over $U$. It never vanishes because by definition, $g$ has no fixed point in the normal bundle $N$ except the zero section. Let us calculate it in more geometrical terms. Since the codimension $q$ of $M_g$ is even, the normal fiber $N_x$ at a point $x\in M_g$ can be decomposed into an orthogonal direct sum of planes,
\be
N_x=\bigoplus_{k=1}^{q/2}N^k_x\ ,\qquad N^k_x=\rr^2\ ,
\ee
and $g_N$ acts on each plane $N^k_x$ by a rotation of angle $\te^k\neq 0$. Therefore, one has
\be
\det(1-g_N)=\prod_{k=1}^{q/2}4\sin^2(\te^k/2)\ .
\ee
Consequently, the differential form $\ch(S_N,g)$ is always invertible. Once again, it may not be globally defined on $M_g$, due to the choice of sign in the spin lift $r^S(g)$ and the choice of orientation for the fibers of $N$. However, $M$ being oriented, one can fix a local orientation on $M_g$ such that the induced orientation of the fibers $TM_g\oplus N$ coincides with the orientation of $TM$ over $U$. A global orientation on $M_g$ may not exist, but with some abuse of notation, we see that the cap-product
\be
C=\widehat{A}(M_g)\frac{\ch(E/S,g)}{\ch(S_N,g)}\cap [M_g]
\ee
is a globally well-defined de Rham current over $M$. Indeed, the local sign ambiguities in defining $\ch(E/S,g)$, $\ch(S_N,g)$ and the fundamental class $[M_g]$ are killed when they are multiplied altogether. Actually $C$ is closed and $Z_g$-invariant. \\
From theorem \ref{tind}, we know that the image of the $K$-homology class $[D]\in K_*^G(M)$ under the assembly map has a Chern character $\ch\circ\mu(D)\in HE_0(\Bc)$ represented by the entire cycle
\be
\ch(tD):=\chi(\Oman^+\Ec,\rho_*,tD)\circ\gamma(\eh)\in X(\Tc\Bc)\cong \Oman\Bc
\ee
for any parameter $t>0$, and $\eh\in\Tc\Ac$ is the idempotent lift (\ref{lif}). We let $\ch_n(tD)\in \Om^n\Bc$ be the $n$th degree component of this cycle. Using the decomposition $\Om^n\Bc=\Bc^{\hotimes n}\oplus \Bc^{\hotimes (n+1)}$, we see that $\ch_n(tD)$ is a function on $G^n\cup G^{n+1}$. The corollary below expresses the pointwise limit as $t\to 0$ of this function, in terms of the (classical) equivariant characteristic classes and the noncommutative Chern character of definition \ref{dchern}.
\begin{corollary}\label{cloc}
Let $D:\cinfc(E)\to\cinfc(E)$ be a $G$-invariant generalized Dirac operator representing a $K$-homology class of Dirac type $[D]\in K_*^G(M)$, and let $e\in \cinfc(M)\rtimes G$ be an idempotent representing the canonical $K$-theory class. For any admissible completion $\Bc$ of the convolution algebra $C_c(G)$, the Chern character of $\mu(D)\in K_*(\Bc)$ is represented by the entire cycle $\ch(tD)\in \Oman\Bc$, $t>0$. Its $n$th degree component $\ch_n(tD)\in \Om^n\Bc$ is a function over $G^n\cup G^{n+1}$ whose pointwise limit as $t\to 0$ is given by the localization formula
\be
\lim_{t\to 0}\ch_n(tD)(\gt)=\sum_{M_g}\frac{(-)^{q/2}}{(2\pi i)^{d/2}}\int_{M_g}\widehat{A}(M_g)\frac{\ch(E/S,g)}{\ch(S_N,g)}\, \ch_n(e)(\gt)\ ,
\ee
where $\gt\in G^n\cup G^{n+1}$, and $g\in G$ is the concatenation product $g_n\ldots g_1$ (resp. $g_n\ldots g_0$) if $\gt=(g_1,\ldots, g_n)$ (resp. $\gt=(g_0,\ldots, g_n)$). The sum runs over the fixed manifolds $M_g$ of all possible dimensions $d$ and codimensions $q=\mbox{dim}\, M-d$. Finally the Chern character $\ch_n(e)$ of definition \ref{dchern} is viewed as a function on $G^n\cup G^{n+1}$ with values in $\Om_c(M)$.
\end{corollary}
{\it Proof:} Recall that the idempotent $\eh\in\Tc\Ac$ is the analytic differential form (\ref{lif}):
$$
\eh=e+ \sum_{k\ge 1} \frac{(2k)!}{(k!)^2} (e-\frac{1}{2}) (dede)^{k}\in \Oman^+\Ac=\Tc\Ac\ ,
$$
whereas its image $\gamma(\eh)$ in the $(b+B)$-complex of entire chains $\Ome\Tc\Ac$ is the entire cycle
$$
\gamma(\eh)= \eh+ \sum_{n\ge 1} (-)^n\frac{(2n)!}{n!} (\eh-\frac{1}{2}) (\dd\eh \dd\eh)^n
$$
with $\dd$ the differential of noncommutative forms over $\Tc\Ac$. When the parity of $[D]$ is even, the composite $\ch(tD)=\chi(\Oman^+\Ec,\rho_*,tD)\circ\gamma(\eh)$ is a cycle of even degree in the $X$-complex $X(\Tc\Bc)$ given by the series
$$
\ch(tD)=\chi^0_0(\eh) +\sum_{n\geq 1}(-)^n\frac{(2n)!}{n!} \chi^{2n}_0 ((\eh-\frac{1}{2}) (\dd\eh \dd\eh)^n) \in \Tc\Bc\ ,
$$
whereas in the odd case, 
$$
\ch(tD)=\chi^0_1(\eh) +\sum_{n\geq 1}(-)^n\frac{(2n)!}{n!} \chi^{2n}_1 ((\eh-\frac{1}{2}) (\dd\eh \dd\eh)^n) \in \Om^1\Tc\Bc_{\nat}\ .
$$
Recall that the components $\chi^n_0:\Om^n\Tc\Ac \to \Oman^+\Bc\cong \Tc\Bc$ and $\chi^n_1:\Om^n\Tc\Ac \to \Oman^-\Bc\cong \Om^1\Tc\Bc_{\nat}$ are defined by equations (\ref{joe3}) and (\ref{joe4}). We will concentrate only on the even case since the formulas are simpler (the odd case is similar). One wishes to calculate the limit $t\to 0$ of the components of the cycle $\ch(tD)\in \Oman^+\Bc$ in all degrees. Let us denote by $p_{2m}$ the projection of $\Oman^+\Bc$ onto the subspace of $2m$-forms $\Om^{2m}\Bc$. Then one has
\be
\ch_{2m}(tD)=p_{2m}\chi^0_0(\eh) +\sum_{n\geq 1}(-)^n\frac{(2n)!}{n!} p_{2m}\chi^{2n}_0 ((\eh-\frac{1}{2}) (\dd\eh \dd\eh)^n) \ ,\label{ser}
\ee
and each term of the series over $n$ is actually an element the subspace of continuous, compactly supported forms $\Om^{2m}_c(G)\subset \Om^{2m}\Bc$ (see section \ref{sbiv}) given by a \emph{finite sum}
$$
p_{2m}\chi^{2n}_0 ((\eh-\frac{1}{2}) (\dd\eh \dd\eh)^n)= \sum_{i_0,\ldots, i_{2n}}p_{2m}\chi^{2n}_0 (x_{i_0}\dd x_{i_1}\ldots \dd x_{i_{2n}})\ .
$$
Here the $x_i$'s are chosen amongst the components $(e-\frac{1}{2}) (dede)^{k}$ of the idempotent lift $\eh$ in degree $\leq 2m$. Hence, all the $x_i$'s lie in the non-completed tensor algebra $\Om^+\Ac=T\Ac$. Using the estimates established in the proof of Proposition \ref{pana}, it is possible to control the supremum norm of $p_{2m}\chi^{2n}_0 ((\eh-\frac{1}{2}) (\dd\eh \dd\eh)^n)$ as a continuous function with compact support on $G^{2m}\cup G^{2m+1}$. Indeed, we know that if $S\subset\Ac$ is a small subset, and $x_i=a^i_0da^i_1\ldots da^i_{2k_i}\in \Om^{2k_i}\Ac$ are such that all $a^i_j$'s are contained in $ S$, then there is a compact subset $K_S\subset G$ and a constant $\zeta_S$ depending only on $S$, such that 
$$
\sup_{G^{2m}\cup G^{2m+1}}|p_{2m}\chi^{2n}_0(x_0\dd x_1\ldots\dd x_{2n})| \leq  \Tr(fe^{-t^2D^2}f)\, \frac{(t\zeta_S\,|K_S|)^{2n}}{(2n)!}\, (4\zeta_S)^{2m+1}\ ,
$$
for some plateau function $f\in\cinfc(M)$. Let us now use $p$-summability. From the proof of Theorem \ref{tind}, one knows that for small values of $t$,
$$
\Tr(fe^{-t^2D^2}f)\sim \frac{1}{t^p}
$$
for any real number $p>$ dim $M$. Consequently, the supremum of $\chi^{2n}_0$ is of order $t^{2n-p}$, and the tail of the series (\ref{ser}), for $2n >$ dim$M$, vanishes at the limit $t\to 0$. We will show now that in any degree 2$n\in [0,\mbox{dim}\, M]$, the pointwise limit of the function $p_{2m}\chi^{2n}_0 ((\eh-\frac{1}{2}) (\dd\eh \dd\eh)^n)\in \Om^{2m}_c(G)$ leads to a local expression when $t\to 0$. Recall the explicit formula (\ref{joe3}) 
\beq
\lefteqn{\chi^{2n}_0(x_0\dd x_1\ldots\dd x_{2n})=t^{2n}\sum_{i=0}^{2n} (-)^i \int_{\Delta_{2n+1}}ds\, \tau(e^{-s_0t^2D^2}[D,\rho_*^{i+1}]e^{-s_1t^2D^2}\ldots}\non\\
&&\ldots  [D,\rho_*^{2n}]e^{-s_{2n-i}t^2D^2} \rho_*^0 e^{-s_{2n-i+1}t^2D^2}[D,\rho_*^1]\ldots [D,\rho_*^i]e^{-s_{2n+1}t^2D^2})\ ,\non
\eeq
for any elements $x_1,\ldots, x_{2n}\in T\Ac$, where $\rho_*^i:=\rho_*(x_i)$ is the image of $x_i$ by the homomorphism $\rho_*: T\Ac \to \End_{\Tc\Bc}(\Oman^+\Ec)$. We perform the computation for the last term ($i=n$). The element under the trace $\tau$ lies in the subalgebra of trace-class endomorphisms $\Om_c^+(G;\ell^1)\subset \End_{\Tc\Bc}(\Oman^+\Ec)$ (see section \ref{sbiv}),
$$
e^{-s_0t^2D^2}\rho_*^0 e^{-s_1t^2D^2}[D,\rho_*^1]\ldots [D,\rho_*^{2n}]e^{-s_{2n+1}t^2D^2} \in \Om_c^+(G;\ell^1)\ ,
$$
where the products are of Fedosov type. Explicitly, if $x=a_0da_1\ldots da_{2k}$ then $e^{-sD^2}\rho_*(x)$ is an element of $\Om^{2k}_c(G;\ell^1)$ whose evaluation on a point $(g_0,\ldots,g_{2k})\in G^{2k+1}$ is a trace-class operator on the Hilbert space $\Hc=L^2(E)$:
$$
e^{-sD^2}\rho_*(x)(g_0,\ldots,g_{2k})=e^{-sD^2}a_0(g_0)r(g_0)\ldots a_{2k}(g_{2k})r(g_{2k})\ .
$$
We recall that $a_i(g_i)\in \cinfc(M)$ is viewed as a bounded operator on $\Hc$, and $r:G\to \Uc(\Hc)$ is the representation of $G$ by pullback on the sections of $E$. One has
$$
a_0(g_0)r(g_0)\ldots a_{2k}(g_{2k})r(g_{2k})=a_0(g_0)a_1(g_1)^{g_0}\ldots a_{2k}(g_{2k})^{g_{2k-1}\ldots g_0}r(g_{2k}\ldots g_0)\ ,
$$
which shows the following equality:
$$
e^{-sD^2}\rho_*(x)(g_0,\ldots,g_{2k})=e^{-sD^2}\psi(x)(g_0,\ldots,g_{2k})r(g_{2k}\ldots g_0)\ ,
$$
where $\psi: T\Ac\to \cinfc(M)\rtimes T\Bc$ is the homomorphism (\ref{homo}), given in terms of differential forms by (\ref{medor1}, \ref{medor2}), and $\psi(x)\in \cinfc(M)\hotimes \Om^{2k}\Bc$ is viewed as a function on $G^{2k+1}$ with values in the algebra $\cinfc(M)$ represented in $\Lc(\Hc)$. In the same way, the fact that $D$ commutes with the group representation $r$ implies
$$
e^{-sD^2}[D,\rho_*(x)](g_0,\ldots,g_{2k})=e^{-sD^2}[D,\psi(x)(g_0,\ldots,g_{2k})]r(g_{2k}\ldots g_0)\ ,
$$
and more generally, the evaluation of the $2m$-form
$$
p_{2m}\tau(e^{-s_0t^2D^2}\rho_*^0 e^{-s_1t^2D^2}[D,\rho_*^1]\ldots [D,\rho_*^{2n}]e^{-s_{2n+1}t^2D^2}) \in \Om_c^{2m}(G)
$$
on a point $\gt\in G^{2m}\cup G^{2m+1}$ yields 
$$
p_{2m}\tau \big(e^{-s_0t^2D^2} \psi_0 e^{-s_1t^2D^2}[D,\psi_1]\ldots [D,\psi_{2n}]e^{-s_{2n+1}t^2D^2}r(g) \big)(\gt)
$$
where $\psi_i=\psi(x_i)$, and $g\in G$ is the concatenation product of $\gt$. The same formula holds for all the terms of the sum (\ref{joe3}), modulo cyclic permutations of $\psi_0$, $[D,\psi_1]$, ... , $[D,\psi_{2n}]$. Now, if $f_0,\ldots,f_{2n}\in \cinfc(M)$ are smooth functions with compact support on $M$ and $g\in G$, one can use Getzler's asymptotic symbol calculus \cite{BGV} to show that
\beq
\lefteqn{\lim_{t\to 0}t^{2n} \tau \big( e^{-s_0t^2D^2}f_0 e^{-s_1t^2D^2}[D,f_1]\ldots [D,f_{2n}]e^{-s_{2n+1}t^2D^2} r(g) \big) =}\non\\
&& \qquad\qquad\qquad \sum_{M_g} \frac{(-)^{q/2}}{(2\pi i)^{d/2}}\int_{M_g} \widehat{A}(M_g)\frac{\ch(E/S,g)}{\ch(S_N,g)} \wedge f_0\delta f_1\ldots \delta f_{2n}\ ,\non
\eeq
where the sum runs over the submanifolds $M_g$ of fixed points of dimension $d$ and codimension $q$, and $\delta$ is the de Rham coboundary. A similar formula holds with the cyclic permutations of $f_0$, $[D,f_1]$, ... , $[D,f_{2n}]$. Finally, taking into account the integration over the simplex $\Delta_{2n+1}$ which yields a factor of $1/(2n+1)!$, we are able to evaluate the function $p_{2m}\chi^{2n}_0(x_0\dd x_1\ldots\dd x_{2n}) \in \Om_c^{2m}(G)$ on a point $\gt$ at the limit $t\to 0$:
\beq
\lefteqn{\lim_{t\to 0}p_{2m}\chi^{2n}_0(x_0\dd x_1\ldots\dd x_{2n})(\gt)=\sum_{M_g} \frac{(-)^{q/2}}{(2\pi i)^{d/2}(2n+1)!}\times}\non\\
&&\int_{M_g}\widehat{A}(M_g)\frac{\ch(E/S,g)}{\ch(S_N,g)}\wedge\sum_{i=0}^{2n}(-)^i p_{2m} \big(\delta\psi_{i+1}\ldots \delta\psi_{2n} \psi_0\delta\psi_1\ldots \delta\psi_i\big)(\gt)\non
\eeq
where the element $p_{2m}(\delta\psi_{i+1}\ldots \delta\psi_{2n} \psi_0\delta\psi_1\ldots \delta\psi_i)\in \Om_c^{2n}(M)\hotimes \Om^{2m}\Bc$ is viewed as a function on $G^{2m}\cup G^{2m+1}$ with values in $\Om_c^{2n}(M)$, and $g\in G$ is the concatenation product of $\gt$. It is therefore just a matter of simple algebraic manipulations, using the idempotent properties of $\eh\in\Tc\Ac$, to show that 
$$
\lim_{t\to 0}\ch_{2m}(tD)(\gt)=\sum_{M_g}\frac{(-)^{q/2}}{(2\pi i)^{d/2}}\int_{M_g}\widehat{A}(M_g)\frac{\ch(E/S,g)}{\ch(S_N,g)}\, \ch_{2m}(e)(\gt)
$$
as wanted. The odd case is obtained in the same way. \cqfd\\
\begin{remark}\textup{It must be noted that since the action of $G$ is proper, an element $g\in G$ can have fixed points on $M$ only if it is contained in a compact subgroup of $G$. This implies that the function $\lim_{t\to 0}\ch_{n}(tD)$ is non-zero only when the concatenation $g$ of the point $\gt$ belongs to a compact subgroup.}
\end{remark}
\begin{example}\textup{Suppose $G$ is compact and $\Bc=L^1(G,dg)$ is the algebra of integrable functions with respect to the Haar measure $dg$. For a compact group the only interesting cyclic cocycles over $L^1(G)$ are just given by bounded traces $\varphi:L^1(G)\to\cc$, or equivalently by cyclic cocycles of degree zero. One can write
\be
\varphi(b)=\int_Gdg\, \phi(g)\, b(g)\quad \forall b\in L^1(G)\ ,
\ee
with $\phi\in L^{\infty}(G,dg)$. $\varphi$ is a trace iff the function $\phi$ is constant along the adjoint orbits in $G$. Since $\varphi$ is a cyclic cocycle of degree zero, its pairing with the Chern character of $\mu(D)$ only involves the degree zero component $\ch_{0}(tD)\in \Bc$. Hence all the information we need is concentrated in degree zero. Next, we note that the manifold $M$ is necessarily compact, and the idempotent representing the canonical $K$-theory class $[e]\in K_0(\Ac)$ may be built out of a constant cut-off function $c=1$ over $M$, provided the Haar measure $dg$ is normalized so that $\int_Gdg=1$. Therefore 
\be
e(g,x)=1\quad \forall g\in G\ ,\ x\in M\ .
\ee
This particular choice of cut-off function simplifies drastically the computation of the Chern character $\ch_*(e)$, whose degree zero component is a constant function on $G$ with values in $\Om_c(M)$:
\be
\ch_{0}(e)(g)=1\quad \forall g\in G\ .
\ee
The localization formula of corollary \ref{cloc} then gives
\be
\lim_{t\to 0}\ch_{0}(tD)(g)=\sum_{M_g}\frac{(-)^{q/2}}{(2\pi i)^{d/2}}\int_{M_g}  \widehat{A}(M_g)\frac{\ch(E/S,g)}{\ch(S_N,g)}
\ee
for any $g\in G$. One recovers the Atiyah-Segal-Singer equivariant index theorem for compact groups \cite{AS,BGV}.}
\end{example}
\begin{example}\textup{The Connes-Moscovici index theorem for coverings \cite{CM90} is also a consequence of \ref{cloc}. Indeed, let $G$ be a countable discrete group acting freely and properly on $M$, with $X=G\backslash M$ a compact manifold of dimension $d$. We may view $M$ as a $G$-principal bundle over $X$, with $G$ acting by deck transformations. Any complex vector bundle $E$ over $X$ can be lifted to a $G$-equivariant vector bundle $\widetilde{E}$ over $M$, and if $D:\cinf(E)\to\cinf(E)$ is a generalized Dirac operator, its lift $\widetilde{D}$ is a $G$-invariant operator over $M$. \\
Since $G$ is discrete, the convolution algebra $C_c(G)$ is isomorphic to the group ring $\cc G$. Let $v\in Z^n(G)$ be a group $n$-cocycle, 
\be
v(g_2,\ldots,g_{n+1})+\sum_{i=1}^n(-)^i v(g_1,\ldots,g_ig_{i+1},\ldots,g_{n+1})+\ldots +(-)^{n+1}v(g_1,\ldots,g_n)=0\ ,
\ee
assumed to be normalized in the sense that $v(g_1,\ldots, g_n)=0$ whenever one of the $g_i$'s is equal to 1 or $g_1\ldots g_n=1$. Following \cite{CM90}, we associate to $v$ a cyclic $n$-cocycle $\varphi_v:\Om^n C_c(G)\to\cc$ over the convolution algebra, defined by the formula
\be
\varphi_v(b_0db_1\ldots db_n)=\sum_{g_1,\ldots, g_n\in G} v(g_n,\ldots, g_1)\, b_0((g_n\ldots g_1)^{-1})b_1(g_1)\ldots b_n(g_n)\ ,
\ee
and $\varphi_v(db_1\ldots db_n)=0$ for any $b_i\in C_c(G)$. Now assume that $v$ has polynomial growth with respect to a right-invariant distance function on $G$, and choose an admissible completion $\Bc$ of $C_c(G)$ as in example \ref{edist}. If the parameter $\al$ defining the norm on $\Bc$ is large enough, then $\varphi_v$ extends to a cyclic $n$-cocycle over $\Bc$. We obtain a higher index of $D$ by pairing the Chern character of $\mu(\widetilde{D})\in K_*(\Bc)$ with $\varphi_v$. The localization formula \ref{cloc} gives
\be
\langle \varphi_v,\mu(\widetilde{D})\rangle=\frac{1}{(2\pi i)^{d/2}}\int_X \widehat{A}(X)\ch(E/S)\, f^*(v)\ ,
\ee 
where $f^*(v)$ is the pullback of the cohomology class $v\in H^*(G)\cong H^*(BG)$ with respect to the classifying map $f:X\to BG$ associated to the $G$-principal bundle $M\stackrel{G}{\longrightarrow} X$. In this case, the $G$-invariant $\widehat{A}$-genus and Chern character $\ch(E/S)$ are considered as closed differential forms over $X$.}
\end{example}



\end{document}